%% file: tameflow.tex
\documentclass[reqno, 11pt]{amsart}

\newtheorem{theorem}{Theorem}[section]
\newtheorem{lemma}[theorem]{Lemma}
\newtheorem{proposition}[theorem]{Proposition}
\newtheorem{corollary}[theorem]{Corollary}

\theoremstyle{definition}
\newtheorem{definition}[theorem]{Definition}
\newtheorem{ex}[theorem]{Example}

\theoremstyle{remark}
\newtheorem{remark}[theorem]{Remark}

\numberwithin{equation}{section}

\usepackage{bbm}
\usepackage{euscript}
\usepackage{pb-diagram}
\usepackage{lamsarrow}
\usepackage{pb-lams}
\usepackage{amsmath}
\usepackage{amsthm}

\usepackage{epsfig}
\usepackage{amssymb}
\usepackage{pifont}
\usepackage{amsbsy}

\newskip\aline \newskip\halfaline
\def\skipaline{\vskip\aline}

\def\qedbox{$\rlap{$\sqcap$}\sqcup$}
\def\qed{\nobreak\hfill\penalty250 \hbox{}\nobreak\hfill\qedbox\skipaline}

\def\proofend{\eqno{\mbox{\qedbox}}}

\newcommand{\one}{\mathbbm{1}}

\newcommand\bD{{\mathbb D}}

\newcommand{\bM}{{{\mathbb M}}}

\newcommand\bR{{\mathbb R}}

\newcommand\bZ{{\mathbb Z}}

\newcommand\bRP{{\mathbb R}{\mathbb P}}

\DeclareMathOperator{\tr}{{\rm tr}}
\DeclareMathOperator{\supp}{{\rm supp}}

\DeclareMathOperator{\dist}{dist}
\DeclareMathOperator{\Int}{\boldsymbol{Int}}

\DeclareMathOperator{\Gr}{\mathbf{Gr}}
\DeclareMathOperator{\diag}{Diag} 
\DeclareMathOperator{\Vect}{Vect} \DeclareMathOperator{\End}{End}

\newcommand{\be}{\boldsymbol{e}}

\newcommand{\bh}{\boldsymbol{h}}
\newcommand{\ii}{\boldsymbol{i}}

\newcommand{\bp}{\boldsymbol{p}}
\newcommand{\br}{\boldsymbol{r}}

\newcommand{\bv}{\boldsymbol{v}}

\newcommand{\bsI}{\boldsymbol{I}}

\newcommand{\si}{{\sigma}}

\newcommand{\de}{{\delta}}

\newcommand{\ve}{{\varepsilon}}
\newcommand{\eps}{{\epsilon}}
\newcommand{\vfi}{{\varphi}}

\newcommand{\eA}{\EuScript{A}}
\newcommand{\eB}{\EuScript{B}}
\newcommand{\eC}{\EuScript{C}}
\newcommand{\eD}{\EuScript{D}}

\newcommand{\eF}{\EuScript{F}}
\newcommand{\eG}{\EuScript{G}}
\newcommand{\eH}{\EuScript H}

\newcommand{\eK}{\EuScript{K}}
\newcommand{\eL}{\EuScript{L}}
\newcommand{\eM}{\EuScript{M}}
\newcommand{\eN}{\EuScript{N}}
\newcommand{\eO}{\EuScript{O}}
\newcommand{\eP}{\EuScript{P}}

\newcommand{\eR}{\EuScript{R}}
\newcommand{\eS}{\EuScript{S}}
\newcommand{\eT}{\EuScript{T}}

\newcommand{\eV}{\EuScript{V}}
\newcommand{\eW}{\EuScript{W}}

\newcommand{\eY}{\EuScript{Y}}

\newcommand{\ra}{\rightarrow}
\newcommand{\hra}{\hookrightarrow}
\newcommand{\Lra}{{\longrightarrow}}

\newcommand{\lan}{\langle}
\newcommand{\ran}{\rangle}

\def\inpr{\mathbin{\hbox to 6pt{\vrule height0.4pt width5pt depth0pt \kern-.4pt \vrule height6pt width0.4pt depth0pt\hss}}}

\newcommand{\pa}{\partial}

\newcommand{\rsq}{\rightsquigarrow}
\newcommand{\rsqf}{\stackrel{f}{\rightsquigarrow}}
\newcommand{\precf}{\preccurlyeq_f}
\newcommand{\ori}{\boldsymbol{or}}

\DeclareMathOperator{\Aff}{\mathbf{Aff}}
\DeclareMathOperator{\St}{\mathbf{St}}

\DeclareMathOperator{\Fr}{\mathbf{Fr}}
\DeclareMathOperator{\clos}{\mathbf{clos}}
\DeclareMathOperator{\ev}{\mathbf{ev}}

\DeclareMathOperator{\cl}{\boldsymbol{cl}}

\DeclareMathOperator{\Cr}{\mathbf{Cr}}
\DeclareMathOperator{\Sym}{Sym}

\DeclareMathOperator{\pfc}{\widehat{\EuScript{S}}_{an}}
\DeclareMathOperator{\conv}{\mathbf{conv}}
\DeclareMathOperator{\Lk}{\mathbf{Lk}}
\DeclareMathOperator{\bd}{\boldsymbol{bd}}
\DeclareMathOperator{\Graff}{\mathbf{Graff}}
\DeclareMathOperator{\Cyl}{\mathbf{Cyl}}
\newcommand{\Lla}{\longleftarrow}
\newcommand{\BXi}{\boldsymbol{\Xi}}
\DeclareMathOperator{\Flg}{\mathbf{Fl_{grad}}}
\begin{document}

\title{Tame Flows}
\date{Started: March 6, 2006. Completed: February 10, 2007. First revision. April 2007.}

\author{Liviu I. Nicolaescu}

\address{Department of Mathematics, University of Notre Dame, Notre Dame, IN 46556-4618.}
\email{nicolaescu.1@nd.edu} \urladdr{http://www.nd.edu/~lnicolae/}

\begin{abstract}
The tame flows are  ``nice'' flows on ``nice'' spaces. The  nice (tame)  sets are the pfaffian sets introduced by Khovanski,  and a  flow  $\Phi: \mathbb{R}\times X\rightarrow X$ on  pfaffian set $X$  is tame if  the graph of $\Phi$ is a pfaffian subset of $\mathbb{R}\times X\times X$.  Any compact  tame set  admits plenty tame  flows. We prove that the  flow  determined by the gradient of a generic  real analytic  function with respect to a generic   real analytic metric is tame.  The typical tame gradient  flow  satisfies the Morse-Smale condition, and we prove that  in the tame  context  the Morse-Smale condition is equivalent to the fact that the stratification  by unstable manifolds is Verdier and  Whitney regular. We explain how to compute the   Conley indices of  isolated  stationary points of tame flows in terms of their unstable varieties, and then use  this technology to   produce a Morse theory on  posets  generalizing R. Forman's discrete Morse theory. Finally, we use the  Harvey-Lawson finite volume flow  technique to produce a homotopy   between the DeRham complex of a smooth manifold and the  simplicial chain complex associated to a triangulation.

\end{abstract}

\maketitle

\tableofcontents


\section*{Introduction}
\addtocontents{section}{Introduction}
\input{tameflowintr.tex}

\section{Tame spaces}
\label{s: 1}
\setcounter{equation}{0}
\input{tameflow1.tex}

\section{Basic properties and examples of tame flows}
\label{s: 2}
\setcounter{equation}{0}
\input{tameflow2.tex}

\section{Some global properties  of tame flows}
\label{s: 3}
\setcounter{equation}{0}
\input{tameflow3.tex}

\section{Tame Morse flows}
\label{s: 4}
\setcounter{equation}{0}
\input{tameflow4.tex}

\section{Tame Morse-Smale flows}
\label{s: 5}
\setcounter{equation}{0}
\input{tameflow5.tex}

\section{The gap between two vector subspaces}
\label{s: 6}
\setcounter{equation}{0}
\input{tameflow6.tex}

\section{The Whitney and Verdier regularity conditions}
\label{s: 7}
\setcounter{equation}{0}
\input{tameflow7.tex}

\section{Morse-Smale $\Longleftrightarrow$ Morse-Verdier}
\label{s: 8}
\setcounter{equation}{0}
\input{tameflow8.tex}

\section{The Conley index}
\label{s: 9}
\setcounter{equation}{0}
\input{tameflow9.tex}

\section{Flips/flops and gradient like tame flows}
\label{s: 10}
\setcounter{equation}{0}
\input{tameflow10.tex}

\section{Simplicial flows}
\label{s: 11}
\setcounter{equation}{0}
\input{tameflow11.tex}

\section{Tame currents}
\label{s: 12}
\setcounter{equation}{0}
\input{tameflow12.tex}

\newpage

\appendix

\section{An ``elementary'' proof of the generalized Stokes formula}
\label{s: a}
\setcounter{equation}{0}
\input{tameflow-a.tex}

\section{On the topology of tame sets}
\label{s: b}
\setcounter{equation}{0}
\input{tameflow-b.tex}

\newpage

\input{tameflow-bib.tex}
\end{document}

%% file: tameflowintr.tex
Loosely speaking, the tame   sets   (respectively tame flows)  are
sets  (respectively  continuous flows)   which display very few
pathologies.   Technically speaking, they are sets or  flows
definable within a tame structure.

The  subject of $o$-minimal or tame geometry is  not as popular as
it ought to be  in geometric circles, although this situation is
beginning to change.  The tame geometry is a  vast generalization
of  the  more classical subject of real algebraic geometry.   One
such extension  of real algebraic geometry   was conceived and
investigated  by A. Khovanski  in \cite{Kho}, and our  tame sets
are the precisely Khovanski's pfaffian sets. In particular, all
the tame sets will be subsets of   Euclidean spaces.

 If we think of a flow as  generated by  a system of ordinary differential equations then,   roughly speaking, the tame flows correspond   to first order ordinary differential equations which we can solve explicitly by quadratures, with one  important caveat:   the resulting    final   description of the solutions should not involve trigonometric functions  because  tame flows  do not  have periodic orbits. For example, a linear system of ordinary differential equations determines a tame flow if and only if the defining matrix has only real eigenvalues.

Given that the tame  sets display  very few pathologies, they
form a much more restrictive class of subsets of Euclidean spaces,
and in particular, one might expect that   the tame flows    are
not as plentiful. In the present paper we set up to convince the
reader that  there is a rather large  supply of such flows,  and
that they are worth investigating due to their  rich structure.

The paper is  structured around three major themes: examples of
tame flows, properties of tame flows, and applications of tame
flows.

To  produce  examples of tame flows  we describe   several general
classes of tame flows, and  several general surgery like
operations on tame  flows which  generate  new flows out of some
given ones.  These operations  have a simplicial flavor: we can
cone a flow,  we can  produce  join of two flows, or we can glue
two flows along a common, closed invariant subset.

The  simplest example of  tame flow is the trivial flow on a set
consisting of single point. An iterated application of the cone
operations produces  canonical tame flows on any  affine
$m$-simplex, and then by gluing, on any triangulated  tame set.
Since any  tame set can be triangulated, we conclude that  there
exist many tame flows on any tame set.

Another class of tame flows,  which cannot be obtained by the
cone  operation,  consists   of  the gradient flows of  ``most''
real analytic functions  on a real analytic manifold equipped with
a real analytic metric.

More precisely, we prove  that, for any real analytic  $f$
function on a real analytic  manifold $M$, there exists   a dense
set  of real analytic metrics $g$ with the property that the  flow
generated by $\nabla^g f$ is  tame. This is a rather nontrivial
result,  ultimately based on the  Poincar\'e-Siegel theorem
concerning the canonical form of a vector field in a neighborhood
of a stationary point. The Poincar\'e-Siegel  theorem plays the
role of the more elementary Morse lemma.

The   usual techniques pioneered by Smale   show that a   tame
gradient flow can be slightly modified to a gradient like tame
flow satisfying the  Morse-Smale regularity conditions.

We  investigate the  stratification  of a manifold given by the
unstable manifolds of the downward gradient  flow of some real
analytic function.  We observe that if this   stratification
satisfies the Verdier  regularity condition, then  the flow
satisfies the Morse-Smale conditions.    We    succeeded in
proving an (almost) converse to this statement. More precisely, we
show   that if the  tame gradient flow  associated  to a real
analytic function $f$ and metric $g$   satisfies the Morse-Smale
condition, and if for  every unstable  critical point $x$ of  $f$,
the spectrum $\Sigma_x$ of the Hessian of $f$  at $x$ satisfies
the \emph{clustering  condition}
\[
\max \Sigma_x < \dist(\Sigma_x^+,0)
+\dist(\Sigma_x^-,0),\;\;\mbox{where}\,\;\Sigma^\pm_x:=\bigl\{
\lambda\in \Sigma_x;\;\;\pm \lambda >0\,\bigr\},
\]
then the stratification   by unstable manifolds  satisfies the
Verdier regularity condition.    Again, the Poincar\'e-Siegel
theorem shows that   set of tame gradient flows satisfying the
spectral clustering condition above is nonempty and  ``open''.

In the tame world, the Verdier condition implies the Whitney
condition. We deduce that the  unstable manifolds  of a tame Morse-Smale
flows satisfying the  spectral clustering condition form a Whitney
stratification.    The  results of F. Laundebach \cite{Lauden}  on
the local conical structure of the stratification by the unstable
manifolds  follow from the general results on the local structure
of a Whitney stratified space.

We  also investigate Morse like  tame flows on  singular spaces,
i.e. tame Morse flows which admit a Lyapunov function.  We explain
how to compute the   (homotopic) Conley index of an isolated
stationary point in terms of the unstable variety of that point.
We achieve this by  proving a singular counterpart of the
classical  result in Morse theory:    crossing a  critical level
of a Morse function  corresponds homotopically to attaching a cell
of a certain dimension.    Since we are working on singular spaces
the change in the homotopy type is a bit more complicated, but
again, crossing a critical  level has a similar homotopic flavor.
The sublevel sets of the Lyapunov function  change by a cone
attachment. The cone has a very precise  dynamical description,
namely it is the cone  spanned by the trajectories of the flow
``exiting''  the stationary point.

   We use the Conley index computation in the study of certain Morse like flows on simplicial  complexes associated to finite posets.   When we specialize to the case of posets of faces of a regular CW decomposition    of a space we obtain, as a very special case , R. Forman's discrete Morse theory, \cite{For}.

We also blend the  tameness with the finite volume  techniques of
Harvey-Lawson to prove that the DeRham complex  of a compact,
orientable smooth manifold  is naturally homotopic to the
simplicial chain complex  (with real coefficients) of a
triangulation of the manifold.   No tameness assumption on $M$ is
needed, by the  tameness sneaks in through the back door, in  the
proof.

Here   is briefly the  organization of  the paper.   Section \ref{s: 1} is
a  crash course in tame geometry  where we define    precisely
the meaning of the attribute ``tame'' and  list without proofs  a
few   geometric consequences of tameness  used throughout the
paper.  In Section \ref{s: 2},\ref{s: 3}   we describe  a large list of examples of
tame flows,    and  prove several elementary properties  of an
arbitrary tame flow.  In particular, in these sections we describe
in detail   some canonical tame flows  on  affine simplices
(Example \ref{ex: simplicial}), and on Grassmannians (Example
\ref{ex: grass}) which will play an important role in the paper.

Sections \ref{s: 4}-\ref{s: 8} are devoted to   gradient flows  determined by   real
analytic functions  on real analytic manifolds equipped with real
analytic metrics.  We  prove that   ``most'' of these flows  are
tame (Theorem \ref{th: PS}), satisfy the Morse-Smale condition
(Theorem \ref{th: MS}), and  moreover,  that the Morse-Smale is
equivalent with the fact  that  the stratification   by unstable
manifolds is Verdier and Whitney regular (Theorem \ref{th: vis},
\ref{th: siv}).

In Section \ref{s: 9} we   describe   how to  compute  the Conley index  of
an isolated stationary point of a tame flow admitting Lyapunov
functions in terms of the unstable variety of that point (Theorem
\ref{th: conley}). In Section \ref{s: 10} we produce an almost complete topological classification  of  gradient like tame flows  with finitely many stationary points on compact tame spaces. In Section \ref{s: 11}  we use  the   Conley index
computations to investigate the homotopy type of posets by using
certain tame flows  associated to  certain discrete Morse like
functions on  posets (Theorem \ref{th: simp-flow}, \ref{th: cminus}). In the last section  we  explain  how to use
the Harvey-Lawson techniques to   produce results about the homotopy
type of the DeRham complex (Theorem \ref{th: derham-simpl}).

\bigskip

\noindent {\bf Acknowledgment.}  \hyphenation{Star-chen-ko} I
learned the basics of tame geometry from my colleague Sergei
Starchenko. I   am grateful to  him for  his  patience, generosity
and expertise   in answering  my  zillion   questions      about
this subject.

I also want to thank Adam Dzedzej for     communicating to me several corrections to the original  version of this work.

%% file: tameflow1.tex
Since the subject of tame  geometry is  not     very familiar to many  geometers we  devote this section to a   brief  introduction  to this topic.   Unavoidably, we will  have to omit many interesting details and  contributions, but  we refer to \cite{Co, Dr1, DrMi2} for  more systematic presentations. For every set $X$ we will denote by $\eP(X)$ the   collection of all subsets of $X$

An  \emph{$\bR$-structure}\footnote{This is a highly condensed and   special version of  the traditional definition of structure.     The model theoretic definition  allows for  ordered fields, other than $\bR$, such as extensions of $\bR$ by ``infinitesimals''. This can come in handy even if  one is interested  only in the  field $\bR$.}  is a collection $\eS=\bigl\{\, \eS^n\,\bigr\}_{n\geq 1}$,  $\eS^n\subset \eP(\bR^n)$, with the following properties.

\begin{description}

\item[${\bf E}_1$]   $\eS^n$ contains all the real algebraic subvarieties of $\bR^n$, i.e., the zero sets of  finite collections of polynomial in $n$ real variables.

\item[${\bf E}_2$]  For every   linear map $L:\bR^n\ra \bR$, the half-plane $\{\vec{x}\in \bR^n;\;\;L(x)\geq 0\}$  belongs to $\eS^n$.

\item[${\bf P}_1$] For every $n\geq 1$, the family $\eS^n$ is closed under  boolean operations, $\cup$, $\cap$ and complement.

\item[${\bf P}_2$]  If $A\in \eS^m$, and $B\in \eS^n$, then $A\times B\in \eS^{m+n}$.

\item[${\bf P}_3$]  If $A\in \eS^m$, and $T:\bR^m\ra \bR^n$ is an affine map, then $T(A)\in \eS^n$.

\end{description}

\begin{ex}[Semialgebraic sets]    Denote by $\eS_{alg}$ the collection of real semialgebraic sets.  Thus,  $A\in \eS^n_{alg}$ if and only if  $A$  is a finite  union of sets,  each of which is described by finitely many polynomial equalities and inequalities. The celebrated Tarski-Seidenberg theorem states that $\eS_{alg}$ is a structure.

For example,  the set
\[
A =\bigl\{ (x,y,z)\in \bR^3;   x^3+y^4+z^5 \geq 6 xyz\;\;\mbox{or}\;\;x^2+z^2\geq 1\,\bigr\}
\]
is  semialgebraic, and Tarski-Seidenberg theorem  implies that its projection on the  $(x,y)$-plane
\[
\bigl\{ (x,y)\in \bR^2;\;\;\exists z\in \bR: \;\; x^3+y^4+z^5 \geq 6 xyz,\;\;\mbox{or}\;\;x^2+z^2\geq 1\,\bigr\}
\]
is also  semialgebraic.\qed
\end{ex}

 Given a structure $\eS$, then an $\eS$-\emph{definable} set is a set that belongs to one of the  $\eS^n$-s. If $A, B$ are $\eS$-definable, then a function $f: A\ra B$ is called $\eS$-\emph{definable} if its graph
 \[
 \Gamma_f :=\bigl\{ (a,b)\in A\times B;\;\;b=f(a)\,\bigr\}
 \]
 is $\eS$-definable. The reason   these sets are called definable has to do with mathematical logic.

 A  \emph{formula}\footnote{We are deliberately vague on the meaning of formula.} is a property  defining a certain set. For example, the two different looking formulas
 \[
 \bigl\{ x\in \bR;\;\;x\geq 0\},\;\;\bigl\{ x\in \bR;\;\;\exists y\in\bR:\;\;x=y^2\},
 \]
 describe the  same set, $[0,\infty)$.

  Given a collection of formulas, we can obtain  new formulas, using the logical operations $\wedge,\vee, \neg$, and  quantifiers $\exists$,  $\forall$. If we start   with a collection of formulas, each describing an $\eS$-definable set, then any formula obtained from them by applying  the above logical  transformations will    describe a definable set.

 To see this, observe that the operators $\wedge,\vee,\neg$ correspond to the boolean operations, $\cap,\cup$, and taking the complement.  The  existential quantifier corresponds to taking a projection. For example,    suppose we are given a formula $\phi(a,b)$,  $(a,b)\in A\times B$, $A,B$ definable, describing a   definable set $C\subset A\times B$.  Then the formula
 \[
\bigl\{\, a\in A;\;\exists b\in B:\;\;\phi(a,b)\,\bigr\}
 \]
 describes the image of the subset $C\subset A\times B$ via the  canonical projection  $A\times B\ra A$. If $A\subset \bR^m$, $B\subset \bR^n$, then the projection $A\times B\ra A$ is the restriction to $A\times B$ of the linear projection $\bR^m\times \bR^n\ra \bR^m$ and  ${\bf P}_3$ implies that the image of $C$ is also definable. Observe that the universal quantifier    can be replaced with the operator $\neg\exists\neg$.

\begin{ex} (a) The composition of two definable functions $A\stackrel{f}{\ra} B\stackrel{g}{\ra} C$ is a definable function because
\[
\Gamma_{g\circ f}=\bigl\{ \,(a,c)\in A\times C; \exists  b\in B: (a,b)\in \Gamma_f,\;\;(b,c)\in \Gamma_g\,\bigr\}.
\]
Note that any polynomial  with real coefficients is a definable function.

\noindent (b) The image  and the preimage of a definable set via a definable  function is a definable set.    Using ${\bf E}_2$ we deduce that any semialgebraic set $\eS$ is definable. In particular, the   Euclidean norm
\[
|\bullet|: \bR^n\ra \bR,\;\; |(x_1,\dotsc,x_n)|=\Bigl(\sum_{i=1}^n x_i^2\Bigr)^{1/2}
\]
is $\eS$-definable.

\noindent (c) Suppose $A\subset \bR^n$ is definable. Then its closure $\cl(A)$ is described by the formula
\[
\bigl\{ x\in \bR^n;\;\;\forall \ve >0,\;\;\exists a\in A:\;\;|x-a|<\ve\,\bigr\},
\]
and we deduce that $\cl(A)$ is also definable. Let us examine the correspondence between the operations on formulas and operations on sets on this example.

We rewrite this formula as
\[
\forall \ve \Bigl(\, (\ve>0) \Rightarrow \exists  a (a\in A) \wedge (x\in\bR^n)\wedge  (|x-a|<\ve)\,\Bigr).
\]
In the above formula we see one free variable $x$, and the set  described by this formula consists of those $x$ for which that formula is a true statement.

The above formula is made of  the ``atomic'' formul{\ae},
\[
(a\in A),\;\;(x\in \bR^n),\;\; (|x-a|<\ve),\;\;(\ve>0),
\]
which all describe definable sets.  The logical connector $\Rightarrow$ can be replaced by $\vee\neg$. Finally, we can replace the  universal quantifier to rewrite the formula   as a  transform of atomic formulas via the basic logical operations.
\[
\neg\Bigl\{\exists \ve \neg\Bigl(\, (\ve>0) \Rightarrow \exists  a (a\in A) \wedge (x\in\bR^n)\wedge  (|x-a|<\ve)\,\Bigr)\Bigr\}.\proofend
\]
\end{ex}

Given an $\bR$-structure $\eS$, and a collection $\eA=(\eA_n)_{n\geq 1}$, $\eA_n\subset\eP(\bR^n)$, we can form a new structure  $\eS(\eA)$, which is the smallest structure containing   $\eS$ and the sets in $\eA_n$. We say that $\eS(\eA)$ is obtained  from  $\eS$ by \emph{adjoining the collection $\eA$}.

\begin{definition} An $\bR$-structure is called \emph{$o$-minimal} (order minimal) or \emph{tame} if  it satisfies the property

\begin{description}
\item[{\bf T}]    Any set $A\in  \eS^1$ is a \emph{finite} union of open intervals $(a,b)$, $-\infty \leq a <b\leq \infty$, and singletons $\{r\}$. \qed
\end{description}

\end{definition}

\begin{ex} (a) The collection  $\eS_{alg}$ of real semialgebraic sets  is a tame structure.

\noindent (b)(Gabrielov-Hironaka-Hardt)   A \emph{restricted} real analytic function is a  function $f:\bR^n\ra \bR$ with the property that there exists a real analytic function $\tilde{f}$ defined in an open  neighborhood $U$ of the cube $C_n:=[-1,1]^n$ such that
\[
f(x)=\begin{cases}
\tilde{f}(x) & x\in C_n\\
0 & x\in \bR^n\setminus C_n.
\end{cases}
\]
we denote by $\eS_{an}$ the structure obtained from $\eS_{alg}$ by adjoining the  graphs of all the restricted real analytic functions. Then $\eS_{an}$ is a tame structure, and the $\eS_{an}$-definable sets are called \emph{globally subanalytic sets}.

\noindent(c)(Wilkie, van den Dries, Macintyre, Marker)   The
structure  obtained  by adjoining to  $\eS_{an}$ the graph of the
exponential function $\bR\ra \bR$, $t\mapsto e^t$,  is an
tame structure.

\noindent(d)(Khovanski-Speissegger)   There   exists  a tame structure  $\eS'$ with the following properties
\begin{enumerate}
\item[($d_1$)] $\eS_{an} \subset \eS'$

\item[($d_2$)]  If $U\subset \bR^n$ is open, connected  and $\eS'$-definable, $F_1,\dotsc, F_n: U\times \bR\ra \bR$ are $\eS'$-definable and $C^1$, and  $f: U\ra \bR$ is a $C^1$ function satisfying
\begin{equation}
\frac{\pa f}{\pa x_i} = F_i(x, f(x)),\;\;\forall x\in \bR,,\;\;i=1,\dotsc, n,
\label{eq: pfaff}
\end{equation}
then  $f$ is $\eS'$-definable.
\end{enumerate}
The smallest structure  satisfying the above two properties, is called the \emph{pfaffian closure}\footnote{Our definition  of pfaffian closure   is more restrictive than the original one in \cite{Kho, Sp1}, but it suffices for the geometrical applications we have in mind.}
  of $\eS_{an}$, and we will denote it by $\pfc$.

Observe that if $f: (a,b)\ra \bR$ is $C^1$,  $\pfc$-definable,  and $x_0\in (a,b)$ then the antiderivative $F:(a, b)\ra \bR$
  \[
  F(x)=\int_{x_0}^x f(t)dt,\;\; x\in (a,b),
  \]
  is also $\pfc$-definable. \qed

\end{ex}

 The  definable sets  and function of a tame structure have  rather remarkable \emph{tame} behavior which prohibits  many pathologies.  It is perhaps   instructive to give an example of function which is not definable in any tame structure. For example, the function $x\mapsto \sin x$ is not definable in a tame structure because the intersection of its graph with the horizontal axis is the  countable set $\pi\bZ$  which  violates  the $o$-minimality condition ${\bf O}$.

 We will list below some of the nice properties of the sets and function definable  in a tame structure  $\eS$. Their proofs can be found in \cite{Dr1}.

\smallskip

 \noindent $\bullet$ (\emph{Piecewise  smoothness of  one variable tame functions.}) If $f:(0,1)\ra \bR$ is an $\eS$-definable function, and $p$ is a positive integer, then there exists
 \[
0=a_0< a_1<a_2<\cdots <a_n=1
\]
such that the restriction of  $f$ to each subinterval  $(a_{i-1},a_i)$ is $C^p$ and   monotone. Moreover $f$ admits right and left limits at any $t\in[0,1]$.

\noindent $\bullet$ (\emph{Closed graph theorem.})  Suppose $X$ is a tame
set and $f: X\ra \bR^n$ is a  tame bounded function.  Then $f$ is
continuous if and only if  its graph is closed in $X\times \bR^n$.

\noindent $\bullet$ (\emph{Curve selection.}) If $A$ is an $\eS$-definable set, and $x\in\cl(A)\setminus A$, then there exists an $\eS$ definable  continuous map
\[
\gamma:(0,1)\ra A
\]
such that $x=\lim_{t\ra 0} \gamma(t)$.

\noindent $\bullet$ Any definable set has finitely many  connected components, and each of them  is definable.

\noindent $\bullet$   Suppose $A$ is an $\eS$-definable  set, $p$ is a positive integer, and $f: A\ra \bR$ is a definable function. Then $A$ can be partitioned into finitely many  $\eS$ definable sets $S_1,\dotsc, S_k$,     such that each  $S_i$ is a $C^p$-manifold, and each of the restrictions $f|_{S_i}$ is a $C^p$-function.

\noindent $\bullet$ (\emph{Triangulability})  For every   compact definable set $A$, and any finite collection of definable  subsets $\{S_1,\dotsc, S_k\}$, there exists  a compact simplicial complex $K$, and a  definable homeomorphism
\[
\Phi: K\ra A
\]
such that  all the sets $\Phi^{-1}(S_i)$ are unions of  relative interiors of faces of $K$.

\noindent $\bullet$ (\emph{Definable selection}.)  Suppose $A,\Lambda$ are $\eS$-definable. Then a \emph{definable} family of subsets of $A$ parameterized by $\Lambda$ is a definable subset
\[
S\subset A\times \Lambda.
\]
We set
\[
S_\lambda:= \bigl\{ a\in A;\;\; (a,\lambda)\in S\,\bigr\},
\]
and we denote by $\Lambda_S$ the projection of $S$ on $\Lambda$. Then there exists a definable function $s:\Lambda_S\ra A$ such that
\[
s(\lambda)\in S_\lambda,\;\;\forall \lambda\in \Lambda_S.
\]

\noindent $\bullet$ (\emph{Dimension.})  The  dimension of an $\eS$-definable  set $A\subset \bR^n$ is the supremum over all the nonnegative integers $d$ such that there exists a $C^1$  submanifold of $\bR^n$ of dimension $d$ contained in $A$.  Then $\dim A <\infty$, and
\[
\dim (\cl(A)\setminus A) <\dim A.
\]
Moreover, if $(S_\lambda)_{\lambda\in\Lambda}$ is a definable family of definable sets then the function
\[
\Lambda\ni\lambda \mapsto\dim S_\lambda
\]
is definable.

\noindent $\bullet$ (\emph{Definable triviality of tame maps.}) We say that a tame map $\Phi: X\ra S$ is \emph{definably trivial}   if there exists a    definable set  $F$, and a definable homeomorphism $\tau: X\ra F\times S$ such that the diagram below is commutative
\[
\begin{diagram}
\node{X}\arrow{se,b}{\Phi}\arrow[2]{e,t}{\tau}\node[2]{S\times F}\arrow{sw,b}{\pi_S}\\
\node{}\node{S}\node{}
\end{diagram}.
\]
If $\Psi: X\ra Y$ is a definable map, and $p$ is a positive integer, then there exists  a partition of $Y$ into definable $C^p$-manifolds  $Y_1,\dotsc, Y_k$ such that     each  the restrictions
\[
\Psi: \Psi^{-1}(Y_k)\ra Y_k
\]
is definably trivial.

\noindent $\bullet$ (\emph{Definability of Euler characteristic.}) Suppose  $(S_\lambda)_{\lambda\in \Lambda}$ is a  definable family  of compact tame sets. Then the map
\[
\Lambda\ni \lambda\mapsto \chi(S_\lambda)=\mbox{the Euler characteristic of $S_\lambda$}\in \bZ
\]
is definable.  In particular, the  set
\[
\bigl\{ \,\chi(S_\lambda);\;\;\lambda\in \Lambda,\bigr\}\subset \bZ
\]
is finite.

\noindent $\bullet$ (\emph{Scissor equivalence principle}) Suppose $S_0, S_1$ are two tame sets. We say that they are \emph{scissor equivalent} if there exist a tame bijection  $F: S_0\ra S_1$. (The bijection $F$ \emph{need not be continuous}.) Then $S_0$  and $S_1$ are scissor equivalent if and only if  they have the same dimension and the same Euler characteristic.

\noindent $\bullet$ (\emph{Crofton formula.}, \cite{BK}, \cite[Thm. 2.10.15, 3.2.26]{Feder}) Suppose $E$ is an Euclidean space, and  denote by $\Graff^k(E)$ the Grassmannian of affine subspaces of codimension $k$ in $E$.  Fix an invariant measure $\mu$ on $\Graff^k(E)$.  $\mu$ is unique up to a multiplicative constant. Denote by $\eH^k$ the $k$-dimensional Hausdorff measure. Then there exists a constant $C>0$, depending only on $\mu$, such that for every compact, $k$-dimensional  tame subset $S\subset E$ we have
\[
\eH^k(S)= C\int_{\Graff^k(E)} \chi(L\cap S) d\mu(L).
\]

\noindent $\bullet$ (\emph{Finite volume.}) Any compact $k$-dimensional tame set has finite $k$-dimensional Hausdorff measure.

\noindent $\bullet$ (\emph{Tame quotients.}) Suppose $X$ is a   tame set, and $E\subset X\times X$ is a tame subset   defining an equivalence relation on $X$.    We assume that the natural projection $\pi: E\ra X$  is  \emph{definable proper}, i.e., for any compact  tame subset $K\subset  X$ the preimage $\pi^{-1}(K)\subset E$ is compact.  Then the quotient space $X/E$  can be realized as a tame set, i.e.,  there exists a tame  set $Y$,  and a tame continuous surjective map $p: X\ra Y$ satisfying the following properties:

\begin{enumerate}
\item[($Q_1$)] $p(x)=p(y)\Longleftrightarrow (x,y)\in E$.
\item[($Q_2$)]  $p$ is definable proper.
\end{enumerate}

The pair $(Y,p)$ is called the \emph{definable quotient} of  $X$ mod $E$. It is a quotient in the category of tame sets  and  tame continuous map in the sense that,  for   any tame continuous function  $f: X\ra  Z$ such that $(x,y)\in E\Longrightarrow f(x)=f(y)$,  there exists a unique  tame continuous map $\bar{f}: Y\ra Z$ such that the diagram below is commutative.
\[
\begin{diagram}
\node{X}\arrow{e,t}{f} \arrow{s,l}{p}\node{Z}\\
\node{Y}\arrow{ne,b,..}{\bar{f}} \node{}
\end{diagram}
\]

\qed

\smallskip

In the sequel we will work \emph{exclusively} with the tame structure  $\pfc$.   We will refer to the $\pfc$-definable sets (functions) as \emph{tame} sets (or functions), or \emph{definable} sets (functions).

%% file: tameflow2.tex
We  can now introduce the   subject of our investigation.

\begin{definition}   A \emph{tame flow} on a tame set $X$  is a continuous flow
\[
\Phi:\bR\times X\ra X,\;\;\bR\times X\ni (t,x)\ra \Phi_t(x),
\]
such that  $\Phi$ is a tame map. \qed

\end{definition}

If $\Phi$ is a tame flow on a tame set $X$, we denote by $\Cr_\Phi$  the set of stationary points of the flow.   Observe that $\Cr_\Phi$ is   a tame subset of $X$.

\begin{definition}  Suppose $\Phi$ is a tame flow on the  tame set $X$. Then a \emph{tame Lyapunov} function for $\Phi$ is a tame continuous function $f: X\ra \bR$ which decreases strictly along the  nonconstant trajectories of $\Phi$, and it is constant on the path  components of $\Cr_\Phi$.
\qed
\end{definition}

\begin{proposition} (a)If $\Phi$ is a tame flow on the tame set $X$, and $F: X\ra Y$ is a tame homeomorphism then the conjugate
\[
\Psi_t= F\circ \Phi_t \circ F^{-1}: Y\ra Y
\]
is also a tame flow.

\noindent (b) If $\Phi$ is a tame flow on the tame set $X$, and
$\Psi$ is a tame flow on the tame set $Y$, then the product  flow on $X\times Y$,
\[
\Phi\times \Psi: \bR\times X\times Y\ra X\times Y,\;\; (t,x,y)\mapsto (\Phi_t(x),\Psi_t(y))
\]
is tame.   Moreover, if $f$ is  a tame Lyapunov function for $\Phi$, and $g$ is a tame Lyapunov function for $\Psi$, then
\[
f\boxplus g: X\times Y\ra \bR,\;\; f\boxplus g(x,y)=f(x)+g(y),
\]
is a tame Lyapunov function for $\Phi\times \Psi$.

\noindent (c) If  $\Phi$ is a tame flow, then  its  opposite $\tilde{\Phi}_t:=\Phi_{-t}$ is also a tame flow.

\noindent (d) If $\Phi$ is a tame  flow on the tame space $X$ and $Y$ is a $\Phi$-invariant tame subspace then the restriction of $\Phi$ to $Y$ is also a tame flow.

\noindent (e) Suppose  $X$ is a tame set, and $Y_1, Y_2$ are compact tame subsets.  Suppose $\Phi^k$ is a tame flow on $Y_k$, $k=1,2$, such that $Y_1\cap Y_2$ is $\Phi^k$ invariant, $\forall k=1,2$, and
\[
\Phi^1|_{Y_1\cap Y_2}=\Phi^2|_{Y_1\cap Y_2}.
\]
Then there exists a unique tame  flow $\Phi$ on $X$ such that
\[
\Phi|_{Y_k}=\Phi^k,\;\;k=1,2.
\]
Moreover, if $f_k: Y_k\ra \bR$, $k=1,2$  is a tame Lyapunov  function for $\Phi^k$  and
\[
f_1|_{Y_1\cap Y_2}=f_2|_{Y_1\cap Y_2}
\]
then   the function
\[
f_1 \# f_2 : X\ra \bR,\;\;\;(f_1\# f_2)(x)=\begin{cases}
f_1(x) & x\in Y_1\\
f_2(x) & x\in Y_2
\end{cases}
\]
is a tame Lyapunov function for $\Phi$.

\label{prop: tame-flow}
\end{proposition}

\noindent {\bf Proof}\hspace{.3cm} We prove only  (a). The  map $\Psi:\bR\ra X\ra X$, $(t,x)\mapsto F\circ \Phi_t(F^{-1}(x))$ can be written as the composition of tame maps
\[
\bR\times Y \stackrel{\one \times F^{-1}}{\Lra} \bR\times  X\stackrel{\Phi}{\Lra} X\stackrel{F}{\ra} Y. \proofend
\]

\begin{ex} The translation flow  on $\bR$ given by
\[
T_t(x)=x+t,\;\;\forall t, x\in \bR
\]
is tame since its graph is the graph of   $+: \bR\times \bR\ra \bR$.  The  identity  $I_{\bR}:\bR\ra \bR$ is a tame Lyapunov function for the opposite flow.\qed
\label{ex: trans}
\end{ex}

\begin{ex} Let $X=[0,1]$, and consider the flow $\Phi$ on $X$ generated by the vector field
\[
\xi=x(x-1)\pa_x.
\]
The function $t\mapsto x(t)=\Phi_t(x_0)$ satisfies the initial value problem
\[
\dot{x}=x(x-1),\;\;x(0)=x_0.
\]
If $x_0\in \{0,1\}$ then $x(t)\equiv x_0$. If $x_0\in (0,1)$     then we deduce
\[
\frac{dx}{x(x-1)}=dt\Longleftrightarrow\frac{dx}{x} -\frac{d(1-x)}{1-x}=-dt
\]
so that
\[
\log \frac{x}{1-x} -\log \frac{x_0}{(1-x_0)}=-t.
\]
Hence
\begin{equation}
\frac{x}{1-x}= r(x_0,t):=e^{-t} \frac{x_0}{1-x_0}\Longleftrightarrow x(t)= \frac{r(x_0,t)}{1+r(x_0,t)}=\frac{e^{-t}x_0}{1-x_0+e^{-t}x_0}.
\label{eq: exact-sol}
\end{equation}
This shows that $\Phi$ is tame and  its restriction to $(0,1)$ is  tamely conjugate to the translation flow.  The identity function $[0,1]\ra [0,1]$ is a Lyapunov function for this flow. We will refer to $\Phi$ as the \emph{canonical downward flow} on $[0,1]$.\qed
\label{ex: 1-flow}
\end{ex}

\begin{ex} Consider the unit circle
\[
S^1=\{ (x,y)\in \bR^2;\;\; x^2+y^2=1\}.
\]
The height  function $h_0: S^1\ra \bR$,  $h_0(x,y)=y$, is a real analytic Morse function on $S^1$.  Define
\[
U^+ := S^1\cap \{x>0\},\;\; U^-:=S^1\cap \{x<0\}.
\]
Along $U^+$ we can use $y$ as  coordinate, and  we have $d(x^2+y^2)=0$, so that
\[
dx=-\frac{y}{x} dy\Longrightarrow  dx^2+dy^2= \frac{1}{1-y^2} dy^2.
\]
The gradient of $h_0$ with respect to the round metric $\frac{1}{1-y^2} dy^2$ is then $\xi_0:=(1-y^2)\pa_y$ so that the descending gradient flow of $h$ (with respect to this  metric) is given in the coordinate $y$ by
\[
\dot{y}= -(1-y^2).
\]
This flow is tamely conjugate to the flow in  Example \ref{ex: 1-flow} via the linear increasing  homeomorphism $[-1,1]\ra [0,1]$.   Thus the gradient flow of the  height function on the  round circle is tame. Note that this flow  is obtained by gluing two copies of the standard decreasing  flow on $[0,1]$.
\qed
\label{ex: height}
\end{ex}

\begin{ex}[\emph{A simple non tame flow}] Consider the rotational flow on the unit circle
\[
R: \bR\times S^1\ra S^1,\;\; R_t(e^{\ii\theta})= e^{\ii(t+\theta)}.
\]
This  flow is not tame because the set
\[
A=\bigl\{ t\in \bR;\;\; R_t (1) =e^{\ii t}=1\,\bigr\}=2\pi\bZ
\]
is not tame.

We deduce from this simple  example that a tame flow cannot have nontrivial periodic orbits because the restriction of  the flow to such an orbit is  tamely equivalent to the rotation flow  which  is not tame. This contradicts Proposition \ref{prop: tame-flow}(d).
\qed
\label{ex: rotation}
\end{ex}

\begin{ex}[\emph{A tame flow with no Lyapunov functions}] Consider the vector field $V$ in the plane given by
\[
V= (x^2+|y|)\pa_y.
\]
Observe that $V$ has a unique  zero located at the origin. The  flow lines are  the solutions of
\[
\dot{x}=0,\;\;\dot{y} =(x^2+|y|) ,\;\; x(0)=x_0,\;\;y(0)=y_0.
\]
Note that $y(t)$ increases along the flow lines. Thus, if $y_0\geq 0$, we deduce
\[
\dot{y}= x_0^2+y\Longrightarrow \frac{d}{dt}(e^{-t}y)= e^{-t}x_0^2\Longrightarrow e^{-t}y(t)-y_0 =x_0^2(1-e^{-t})
\]
so that
\[
y(t) =e^ty_0+x_0^2(e^t-1).
\]
If $y_0<0$ then   while $y<0$ we have
\[
\dot{y}+y= x_0^2\Longrightarrow e^ty +|y_0|=x_0^2(e^t-1).
\]
Thus
\[
y(t)=0\Longleftrightarrow e^tx_0^2=x_0^2+|y_0|\Longrightarrow t=T(x_0, y_0):=\log(x_0^2+|y_0|)-\log x_0^2
\]
We deduce that if $y_0<0$ we have
\[
y(t)=\left\{
\begin{array}{rcl}
x_0^2(1-e^{-t})+y_0e^{-t} & {\rm if} & t\leq T(x_0,y_0)\\
x_0^2(e^{t-T(x_0,y_0)} -1) & {\rm if}  & t> T(x_0,y_0)
\end{array}
\right..
\]
The trajectories of this  flow are depicted in the top half of Figure \ref{fig: tame5}
\begin{figure}[h]
\centerline{\epsfig{figure=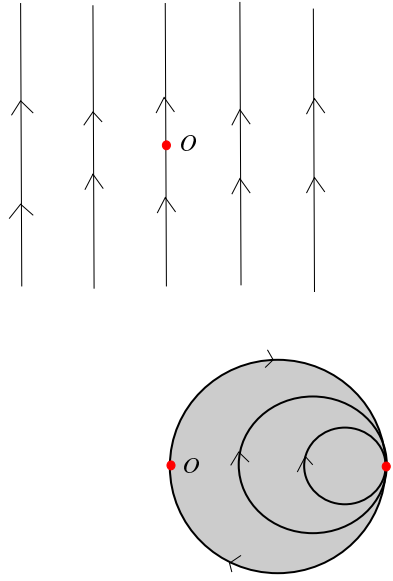,height=3in,width=2in}}
\caption{\sl  A tame  flow   with lots of homoclinic orbits.}
\label{fig: tame5}
\end{figure}

From the above description it follows immediately that this flow is tame  and extends to a tame flow  on $S^2$,  the one-point compactification of the plane.     The  flow on the eastern hemisphere ($X\geq 0$) is depicted at the bottom of Figure \ref{fig: tame5}. Observe that all but two of the orbits of this flow are homoclinic so that this flow does not admit Lyapunov functions. \qed
\end{ex}

\begin{ex}[\emph{The cone construction}]  Suppose  $\Phi$ is a   tame flow over a compact tame set $X\subset \bR^N$. We form the \emph{cone over $\Phi$} as follows.

First, define  the cone over $X$ to be the tame space $C(X)\subset [0,1]\times \bR^N\subset \bR^{N+1}$ defined as the  definable quotient
\[
[0,1]\times X/\{1\}\times X.
\]
The time  $1$-slice $\{1\}\times X$ is mapped to  the   \emph{vertex} of the cone,  denoted by $\ast$.   Denote  by $\pi: [0,1]\times X \ra C(X)$ the natural projection.  Observe that   $\pi$ is a bijection
\[
[0,1)\times X\ra C(X)^*=C(X)\setminus \{\ast\}.
\]
We thus have two maps
\[
\si: C(X)^*\ra X,\;\;\alpha: C(X)\ra [0,1].
\]
called the \emph{shadow}, and respectively \emph{altitude}.    Any point  on the cone, other than the vertex, is uniquely determined  by its shadow and altitude.

The product of the standard decreasing flow $\Psi$ on $[0,1]$ with the flow $\Phi$ on $X$ produces a flow on $[0,1]\times X$ which descends to a flow on the cone $C(X)$ called the \emph{downward cone} of $\Phi$ which we denote  by $C^\Phi$. The vertex is a stationary point of this  flow.   If $p\in C(X)^*$ then,   to  understand  the flow line $t\mapsto C^\Phi_t(p)$, it suffices to keep track of the evolution of  its shadow and its  height. The shadow of $C^\Phi_t(p)$ is $\Phi_t\si(p)$, while the height of $ C^\Phi_t(p)$ is $\Psi_t h(p)$.

Observe that if  $f$ is a     Lyapunov function for  $\Phi$ on $x$, then for every positive constant $\lambda$ the function
\[
f_\lambda: C(X)\ra \bR,\;\; f_\lambda(x)=\begin{cases}
\lambda & x=\ast\\
\lambda\alpha(x)+ \bigl(1-\alpha(x)\,\bigr) f(\si(x)) & x\neq \ast
\end{cases}
\]
is a Lyapunov function for $C^\Phi$.\qed

\label{ex: cone}
\end{ex}

\begin{ex}[\emph{The  canonical tame flow on an affine  simplex}]  We want to investigate the cone construction in a very special case.  Suppose $E$ is a finite dimensional  affine space. For every  subset $V\subset E$ we denote by $\Aff(V)$ the affine subspace spanned by $V$.  The set $V$ is called \emph{affinely independent } if $\dim \Aff (V)=\# V -1$.

If $V=\{v_0,\dotsc, v_k\}$,   and $\dim \Aff(V)=k$ we define
\[
[V]=[v_0,\dotsc , v_k] := {\rm conv}\,(\{v_0,\dotsc, v_k\}),
\]
where "conv" denotes the convex hull operation.   We  will refer to $[v_0,\dotsc, v_k]$ as the affine $k$-simplex with vertices $v_0,\cdots, v_k$. Its \emph{relative interior}, denoted by $\Int[v_0,\dotsc,v_k]$ is defined by
\[
\Int[v_0,\dotsc,v_k]:=\bigl\{\,\sum_{i=0}^k t_iv_i;\;\;t_i>0,\;\;\sum_{i=0}^kt_i=1\,\bigr\}.
\]
Given  a linearly ordered, affinely independent  finite subset  of $E$   we can associate in a canonical fashion a tame flow on the affine simplex spanned by this set.

Fix an affine $k$-simplex $S$ in the affine space $E$ with vertex set $V$.  A linear ordering   on   $V$ is equivalent to a bijection
\[
\{0,1,\cdots, k\}\ra V,\;\; i\mapsto v_i\;\;\mbox{so that}\;\; v_i<v_j \Longleftrightarrow i<j.
\]
Recall the \emph{affine cone construction}.

Let $Y$ be a subset in an affine space $E$, and $v\in E\setminus \Aff(Y)$. The \emph{cone on $Y$ with  vertex $v$} is the set
\[
C_v(Y):=\bigl\{ x=(1-t)v+ t y=v+t(y-v);\;\; t\in[0,1] ,\;\; y\in Y\,\bigr\}.
\]
In other words,   $C_v(Y)$ is the union of all segments joining $v$ to a point $y\in Y$.    Note that since $v\in E\setminus\Aff(Y)$ two such segments have only the vertex $v$ in common. This means that any point $p=C_v(Y)\setminus \{v\}$ can be written uniquely  as an affine combination
\[
p=v+t(y-v),\;\; t\in(0,1],\;\; y\in Y.
\]
If $S=[v_0,\dotsc, v_k]$ is an affine $k$-simplex, then
\[
[v_0,\dotsc,v_i,v_{i+1}]=C_{v_{i+1}}([v_0,\cdots,v_i])
\]
so that
\[
S_k= C_{v_k}\circ \cdots\circ C_{v_1}(\{v_0\}):= C_{v_k}\bigl(\cdots  C_{v_1}(\{v_0\})\cdots \bigr).
\]
The cone construction extends to  sets equipped with vector fields.

Suppose  $Y\subset E$, $v\in E\setminus \Aff(Y)$, and $Z:Y\ra TE$ is a  vector field on $Y$.   Temporarily, we impose no regularity constraints on $Z$. A priori, it could even be discontinuous. Define
\[
\hat{Z}=C_v(Z): C_v(Y)\ra TE,
\]
by setting for $t\in[0,1]$, and $y\in Y$,
\[
\hat{Z}\bigl(\, v+ t(y-v)\,\bigr)= (1-t)\cdot t(y-v)+ t Z(y).
\]
Note that $\hat{Z}(v)=0$ and $\hat{Z}(y)=Z(y)$, $\forall y\in Y$.   We let
\[
 S_i:=[v_0,\dotsc, v_i],
 \]
 and  define
\[
Z_i:= C_{v_i}\circ \cdots C_{v_1}(\vec{0}),
\]
where $\vec{0}$ denotes the trivial vector field on the set $\{v_0\}$.   By construction we have
\[
Z_{i+1}|_{S_i}= Z_i,\;\; Z_i(v_j)=0,\;\;\forall j\leq i.
\]
Observe that along the segment $[v_0,v_1]$ we have
\[
Z_1(v_1+t(v_1-v_0))=  -t(1-t)\overrightarrow{v_0v_1}.
\]
Its flow is  the canonical  downward flow on a segment and it is depicted in  Figure \ref{fig: 1}.
\begin{figure}[h]
\centerline{\epsfig{figure=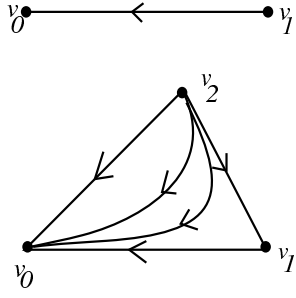,height=1.5in,width=1.5in}}
\caption{\sl  Morse structures on low dimensional simplices.}
\label{fig: 1}
\end{figure}

To understand the nature of  the vector fields $Z_i$ we argue   inductively.  Let $p\in[v_0,\cdots, v_k,v_{k+1}]$, $p\neq v_{k+1}$, and denote  by $q$ the  intersection of the line $v_{k+1}p$ with $[v_0,\cdots, v_k]$  (see Figure \ref{fig: 2}). If $(t_0,\dotsc, t_{k+1})$ denote the barycentric coordinates of $p$ in $S_{k+1}$, and $(s_0,\dotsc,s_k)$ denote the barycentric coordinates of $q$ in $S_k$, then
\[
s_i=\frac{t_i}{t_0+\cdots+t_k}=\frac{t_i}{1-t_{k+1}},\;\;0\leq i \leq k
\]
and
\[
(p-v_{k+1})= (1-t_{k+1}) (q-v_{k+1}).
\]
Then
\[
Z_{k+1}(p)=  t_{k+1}(1-t_{k+1}) (q-v_{k+1}) + (1-t_{k+1}) Z_k(q)
\]
Since $S_k$ is described in $S_{k+1}$ by $t_{k+1}=0$ and $Z_{k+1}|_{S_k}=Z_k$ we can rewrite the last equality as
\[
Z_{k+1}(t_0,\dotsc, t_k, t_{k+1})= t_{k+1}(1-t_{k+1}) \Biggl\{\,\Bigl(\sum_{i=0}^k\frac{t_i}{1-t_{k+1}} v_i\Bigr)-v_{k+1}\,\Biggr\}
\]
\[
+ (1-t_{k+1})Z_{k}\Bigl(\frac{t_0}{1-t_{k+1}},\cdots , \frac{t_k}{1-t_{k+1}}\,\Bigr).
\]
This  shows inductively that $Z_k$ is Lipschitz continuous, and even smooth on the relative interiors of the faces of $S_k$.

Denote by  $\Phi^k_t$ the (local) flow defined by $Z_k$. For any vector $\vec{\lambda}=(\lambda_0,\dotsc,\lambda_k)\in \bR^{k+1}$ such that
\[
\lambda_0<\lambda_1<\cdots < \lambda_k,
\]
we define $f_{\vec{\lambda}}: S_k\ra \bR$ to be the unique affine function on $S_k$ satisfying
\[
f_{\vec{\lambda}}(v_i)=\lambda_i,\;\;\forall i=0,1,\dotsc, k.
\]
We want to show that   for  every $k\geq 0$ the following hold.

\begin{itemize}

\item[{\bf Fact 1.}] The flow $\Phi_t^k$  exists for all $t$ on $S_k$, it is tame, and it is of the triangular type (\ref{eq: pfaff}).

\item[{\bf Fact 2.}] The linearization of $Z_k$ at   a vertex  $v_{\ell}$, $\ell= 0,1,\dotsc, k$   is diagonalizable, its spectrum  is $\{-1,1\}$ and the eigenvalue $1$ has multiplicity $\ell$.

\item[{\bf Fact 3.}] The function $f_{\vec{\lambda}}$ is a Lyapunov function for $\Phi^k$.

\end{itemize}

\noindent {\bf Fact 1.} To show that the flow $\Phi_t^k$ is tame we argue by induction over $k$. The case $k=1$ follows from Example \ref{ex: 1-flow}.    For the inductive step  we fix a vertex $u$ of $S_k$, and relabel the other $u_1,\dotsc, u_k$.

We think of $u$ as the origin  of the  affine space $\Aff(S_{k+1})$, and we introduce  the  vectors
\[
\vec{e}_i:=\overrightarrow{uu_i}=u_i-u,\;\;\vec{e}_{k+1}:=\overrightarrow{uv_{k+1}}=v_{k+1}-u.
\]
 These define   linear coordinates $(x_1,\dotsc, x_k,x_{k+1})$  on $\Aff(S_{k+1})$ so that
 \[
 \Aff(S_k)=\{x_{k+1}=0\}.
 \]
 We say that these are the \emph{linear coordinates determined by the vertex $u$}.

  Consider the point $p\in S_{k+1}\setminus v_{k+1}$ with linear coordinates
\[
p\longleftrightarrow(x_1,\dotsc, x_k, x_{k+1}).
\]
 Denote by $p'\in S_k$ the projection of $p$ on $S_k$ parallel to $e_{k+1}$,  and by $q$ the intersection of the line $v_{k+1}p$ with $S_k$ (see Figure \ref{fig: 2}). We say that $q$ is the \emph{shadow} of $p$. Then $p'$ has coordinates
 \[
 p'\longleftrightarrow(x_1,\dotsc, x_k,0),
 \]
while the shadow $q$ has coordinates
  \[
 q\longleftrightarrow \bigl(\, \frac{x_1}{1-x_{k+1}},\cdots, \frac{x_k}{1-x_{k+1}},0\,\bigr).
  \]
\begin{figure}[h]
\centerline{\epsfig{figure=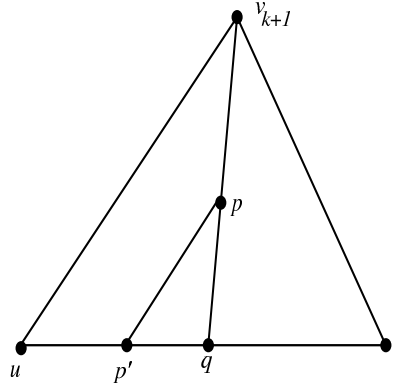,height=2in,width=2in}}
\caption{\sl  Dissecting the cone construction.}
\label{fig: 2}
\end{figure}
 Since  $\overrightarrow{v_{k+1}p}= (1-x_{k+1})\overrightarrow{v_{k+1}q}$  we deduce
\[
Z_{k+1}(x_1,\cdots, x_k,x_{k+1})= x_{k+1}(1-x_{k+1}) \overrightarrow{v_{k+1}q} +(1-x_{k+1}) Z_k(q)
\]
\[
=x_{k+1}(1-x_{k+1}) \Biggl\{\,\Bigl(\sum_{i=1}^k \frac{x_i}{1-x_{k+1}}\vec{e_i}\,\Bigl)-\vec{e}_{k+1}\,\Biggr\}
\]
\[
+(1-x_{k+1})Z_k\bigl(\, \frac{x_1}{1-x_{k+1}},\cdots, \frac{x_k}{1-x_{k+1}}\,\bigr)
\]
\[
=-x_{k+1}(1-x_{k+1})\vec{e}_{k+1}+x_{k+1}\sum_{i=1}^kx_i \vec{e}_i +(1-x_{k+1})Z_k\bigl(\, \frac{x_1}{1-x_{k+1}},\cdots, \frac{x_k}{1-x_{k+1}}\,\bigr).
\]
If we write
\[
Z_k =\sum_{i=0}^k Z_k^i \vec{e}_i
\]
then we deduce
\[
Z_{k+1}^{k+1}= x_{k+1}(x_{k+1}-1),\;\; Z_{k+1}^i= x_{k+1}x_i + (1-x_{k+1})Z_k^i\Bigl(\frac{x_1}{1-x_{k+1}},\cdots, \frac{x_k}{1-x_{k+1}}\,\Bigr).
\]
This  shows inductively that,  away from the vertex $v_{k+1}$  of the simplex $S_{k+1}$, the  vector field $Z_{k+1}$ has upper triangular form in the coordinates $x_1,\cdots, x_{k+1}$, i.e.,  the $i$-th component  $Z_{k+1}^i(x_1,\cdots, x_{k+1})$ depends only on the variables  $x_i,\cdots, x_{k+1}$.  This  defines a system   differential equations of the pfaffian  type  (\ref{eq: pfaff}).

We want to prove that  the vector field  $Z_{k+1}$  determines a \emph{globally defined} flow on the simplex $S_k$.    From the inductive assumption we deduce  that  for any $k$ simplex   with linearly oriented   vertex set the  corresponding vector field determines a globally    defined tame flow. Consider a point
\[
p \in S_{k+1}\setminus   \{v_{k+1}\}.
\]
We will use the linear coordinates $(x_1,\dotsc, x_{k+1})$ determined by the vertex $v_0$.  Assume that the linear coordinates of $p$ are
\[
p=(a_1,\dotsc, a_{k+1})
\]
The flow line of $Z_{k+1}$  through $p$  is a path
\[
t\stackrel{\gamma}{\longmapsto} (x_1(t),\dotsc, x_{k+1}(t))
\]
 satisfying the initial value problem
\begin{equation}
\dot{x}_{k+1}=-x_{k+1}(1-x_{k+1}),\;\; x_{k+1}(0)=a_{k+1}
\label{eq: last}
\end{equation}
\begin{equation}
\dot{x}_i = x_{k+1}x_i + (1-x_{k+1}) Z_k^i \bigl(\frac{x_i}{(1-x_{k+1})},\cdots, \frac{x_k}{(1-x_{k+1})} \,\bigr),\;\;x_i(0)=a_i.
\label{eq: lasti}
\end{equation}
For simplicity we write $x:=x_{k+1}$. We  introduce the  \emph{shadow}  coordinates
\[
s_i=\frac{x_i}{(1-x_{k+1})}\Longleftrightarrow x_i=s_i(1-x_{k+1}),\;\;i=1,2,\cdots, k.
\]
The projection of the path $\gamma(t)$  from the vertex $v_{k+1}$ onto the face  $[v_0,\dotsc, v_k]$ is given  in linear coordinates by the shadow path $t\mapsto(\, s_1(t),\dotsc, s_k(t)\,)$.

Since $\dot{x}=-x(1-x)$ we  deduce
\[
\frac{d}{dt} x_i =\frac{d}{dt} (s_i(1-x)) = \dot{s}_i(1-x)-s_i\dot{x}=\dot{s}_i(1-x)+s_ix(1-x).
\]
Using this in (\ref{eq: last}) and (\ref{eq: lasti})  we deduce
\begin{equation}
\dot{x}=-x(1-x)\;\;\dot{s}_i= Z_k^i(s_i,\cdots, s_k),\;\;j=1,\dotsc, k.
\label{eq: last2}
\end{equation}
This   computation, coupled with the inductive assumption   show that $\Phi^{k+1}_t(p)\in S_{k+1}$, $\forall t\in \bR$.

The flow $\Phi_t^{k+1}$ can be given the following simple interpretation.  Denote by $\Phi^k_t$ the flow on $[v_0,\cdots, v_k]$, and  by $s$ the   shadow map
\[
s: {\rm Int}\,[v_0,\dotsc, v_{k+1}]\ra {\rm Int}\, [v_0,\dotsc, v_k],\;\; s(p):= v_{k+1}p\cap [v_0,\dotsc, v_k],
\]
where $v_kp$ denotes the line  passing through $v_{k+1}$ and $p$.  We  assume $v_0$ is the origin  of our affine space so   we can describe a point in the simplex $[v_0,v_1,\dotsc, v_i]$  by its linear coordinates $(x_1,\dotsc, x_i)$. Given
\[
p_0=(a_1,\dotsc, a_{k+1})\in {\rm Int}\,[v_0,\dotsc, v_{k+1}],
\]
 we set $q_0=s(p_0)$ and then we have
 \[
 \Phi^{k+1}_t(p_0)=   x(t) v_{k+1} + \bigl(\, 1-x(t)\,\bigr) \Phi_t^k(s(p_0)),\;\; x(t)= \frac{e^{-t} a_{k+1}}{1-a_{k+1} +e^{-t}a_{k+1}}.
 \]
The path $\Phi^{k+1}_t(p_0)$ can be  visualized using a natural moving frame.

Denote by  $s_0$ the shadow of $p_0$.   Now let $s_0$  go with the flow $\Phi^k$, $s(t)=\Phi_t^k(s_0)$. The point $p(t)=\Phi^{k+1}(p_0)$ lies on the segment $[s(t), v_{k+1}]$. If we   affinely identify this segment with the segment $[0,1]$ so that $s(t)\longleftrightarrow 0$ and $v_{k+1}\longleftrightarrow 1$, then   the motion of the point $p(t)$ along the (moving) segment $[s(t),v_{k+1}]$ is mapped  to the motion on the unit segment $[0,1]$ governed by the canonical downward flow on $[0,1]$. In other words,   the flow $\Phi^{k+1}$ is the negative  cone on the flow $\Phi^k$. This proves that $\Phi^{k+1}$ is tame.

To see how this works in  concrete examples, suppose   $S_2$ is the $2$-simplex
\[
\{ (x,y)\in \bR^2 ;\; 0\leq x,y,\;\; x+y\leq 1\}
\]
with vertices $v_0=(0,0)$, $v_1=(1,0)$, $v_2=(0,1)$. Consider the point $p_0=(x_0,y_0)$ in the interior of this simplex.   If $\Phi$ is the  flow defined by the  vector field $Z_2$ then
\[
\Phi_t(x_0,y_0)=\Bigl(\,\bigl(1-y(t)\,\bigr)\frac{e^{-t}x_0}{ 1-x_0-y_0+e^{-t}x_0},\; y(t)\,\Bigl),\;\;y(t)= \frac{e^{-t}y_0}{1-y_0+e^{-t}y_0}.
\]

\noindent {\bf Fact 2.} The statement about the linearization of $Z_k$ at the vertices of $S_k$ is again proved by induction. The statement is obvious for $k=1$.  For the inductive step, denote by $u$ one of the vertices of $S_k$, and  label the remaining ones by $u_1,\dotsc, u_k$.  We again think of $u$ as the origin of $\Aff(S_{k+1})$ and as such we  obtain  a  basis
\[
\vec{e}_i =u_i-u,\;\;\vec{e}_{k+1}=v_{k+1}-u,
\]
and linear coordinates $(x_1,\dotsc, x_{k+1})$. The point $u$  has linear coordinates  $x_i=0$, $0\leq i\leq k+1$ in $S_{k+1}$.  Denote by $\nabla$ the trivial connection on the tangent bundle of $\Aff(S_{k+1})$. For $i=0,\dotsc ,k$ we have
\[
\nabla_{e_i} Z_{k+1}(x_1,\dotsc, x_k) =x_{k+1} \vec{e}_i+ \nabla_{e_i} Z_k\bigl(\, \frac{x_1}{1-x_{k+1}},\cdots, \frac{x_k}{1-x_{k+1}}\,\bigr)
\]
and
\[
\nabla_{e_{k+1}}Z_{k+1}=\sum_{i=1}^k\bigl(\, x_i\vec{e}_i +(2x_{k+1}-1)\vec{e}_{k+1}\,\bigr)-Z_k\bigl(\, \frac{x_1}{1-x_{k+1}},\cdots, \frac{x_k}{1-x_{k+1}}\,\bigr)
\]
\[
+\frac{1}{1-x_{k+1}} \sum_{i=1}^k x_i\nabla_{e_i}Z_k\bigl(\, \frac{x_1}{1-x_{k+1}},\cdots, \frac{x_k}{1-x_{k+1}}\,\bigr)
\]
Observe that at $(x_1,\cdots, x_{k+1})=0$ we have
\[
\nabla_{e_i} Z_{k+1}= \nabla_{e_i} Z_k,\;\;\;1\leq i\leq k,
\]
and
\[
\nabla_{e_{k+1}}Z_{k+1}= -e_{k+1}
\]
This proves the statement about the linearization of $Z_{k+1}$ at $u\in \{v_0,\dotsc, v_k\}$.

Finally, we want to prove that the linearization of $Z_{k+1}$ at $v_{k+1}$ is the identity.  Since $Z_k(v_i)=\vec{0}$, $\forall i=0,1,\dotsc, k$ we deduce that   at a point  $p$ on the line segment $[v_{k+1},v_i]$ given by
\[
p=v_{k+1}+(1-t) (v_i-v_{k+1}),
\]
the vector field $Z_{k+1}$ is described by
\[
Z_{k+1}(p)= t(1-t)(v_i-v_{k+1}).
\]
If we fix the origin of $\Aff(S_{k+1})$ at $v_{k+1}$, and we set $\vec{f}_i:=\overrightarrow{v_{k+1}v_i}$, then
\[
Z_{k+1}(v_{k+1} + s \vec{f}_i)= s(1-s)\vec{f}_i,\;\;\nabla_{f_i}Z_{k+1}(v_k)=  \vec{f}_i,
\]
so that the linearization of $Z_{k+1}$ at $v_{k+1}$ is the identity operator
\[
\bsI: T_{v_{k+1}} \Aff(S_{k+1})\ra T_{v_{k+1}} \Aff(S_{k+1}).
\]

\noindent {\bf Fact 3.}     We again argue by induction. The statement is  true for $k=1$. For the inductive step, denote by $(x_1,\dotsc,x_{k+1})$ the linear coordinates  on $S_{k+1}$ determined by the vertex $v_0$. In these coordinates we have
\[
f_{\vec{\lambda}} = \lambda_0 +\sum_{j=1}^{k+1}(\lambda_j-\lambda_0) x_j.
\]
If we write $x=x_{k+1}$, and again we introduce the shadow coordinates $s_j=\frac{x_j}{1-x}$, we deduce
\[
f_{\vec{\lambda}}(s_1,\dotsc, s_k,x)= \lambda_0+(\lambda_{k+1}-\lambda_0)x+ (1-x)\sum_{j=1}^k (\lambda_j-\lambda_0) s_j.
\]
If we differentiate $f_{\vec{\lambda}}(s_1,\dotsc, s_k,x)$ along  a flow line  we deduce
\[
\frac{d}{dt}f_{\vec{\lambda}}(s_1,\dotsc, s_k,x)= (\lambda_{k+1}-\lambda_0)\dot{x} -\dot{x} \sum_{j=1}^k (\lambda_j-\lambda_0) s_j+(1-x)\sum_{j=1}^k (\lambda_j-\lambda_0)\dot{s}_j.
\]
 Using (\ref{eq: last2})  we deduce
\[
\frac{d}{dt}f_{\vec{\lambda}}(s_1,\dotsc, s_k,x)
\]
\[
=-(\lambda_{k+1}-\lambda_0)x(1-x) +x(1-x)\sum_{j=1}^k (\lambda_j-\lambda_0) s_j + (1-x)\sum_{j=1}^k (\lambda_j-\lambda_0)\dot{s}_j
\]
\[
=x(1-x)\left( \sum_{j=1}^k (\lambda_j-\lambda_0) s_j- (\lambda_{k+1}-\lambda_0)\,\right)+ (1-x)\sum_{j=1}^k (\lambda_j-\lambda_0)\dot{s}_j
\]
The  first term is strictly negative because $\lambda_j  <\lambda_{k+1}$, and  on $S_k$  we have
\[
\sum_js_j\leq 1,\;\; s_j\geq 0,
\]
where at least one of these inequalities is strict. The  second  term is negative since the restriction of $f_{\vec{\lambda}}$ to the face $S_k$ is a Lyapunov function for $\Phi^k$.  \qed
\label{ex: simplicial}
\end{ex}

The above example has an important consequence.

\begin{proposition} On any compact tame space there exist tame flows with finitely many stationary points  and admitting tame  Lyapunov functions, i.e. continuous tame  functions decreasing strictly along the nonconstant orbits.
\label{prop: exist}
\end{proposition}

\proof   Suppose $X$ is a compact tame space.  Choose an affine simplicial complex $Y$  and a tame homeomorphism $F: Y\ra X$.        Denote by $V(Y)$ the vertex set of $Y$ and choose a map  $\ell: V(Y)\ra \bR$ which is  injective when restricted to the vertex set of any  simplex of $Y$.    We can now use the map $\ell$ to linearly order the vertex set of any simplex $\si$ of $Y$ by declaring
\[
u< v\Longleftrightarrow \ell(u) <\ell(v).
\]
This ordering induces as in Example \ref{ex: simplicial}  a tame flow $\Phi^\si=\Phi_t^{\si,\ell}$ on any  face $\si$ of $Y$ such that
\[
 \Phi^\tau|_\si=\Phi^\si,\;\;\forall \si\prec\tau.
 \]
Thus  the tame flows on  the faces  of $Y$    are compatible on overlaps and thus define  a tame flow  on $Y$. Note that the function  $\ell$ defines a piecewise linear function $\ell: Y\ra \bR$ which decreases strictly along the trajectories of  $\Phi$. Using the homeomorphism $F$ we obtain a tame flow $F\circ \Phi\circ F^{-1}$ on $X$. Its stationary points  correspond  via $F$ with the vertices of $Y$, and $F\circ \ell\circ F^{-1}$ is a tame Lyapunov function. \qed

\begin{ex}  Suppose $E$ is a finite dimensional real Euclidean space, and $A\in \End(E)$ is a  symmetric endomorphism. Then the linear flow
\[
\Phi^A: \bR\times E\ra E,\;\;  \Phi_t^A(x)= e^{At}x,\;\;x\in E,
\]
is a tame flow. Similarly,  the flow
\[
\Psi^A: \bR\times \End E\ra \End E,\;\; \Psi^A_t(S)=e^{At}Se^{-At},\;\; S\in \End E
\]
is a tame flow. \qed

\end{ex}

\begin{ex}   Suppose $E$ is a finite dimensional real Euclidean space, and $A\in \End(E)$ is a  symmetric endomorphism. Denote by $\Gr_k(E)$ the Grassmannian of $k$-dimensional subspaces of $E$.  For every $L\in \Gr_k(E)$ we denote by $P_L$ the orthogonal projection  onto $L$.    The map
\[
\Gr_k(E)\ni L\mapsto P_L\in \End E
\]
embeds  $\Gr_k(E)$ as a real algebraic submanifold of  $\End E$.

On $\End E$  we have  and inner product given by
\[
\lan S, T\ran =\tr (ST^*),
\]
and we denote by $|\bullet|$ the corresponding Euclidean norm on $\End E$. This inner product induces a   smooth Riemann metric on $\Gr_k(E)$.

The flow
\[
\Gr_k(E) \ni L\mapsto e^{At} L\in \Gr_k(E)
\]
is tame.  To see  this,  consider an orthonormal basis of eigenvectors of $A$, $e_1,\dotsc, e_n$, $n=\dim E$ such that
\[
Ae_i= \lambda_i e_i,\;\;\lambda_1\geq \lambda_2\geq \cdots \geq \lambda_n.
\]
For every subset $I\subset \{1,\dotsc, n\}$ we write
\[
E_I:={\rm span}\,\{e_i,\;\;i\in I\},\;\;I^\perp:=\{1,\dotsc, n\} \setminus I.
\]
For $\# I=k$ we set
\[
\Gr_k(E)_I=\bigl\{ L\in \Gr_k;\;\;L\cap E_I^\perp=0,\bigr\}.
\]
$\Gr_k(E)_I$ is a   semialgebraic open subset  of $\Gr_k(E)$ and
\[
\Gr_k(E)=\bigcup_{\# I=k}\Gr_k(E)_I.
\]
A subspace $L\in \Gr_k(E)_I$ can be represented as the graph of a linear map $S=S_L: E_I\ra E_I^\perp$, i.e.,
\[
L = \bigl\{  x+Sx;\;\; x\in E_I\,\bigr\}.
\]
Using the basis $(e_i)_{i\in I}$ and $(e_\alpha)_{\alpha\in I^\perp}$  we can represent  $S$   as a $(n-k)\times k$ matrix
\[
S= [s_{\alpha i}]_{i\in I,\;\; \alpha\in I^\perp}.
\]
The subspaces $E_I$ and $E_I^\perp$ are $A$ invariant. Then $ e^{At}L\in \Gr_k(E)_I$,  and it is represented as the graph of the operator $ S_t= e^{At} Se^{-At}$ described by the matrix
\[
\diag(e^{\lambda_\alpha t},\;\;\alpha\in I^\perp)\cdot  S \cdot\diag(e^{-\lambda_i t},\;\;i\in I)=[e^{(\lambda_\alpha-\lambda_i)t}s_{\alpha i}]_{i\in I,\;\; \alpha\in I^\perp}.
\]
This proves that   the flow is tame.

Let us point out that this flow is the  gradient flow of  the function
\[
f_A: \Gr_k(E)\ra \bR,\;\;f_A(L)=\tr AP_L=\lan A, P_L\ran.
\]
This is a Morse-Bott function. We want to describe  a simple consequence of this fact which we will need later on.

Suppose $U$ is a subspace of $E$, $\dim U \leq k$, and define
\[
A:=P_{U^\perp}=\one_E-P_U.
\]
Then
\[
f_A(L)= \tr (P_L-P_LP_U)=\dim L - \tr(P_L P_U).
\]
On the other hand, we have
\[
|P_U-P_UP_L|^2 =\tr (P_U-P_UP_L)(P_U-P_LP_U)= \tr(P_U-P_UP_LP_U)
\]
\[
= \tr P_U-\tr P_UP_LP_U=\dim U- \tr P_U P_L.
\]
Hence
\[
f_A(L)=|P_U-P_UP_L|^2 +\dim L-\dim U,
\]
so that
\[
f_A(L) \geq  \dim L-\dim U,
\]
with equality if and only if $L\supset U$. Thus, the  set of minima of $f_A$ consists of all $k$-dimensional subspaces containing $U$. We denote this set with $\Gr_k(E)_U$. Since $f_A$ is a Morse-Bott function we deduce that
\begin{equation}
\begin{split}
\forall  j\leq k,\;\;\forall U\in\Gr_j(E)\;\;\exists C=C(U)>1,\;\; \forall L\in \Gr_k(E):\\
 \frac{1}{C}\dist(L,\Gr_k(E)_U\,)^2 \leq |P_U-P_UP_L|^2 \leq C  \dist(L,\Gr_k(E)_U\,)^2.
 \end{split}
\label{eq: dist}
\end{equation}
In a later  section we will  prove   more precise results concerning the asymptotics of this  Grassmannian flow. \qed
\label{ex: grass}
\end{ex}

%% file: tameflow3.tex
  We would like to present a few  general results concerning the long time behavior   of a tame flow.

  \begin{definition} Suppose $\Phi: \bR\times X\ra X$ is a  continuous flow on a topological space $X$. Then for every set $A\subset X$ we define
\[
\Phi_+(A)= \bigcup_{t\geq 0} \Phi_t(A)=\Phi([0,\infty)\times
A),\;\;\Phi_-(A)= \bigcup_{t\leq 0}
\Phi_t(A)=\Phi((-\infty,0]\times A)
\]
\[
\Phi(A)=\Phi(\bR\times A)=\Phi_+(A)\cup \Phi_-(A).
\]
We will say that $\Phi_\pm(A)$ is  the  forward/backward
\emph{drift} of $A$ along $\Phi$, and  that $\Phi(A)$ is the
\emph{complete drift}. \qed
\end{definition}
The next result    follows immediately from the definitions.

\begin{proposition} If $\Phi$ is a tame flow on $X$ then  for every tame subset  $A\subset X$   the sets  $\Phi_\pm(A)$ and $\Phi(A)$ are tame.\qed
\label{prop: tame-drift}
\end{proposition}

\begin{theorem} Suppose  $\Phi$ is a continuous flow  on the  tame set $X$.  Consider the flow $G_\Phi:=T\times \tilde{\Phi}\times \Phi$ on $\bR\times X\times X$, where $T$ denotes the translation flow on $\bR$ and $\tilde{\Phi}_t=\Phi_{-t}$. Denote by $\Delta_0$   the initial diagonal
\[
\Delta_0= \bigl\{ (0,x,x)\in \bR\times X\times X\,\bigr\}.
\]
The following  conditions are equivalent.

\smallskip

\noindent (a) $\Phi$ is a tame flow.

\noindent (b) The  complete drift of $\Delta_0$ along $G_\Phi$ is
a tame subspace of $\bR\times X\times X$.

\label{th: tame}
\end{theorem}

\proof (a) $\Longrightarrow$ (b).  Since $\Phi$ is tame we deduce
that $G_\Phi$ is tame and we conclude using Proposition \ref{prop:
tame-drift}.

\smallskip

\noindent (b) $\Longrightarrow$ (a). Observe that
\[
G_\Phi(\Delta_0)=\bigl\{ (t,x_0,x_1)\in \bR\times X\times X;\;\;
\exists x\in X:\;\;x_0=\Phi_{-t}(x),\;\;x_1=\Phi_t(x)\,\bigr\}.
\]
Consider the tame  homeomorphism
\[
F:\bR\times X\times X\ra \bR\times X\times X,\;\; (t,x_0,
x_1)\longmapsto (s,y_0,y_1):= (2t, x_0,x_1)
\]
and observe that $F$ maps $G_\Phi(\Delta_0)$ onto the graph of the
flow $\Phi$.  \qed

\begin{corollary} Suppose  that the flow $\Psi$ on the tame space $S$  is tamely conjugate to the translation flow on  $\bR$.  Then a flow  $\Psi$ on the tame space $X$ is tame if and only if  there exists $s_0\in S$ such that the total drift of the   diagonal
\[
\Delta_{s_0}=\{ (s_0,x,x)\in S\times X\times X\}
\]
with respect to the flow $\Psi\times \tilde{\Phi}\times \Phi$ is
tame.  \qed \label{cor: tame}
\end{corollary}

In applications $S$ will be the open semi-circle
\[
S=\{ (x,y)\in \bR^2;\;\;x^2+y^2=1,\;\;x>0\}
\]
equipped with the negative gradient   flow of the height function
$h(x,y)=y$.  As origin of $s$ we take $s_0=(1,0)$. As explained in
Example \ref{ex: height} this flow is tamely conjugate to the
translation flow on $\bR$. The following result is an immediate
consequence of  tameness.

\begin{proposition} Suppose $X$ is a tame compact set  of dimension $d$,   $S$ is the open semi-circle  equipped with the flow $\Psi$ described above, and $\Phi$ is a tame flow on $X$.  We set $s_t:=\Psi_t(s_0)$, $s_0=(1,0)\in S$. Then $\Phi$ has  finite volume, i.e.     the image of the graph of $\Phi$ via the tame  diffeomorphism
\[
\bR\times X\times X\ra S\times X\times X,\;\; (t,x_0,x_1)\mapsto
(s_t, x_0, x_1)
\]
has finite $(d+1)$-dimensional Hausdorff measure. \label{prop:
finite-volume}
\end{proposition}

\begin{proposition} Suppose $\Phi$ is a tame flow on the compact space $X$. Then there exists a positive constant  $L= L(X,\Phi)$ such that every orbit of $\Phi$ has length $\leq L$.
\end{proposition}

\noindent \proof  Consider the  roof
\[
\begin{diagram}
\node{}\node{\Gamma_\Phi\subset \bR\times X\times X}\arrow{sw,t}{\ell}\arrow{se,t}{r}\node{}\\
\node{X}\node{}\node{X}
\end{diagram},
\]
where
\[
\ell(t,x_0,x_1)=x_0,\;\;\;r(t,x_0,x_1)=x_1.
\]
This roof describes the  family of   subspaces of $X$,
$(\eO_x)_{x\in X}$, where
\[
\eO_x= r\bigl(\ell^{-1}(x)\,\bigr)=\bigl\{
\Phi_t(x);\;\;t\in\bR\,\bigr\}.
\]
We see that $\eO_x$ is the orbit of the flow through $x$, and thus
the family  of orbits is a definable  family   of tame subsets
with diameters bounded from.      The claim in the proposition now
follows  from the  Crofton formula  and the definability of Euler
characteristic. \qed

\begin{proposition} Suppose  that  $\Phi$ is a tame flow on the compact tame space $X$.       Then for every $x\in X$ the limits  $\lim_{t\ra\pm \infty} \Phi(x)$ exist and are  stationary points of the flow denoted by $\Phi_{\pm\infty}(x)$. Moreover, the resulting  maps
\[
\Phi_{\pm\infty}: X\ra X
\]
are tame.
\end{proposition}

\proof    Clearly the limits exist if $x$ is a stationary point.
Assume $x$ is not a stationary point. Then the orbit $\Phi(x)$ is
a one-dimensional tame  subset  and its frontier
\[
\Fr\Phi(x)= \clos\Phi(x) \setminus \Phi(x)
\]
is a tame,  zero dimensional, $\Phi$-invariant subset. In
particular it must be finite collection of stationary points,
$\{x_1,\dotsc, x_\nu\}$.  Choose  small, disjoint, tame, open
neighborhoods $U_1,\dotsc, U_\nu$ of $x_1,\dotsc, x_\nu$,  and set
\[
U=\bigcup_{k=1}^\nu U_k.
\]
Then the set $S=\bigl\{ t\in \bR;\;\; \Phi_t(x)\in U\}$ is a tame
open subset of $\bR$, and thus it consists of finitely many,
disjoint open  intervals, $I_1,\cdots, I_N$.  Since  the set
$\{x_1,\dotsc, x_N\}$ consists of limit points of the orbit,   one
(and only one) of these intervals,  call it $I_+$, is unbounded
from above,  and one and  only one of these intervals, call it
$I_-$,  is unbounded from  below.    Then there exist $x_\pm \in
\{x_1,\dotsc, x_\nu\}$ such that $\Phi_t(x)$ is near $x_\pm$ when
$t\in I_\pm$.  We  deduce that
\[
\Fr \Phi(x)=\{x_\pm\}\;\;{\rm and}\;\;\lim_{t\ra \pm
\infty}=x_\pm.
\]
Denote by $\Gamma^\pm$ the graph of $\Phi_{\pm\infty}$. We deduce
that
\[
(x,y)\in \Gamma^+ \Longleftrightarrow  \forall \ve >0,\;\;\exists
T>0:\;\;  {\rm dist}\, (\Phi_t(x),y) <\ve,\;\;\forall  t> T.
\]
This shows  $\Gamma^+$ is definable.  A similar argument shows
that $\Gamma^-$ is tame. \qed

\begin{definition} Suppose $\Phi$ is a tame flow on the compact tame space $X$.

\smallskip

\noindent (a) We denote by $\Cr_\Phi$ the set of stationary points
of $\Phi$  and for every $x\in \Cr_\Phi$ we  set
\[
W^+(x,\Phi):=\Phi_\infty^{-1}(x),\;\;
W^-(x,\Phi):=\Phi_{-\infty}^{-1}(x)
\]
and we say that $W^\pm(x,\Phi)$ is the \emph{stable} (resp.
\emph{unstable}) variety of $x$.

\noindent (b) For $x_0,x_1\in X$ we set
\[
C_\Phi(x_0,x_1):= W^-(x_0,\Phi)\cap W^+(x_1,\Phi)=\bigl\{ z\in
X;\;\; x=\Phi_{-\infty}(z),\;\;y=\Phi_\infty(z)\,\bigr\}.
\]
We say that   $C_\Phi(x_0,x_1)$ is the $\Phi$-\emph{tunnel} from
$x_0$ to $x_1$. Observe that all the   spaces $\Cr_\Phi$, $C_\Phi$
and $W^\pm (-,\Phi)$ are tame subspaces.\qed

\end{definition}

\begin{ex}[\emph{Tame homoclinic orbits}]       It might be tempting to believe that for a  tame flow $\Phi$ on a compact tame space $X$, and a stationary point $x$  of $\Phi$ the tunnel $C_\Phi(x,x)$ consists only of the stationary point. However it is easy to produce an example  when $C_\Phi(x,x)$ contains nontrivial orbits as well.

\begin{figure}[h]
\centerline{\epsfig{figure=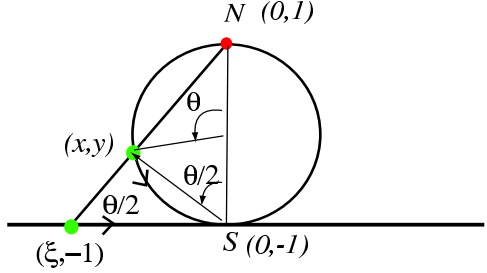,height=1.4in,width=2.5in}}
\caption{\sl  The stereographic projection of the  translation
flow.} \label{fig: tame3}
\end{figure}

Consider a round circle $C$ of radius $1$ in the Euclidean plane
and a line tangent to it (see Figure \ref{fig: tame3}) at a point
$S$. Denote by $N$ the point on the circle diametrically opposed
to  $S$.  We obtain a  flow on the line by linearly identifying it
with the real axis. Via stereographic projection from $N$ we
obtain a tame flow on the complement of $N$  in the circle. This
extends to a tame flow $\phi_t$ on the circle by declaring $N$  to
be a stationary point. The tunnel from $N$ to itself contains one
nontrivial homoclinic orbit and in fact $C(N,N)$ is the entire
circle.

Let us express  this flow in coordinates. Assume   the circle is
given by $x^2+y^2=1$ and the line is given by $y=-1$. The North
Pole has coordinates $N=(0,1)$. We use angular coordinate
$\theta=\theta(x,y)$ so  that $\theta(0,1)=0$,
$\theta(1,0)=\frac{\pi}{2}$ (see Figure \ref{fig: tame3}. Then
\[
y
=\cos\theta,\;\;x=\sin\theta,\;\;\tan(\frac{\theta}{2})=-\frac{2}{\xi}\Longrightarrow\xi=-2\cot(\theta/2)
\]
and therefore
\[
d\xi=
\frac{d\theta}{\sin^2(\theta/2)}=(1+\cot^2(\theta/2))d\theta=
(1+\xi^2/4)d\theta.
\]
Hence
\[
d\theta=\frac{4d\xi}{4+\xi^2},\;\; \frac{d}{d\xi}=
2\sin^2(\theta/2)\frac{d}{d\theta}
\]
This shows the flow is generated by  real analytic vector field on
$S^1$ with  a degenerate zero at the North pole.

Let us  notice, that  this  flow      has the \emph{Harvey-Lawson
property}:  for every  differential form $\alpha$ on $S^2$, the
pullbacks $\phi_t^*\alpha$ converge as currents to a real flat
chain as $t\ra \infty$.

Suppose  $\alpha$ is a smooth $1$-form on the circle.  We regard
$\alpha$ as a current     and we set
\[
\lan\alpha,\beta\ran =\int_C \beta\alpha.
\]
We would like to investigate the behavior of the current
$\phi_t^*\alpha$ as $t\ra \infty$. We  use the coordinates $\xi$
on $C\setminus \{N\}$.    The form $\alpha=a(\theta)d\theta$
changes to
\[
\alpha=\frac{4a(\xi)}{4+\xi^2} d\xi.
\]
In the coordinates $\xi$ the flow is described by
$\Phi_t(\xi)=\xi+t$. Then
\[
\int_C \phi_t^*\alpha\beta= \int_C \alpha \phi_{-t}^*\beta=
\int_{-\infty}^\infty \frac{4a(\xi)\beta(\xi-t)}{4+\xi^2} d\xi
\]
Using the dominated convergence theorem and the fact that
$\lim_{t\ra \infty}\beta(\xi-t)=\beta(N)$ we deduce
\[
\lim_{t\ra\infty} \int_C\phi_t^*\alpha\beta=\beta(N)\int_C\alpha
\]
so that $\phi_t^*\alpha$ converges as a current to
$\Bigl(\int_C\alpha\Bigr)\delta_N$ where $\delta_N$ is the Dirac
$0$-current concentrated at $N$. \qed

\label{ex: homoclinic}
\end{ex}

Suppose $\Phi$ is a tame flow on the compact tame space $X$.
Assume that $\Cr_\Phi$ is finite. Observe that we have a natural
action of $\bR^2$ on $\bR\times X\times X$ given by
\[
(s_0,s_1)\cdot (t,x_0,x_1):= \bigl(\, t+s_1-s_0,
\Phi_{s_0}(x_0),\Phi_{s_1}(x_1)\,\bigr).
\]
 We denote by $\Gamma\subset \bR\times X\times X$ the graph of $\Phi$ and we observe that $\Gamma$ is invariant with respect to the above action of $\bR^2$. We denote by $\Gamma^t\subset X\times X$ the  graph of $\Phi_t$, by $\bar{\Gamma}$ the closure of $\Gamma$ in $[-\infty]\times X\times X$, and by $\Gamma^{\pm \infty}$ the part of $\bar{\Gamma}$ over $\pm\infty$. Extend the  above $\bR^2$-action to $[-\infty,\infty]\times X\times X$ by setting
 \[
(s_0,s_1)\cdot (\pm\infty ,x_0,x_1):= \bigl(\, \pm\infty,
\Phi_{s_0}(x_0),\Phi_{s_1}(x_1)\,\bigr).
\]
For every subset $S\subset X\times X$ we denote by $^\ast\!S$ the
reflection of $S$ in the diagonal, i.e.
\[
^\ast\!S=\{ (x_0,x_1)\in X\times X;\;\;(x_1,x_0)\in S\}
\]

\begin{proposition}   $\bar{\Gamma}$ and $\Gamma^{\pm \infty}$  are   tame, $\bR^2$-invariant subsets of $[-\infty,\infty]\times X \times X$. Moreover,
\[
\Gamma^{-\infty}=^\ast\!\Gamma^\infty,
\]
\[
\{x\}\times W^-(x,\Phi), \;\;W^+(x)\times \{ x\}\subset
\Gamma^\infty,
\]
\[
\dim\Gamma^{\pm\infty}=\dim X = \dim \bar{\Gamma} -1.
\]
\end{proposition}

\noindent {\bf Proof}\hspace{.3cm}  The first part follows from
the tameness of $\Gamma$ and the continuity of the $\bR^2$-action.
Suppose $(x_0,x_1)\in \Gamma_\infty$. Then there exist sequences
$(x_n)\subset X$, $t_n\ra \infty$ such that
\[
(x_n,\Phi_{t_n}x_n))\ra (x_0,x_1).
\]
Let $y_n:= \Phi_{t_n}(x_n)$. Then $x_n =\Phi_{-t_n}(y_n)$, and we
deduce
\[
(x_1,x_0)=\lim_{t_n\ra \infty} (y_n,\Phi_{-t_n}(y_n)) \in
\Gamma_{-\infty}.
\]
Let  $y\in W^-(x,\Phi)$. Then
\[
(y,x)=\lim_{t\ra\infty} (y,\Phi_t(y)) \in
\Gamma_{-\infty}\Longrightarrow (x,y)\in
^\ast\!\Gamma^{-\infty}=\Gamma^\infty.
\]
Hence $\{x\}\times W^-(x,\Phi)\subset \Gamma^\infty$.     The
inclusion $W^+(x,\Phi)\times \{x\}\subset  \Gamma^\infty$ is
proved in a similar fashion.

From the equality $\Gamma^\infty\cup
\Gamma^{-\infty}=\bar{\gamma}\setminus \Gamma$ we deduce
\[
\dim \Gamma_{\pm\infty}\leq \dim\bar{\Gamma}-1=\dim X.
\]
On the other hand,
\[
\dim \Gamma^\infty \geq  \max_{x\in \Cr_\Phi} W^\pm(x,\Phi).
\]
If we observe that
\[
X\setminus \Cr_\Phi=\bigsqcup_{x\in\Cr_\Phi} W^+(x,\Phi)\setminus
\{x\}= \bigsqcup_{x\in\Cr_\Phi} W^-(x,\Phi)\setminus \{x\},
\]
we deduce from the scissor equivalence principle that
\[
\dim X=  \max_{x\in \Cr_\Phi} W^+(x,\Phi)= \max_{x\in \Cr_\Phi}
W^-(x,\Phi),
\]
which proves that $\dim \Gamma^\infty =\dim X$.  \qed

\begin{ex} Consider the flow  $\Phi$ on $S^1$ described in Example \ref{ex: homoclinic} consisting of one stationary point $N$ and a single homoclinic orbit. Then
\[
\Gamma_\infty = \{N\}\times S^1 \cup S^1\times \{N\}\subset
S^1\times S^1.
\]
Consider the  gradient flow on $S^1$    described in  Example
\ref{ex: height}. It consists of a  repeller $N$, an attractor $S$
and   two orbits going from $N$ to $S$. Then
\[
\Gamma_\infty= \{N\}\times S^1\cup S^1\times \{S\}.\proofend
\]
\end{ex}

%% file: tameflow4.tex
For any smooth manifold $M$, and any  differentiable  function
$f:M\ra \bR$, we denote by $\Cr_f\subset M$ the set of critical
points of $f$, and by $\Delta_f\subset \bR$, the
\emph{discriminant set} of $f$, i.e. the set of critical values of
$f$.  For   every positive integer $\lambda$  and every positive
real number $r$ we denote by $\bD^\lambda(r)$ the open  Euclidean
ball in $\bR^\lambda$  of  radius $r$ centered at the origin.
When  $r=1$ we write simply $\bD^\lambda$.

\begin{definition} Suppose $M$ is a compact, real analytic manifold of dimension $m$.

\noindent (a) A \emph{Morse pair} on $M$ is a pair $(\xi,f)$,
where $\xi$ is a $C^2$ vector field  on $M$ and $f: M\ra \bR$ is a
$C^2$,  Morse function on $M$ satisfying the following conditions.

\begin{itemize}

\item $\xi\cdot f <0$ on $M\setminus \Cr_f$.

\item  For every critical point $p$  of $f$ of index $\lambda$
there exist an open neighborhood $U_p$  of $p\in M$,  a
$C^2$-diffeomorphism,
\[
\Psi: U_p\ra \bD^m,
\]
and real numbers $\mu_1,\dotsc, \mu_m >0$ such that $\Psi(p)=0$,
and
\[
\Psi_*(\xi)= \sum_{i\leq \lambda }
\mu_iu^i\pa_{u^i}-\sum_{j>\lambda }\mu_ju^j\pa_{v^j},
\]
where $(u^i)$ denote the Euclidean coordinates on $\bD^m$.
\end{itemize}

\noindent(b) The Morse pair $(f,\xi)$ is called \emph{tame} if the
function $f$  is tame, and the changes of coordinates $\Psi$ are
tame.

\noindent(c) We say that the  coordinate chart $(U_p,\Psi)$ is
\emph{adapted to $(\xi,f)$ at $p$},  and we define
\[
|u_-|^2 =\sum_{j\leq
\lambda}|u^j|^2,\;\;\;|u_+|^2=\sum_{j>\lambda}|u^j|^2,\;\;\mu:=
2\max (\mu_i)+1
\]
\[
E(u):=\left\{
\begin{array}{ccl}
0 & {\rm if} & \mbox{$p$ is a local max or local min}\\
&&\\
\underbrace{\Bigl(\,\sum_{i\leq
\lambda}(u^i)^{\mu/\mu_i}\,\Bigr)}_{E^-(u)}\,\cdot \,
\underbrace{\Bigl(\,\sum_{j>\lambda}(u^j)^{\mu/\mu_j}\,\Bigr)}_{E^+(u)}
& {\rm if} & \mbox{$p$ is a saddle point, $0<\lambda <m$}.
\end{array}
\right.
\]
(d) For  every triplet of  real numbers $\ve,\delta, r>0$ we
define the \emph{block}
\[
\eB_p(\ve,\delta,r):= \bigl\{\, u\in U_p;\;\; |f(u)-f(p)|<\ve,\;\;
E(u)<\delta,\;\;|u_-|^2+|u_+|^2 <r^2\bigr\}.
\]
(e) A \emph{Morse flow} on a compact real analytic manifold $M$ is
the  flow generated   by  a $C^2$-vector field $\xi$,   where  for
some  $C^2$-function $f: M\ra \bR$ the pair $(f,\xi)$ is a  Morse
pair on $M$. \qed
\end{definition}

In the above definition note that $\mu/\mu_i>2$ for any $i$, and
thus, if the pair $(f,\xi)$ is tame, the  function $E(u)$ is a
tame $C^2$-function. Observe also that in $U_p$ we have
\[
\xi\cdot E^\pm = \pm \mu E^\pm,
\]
so that $\xi\cdot E=0$, i.e., the quantity $E$ is conserved along
the trajectories of $\xi$. In the sequel we will need the
following technical result.

\begin{proposition} Suppose $(\xi,f)$ is a   Morse pair, and $p\in \Cr_f$. Then there exists $r_0= r_0(f)>0$ such that for every $r >0$ there exist $\ve_r,\delta_r >0$ such that
\[
\overline{\eB_p(\ve,\delta, r_0)}\subset
\eB_p(\ve,\delta,r)\subset \bD^m(r),\;\;\forall 0<\ve
<\ve_r,\;\;0<\delta<\delta_r.
\]
In other words, no mater how small $r$ is we can choose
$\ve,\delta>0$ sufficiently small so that  the isolating block
$\eB_p(\ve,\delta, r_0)$, a priori  contained in $\bD^m(r_0)$ is
in fact contained in a much smaller ball $\bD^m(r)$. \label{prop:
apriori}
\end{proposition}

\proof  Assume $f(p)=0$.  The statement is  obviously true if $p$
is a local min or a local max. We assume $p$ is a saddle point and
we denote  by
\[
H: T_pM\times T_pM\ra \bR
\]
 the Hessian of $f$ at $p$. There exist $C=C(f)>0$,  and $\alpha=\alpha(\mu_1,\dotsc,\mu_m)>1$ such that,
\[
|f(u)-\frac{1}{2} H(u)| \leq  C(|u_-|^3+ |u_+|^3),\;\;E_\pm(u)
\geq |u_\pm|^{2\alpha},\;\;\forall |u_+|+|u_-|\leq 1.
\]
We deduce that if $u\in \overline{\eB_p(\ve, \delta, r)}$   and
$r<1$ then
\[
-\ve + C(|u_-|^3+ |u_+|^3)\leq  \frac{1}{2} H(u) \leq \ve +
C(|u_-|^3+ |u_+|^3),
\]
\[
 |u_\pm|\leq  r,\;\;  |u_-|\cdot |u_+| \leq  \delta^{1/\alpha}.
\]
Now observe that there exist constants $0<a <1 <b$ such that
\[
-\frac{1}{a} |u_-|^2 + \frac{1}{b} |u_+|^2 \leq \frac{1}{2}
H(u)\leq -a|u_-|^2 +b|u_+|^2.
\]
Putting all of the above together we deduce that there exists
$C_1=C_1(f)>1$ such that  if
\[
 u\in \eB_p(\ve, \delta, r)\;\;\mbox{and}  \;\; r<1
\]
then,
\[
\frac{1}{C_1}(-\ve +|u_-|^3+ |u_+|^3)\leq   -|u_-|^2+ |u_+|^2 \leq
C_1(\ve +|u_-|^3+ |u_+|^3),
\]
and
\[
 |u_-|\cdot |u_+|\leq \delta^{1/\alpha}.
\]
Now fix $r_0=\frac{1}{2C_1} <1$. We want to show that for every $r
<r_0$ there exist $\ve,\delta
>0$ such that  $\eB_p(\ve,\delta, r_0)\subset \bD^m(r)$.

We argue by contradiction, and we assume there exists
$0<\bar{r}<r_0$  such that, for any $\ve,\delta>0$, we have
\[
\eB_p(\ve,\delta, r_0)\not\subset \bD^m(\bar{r}).
\]
We deduce that we can find    a sequence $u_n\in \eB_p(1/n,1/n,
r_0)$ such that $|u_n|\geq \bar{r}$.  Set
\[
s_n:= |(u_n)_-|,\;\;t_n:=|(u_n)_+|.
\]
We deduce that $s_n^2+r_n^2\geq \bar{r}^2$, and
\[
\frac{1}{C_1}\bigl(\,-\frac{1}{n} +s_n^3+t_n^3\,\bigr) <
-s_n^2+t_n^2\leq  C_1\bigl(\, \frac{1}{n} +s_n^3+t_n^3\,\bigr)
,\;\; s_nt_n \leq n^{-1/\alpha},\;\;s_n,\;t_n<r_0.
\]
The condition
\[
s_nt_n \leq n^{-1/\alpha},\;\;0\leq s_n,\;t_n<r_0=\frac{1}{2C_1}
\]
implies that a subsequence of $s_n$  converges to
$s_\infty\in[0,r_0]$, a subsequence of $t_n$ converges to
$t_\infty\in[0,r_0]$ and
\[
 s_\infty t_\infty=0,\;\; s_\infty^2+t_\infty^2\geq \bar{r}^2.
 \]
 We observe that $t_\infty\neq 0$ because, if that were the case,  we would  have
 \[
 0<\frac{1}{C_1}s_\infty^3 \leq -s_\infty^2 <0.
 \]
 Hence we must have  $s_\infty=0$ and $t_\infty \neq 0$. We deduce
 \[
 t_\infty^2 \leq C_1t_\infty^3,\;\; t_\infty >\bar{r}\Longrightarrow  \frac{1}{C_1}\leq t_\infty \leq r_0 =\frac{1}{2C_1}.
 \]
We have reached a contradiction. This concludes the proof of
Proposition \ref{prop: apriori}.  \qed

The above proposition implies that for $r<r_0$, and   any
$\ve,\delta$ sufficiently small,  we have
\[
\pa \eB_p(\ve,\delta, r_0)\cap \pa \bD^m(r) =\emptyset.
\]
When this happens we  say that  $\eB_p(\ve,\delta, r_0)$  is an
\emph{isolating block} of $p$.  The boundary of such an isolating
block  has a decomposition
\[
\pa \eB_p(\ve,\delta, r_0)=\pa_+ \eB_p(\ve,\delta, r_0)\cup \pa_-
\eB_p(\ve,\delta,r_0) \cup \pa_0 \eB_p(\ve,\delta, r_0),
\]
where
\[
\pa_\pm \eB_p(\ve,\delta, r_0 ) =\overline{\eB_p(\ve,\delta_,
r_0)}\cap \{ f=f(p)\pm\ve\},
\]
and
\[
\pa_0 \eB_p(\ve,\delta, r_0)= \overline{\eB_p(\ve,\delta, r_0)}
\cap\{ E(u)=\delta\}.
\]
The function $E(u)$ is  twice differentiable (since $\mu/\mu_i
>2$),  and it  is constant along the trajectories of $\xi$ while
$f$ decreases along these trajectories.    This implies that  no
trajectory of $\xi$ which starts at a point
\[
q\in \{    f(p)-\ve < f < f(p)+\ve\}\setminus \eB_p(\ve,\delta,
r_0)
\]
can intersect the block $\eB_p(\ve,\delta, r_0)$.

\begin{theorem}  Suppose $(\xi,f)$ is a tame  Morse pair on $M$ such that $\xi$ is real analytic. Then the  flow generated by $\xi$ is tame.
\label{th: morse-tame}
\end{theorem}

\proof  First, let us introduce some terminology.    Suppose
$(\xi, f)$ is a tame Morse pair on the real analytic manifold and
$\Phi:\bR\times M\ra M$ is the flow generated by $\Phi$. For every
set   subset $A\in M$ we set
\[
A^\xi:=\Bigl\{ y\in M;\;\; \exists  t\geq 0,\;\;x\in
A:\;\;y=\Phi_t(x)\,\Bigr\}.
\]
In other words, $A^\xi$ is the forward drift of $A$, i.e. the
region of $M$ covered by the forward  trajectories of $\xi$ which
start at a point in $A$. We define similarly
\[
 A^{-\xi}:=\Bigl\{ y\in M;\;\; \exists  t\leq 0,\;\;x\in A:\;\;y=\Phi_t(x)\,\Bigr\}.
 \]
\noindent {\bf Step 1.} Let $c\in \Delta_f$. We will show that
there exists $\si=\si(c)  >0$ such that, for  any $\ve \in
(0,\si)$, and any tame set $A\subset \{c-\ve < f < c+\ve\}$, the
intersection
\[
 A^\xi(\ve):= A^\xi\cap \{c-\ve < f < c+\ve\}
 \]
 is a tame set.

Let $\gamma >0$ such that   the only critical  value  of $f$  in
the interval $(c-\gamma,c+\gamma)$ is $c$, and set
\[
\Cr_f^c:=\Cr_f\cap \{f=c\}.
\]
$\Cr_f^c$ is a finite set.    We can find $\ve_0, r_0>0$  such
that,  for  every $\ve<\ve_0$, and  every $p\in \Cr_f$,  the
blocks $\eB_p:=\eB_p(\ve,\ve, r_0)$  are isolating, and   their
closures are pairwise disjoint.  Set
\[
\si=\min(\gamma,\ve_0).
\]
For $0<\ve <\si$,  $p\in \Cr_f^c$, and any  tame  subset
\[
A\subset\{ c-\ve < f <c+\ve \}
\]
we set,
\[
A_p:=A\cap \eB_p(\ve,\ve,r_0),\;\; \eB_\ve:=\bigcup_{p\in\Cr_f^c }
\eB_p(\ve,\ve,r_0),
\]
\[
Z_\ve= \{ c-\ve <f <c+\ve\}\setminus \eB_\ve, \;\;A_*=A\cap Z_\ve.
\]
Since
\[
A =A_*\cup \Bigl(\bigcup_{p \in\Cr_f^c}A_p)
\]
it suffices to show that  each of the subsets $A_*(\ve)^\xi$ and
$A_p(\ve)^\xi$ is  definable.

Note first that,  since the isolating blocks $\eB_p$ are definable
sets,   each $A_p$ is definable.

 For $p\in \Cr_f^c$   we denote by $\lambda_p$ its index, and we choose a   coordinate chart $(U_p,\Psi_p)$ adapted to $(\xi,f)$ near $p$ such that
 \[
 \Phi_t(u)= (e^{\mu_1t}u^1,\cdots, e^{\mu_\lambda t} u^\lambda, e^{-\mu_{\lambda+1} t} u^{\lambda+1}, \cdots , e^{-\mu_m t} u^m).
 \]
We  deduce  that
\[
A_p^\xi\cap \{ c-\ve <f <c+\ve \}
\]
\[
 =\bigl\{\,  u\in A_p;\;\; \exists t\geq 0:\;\; (e^{\mu_1t}u^1,\cdots, e^{\mu_\lambda t} u^\lambda, e^{-\mu_{\lambda+1} t} u^{\lambda+1}, \cdots , e^{-\mu_m t} u^m)\in \eB_p,\,\bigr\}
\]
This shows $A_p^\xi$ is definable.

Note that  no trajectory of $\xi$ starting  on $Z_\ve $ will
intersect  the  neighborhood $\eB_\ve$ of $\Cr_f^c$.   Let
\[
m:= \inf\Bigl\{\, |\xi\cdot f(x)|;\;\; x\in Z_\ve\,\Bigr\}.
\]
Observe that $m>0$. Fix $T> \frac{2\ve}{m}$.  We deduce that
\[
\forall x\in Z_\ve,\;\;\Phi_T(x)\in \{f<c-\ve\}.
\]
Since the vector  field  $\xi$ is real analytic we deduce from the
Cauchy-Kowalewski theorem (in the general form  proved in
\cite[I.\S 7]{CH}) that the flow map
\[
\Phi:[0,T]\times M\ra M
\]
is real analytic.

Observe that
\[
A_*(\ve)^\xi=\Bigl\{ y\in  M;\;\;  |f(y)-c| < \ve,\;\;  \exists
t\in [0,T],\;\;\exists x\in A_*: \Phi_t(x)=y\,\Bigr\}.
\]
This shows A$_*^\xi(\ve)$ is definable.  In particular, we deduce
that for $\ve <\si(c)$, and every  definable
\[
A\subset \{ c-\ve \leq f\leq c+\ve\}
\]
the set $A^\xi\cap \{ c-\ve \leq f\leq c+\ve\}$ is also definable.

\medskip

\noindent {\bf Step 2.} Suppose the interval $[a,b]$ contains no
critical values of $f$. Then for every tame set $A\subset \{a\leq
f\leq b\}$ the set
\[
A^\xi\cap \{ a\leq f\leq b\},\;\;{\rm and} \;\; A^{-\xi} \cap  \{
a\leq f\leq b\}
\]
are also tame.   Indeed, let
\[
m:=\inf\bigl\{ |\xi\cdot f(x)|;\;\; a\leq f(x)\leq b,\,\bigr\}.
\]
 Since the interval $[a,b]$ contains no critical values  we deduce that $m>0$. Fix $T>\frac{b-a}{m}$. Then
 \[
 \forall x\in \{f=b\},\;\;\Phi_T(x)\in\{f<a\}.
 \]
 Observe that
 \[
 A^\xi\cap f^{-1}([a,b])=\bigl\{ y\in M;\;\; a\leq f(y)\leq b,\;\; \exists t\in[0,T],\;\;x\in A:\;\;y=\Phi_T(x)\,\bigr\}.
 \]
 We deduce from the above description that $A^\xi\cap f^{-1}([a,b])$ is definable since $A$ is so and the map $\Phi:[0,T]\times M\ra M$ is real analytic.

\medskip

\noindent {\bf Step 3.}   Suppose $A$ is a tame subset  of $N$.
Then  $A^\xi$ and $A^{-\xi}$ are also tame. To prove this we must
first consider an $f$-\emph{slicing}.  This is a finite collection
of real numbers
\[
a_0<a_1<\cdots <a_n
\]
with the following properties.

\begin{itemize}

\item $f(M)\subset [a_0,a_n]$.

\item  $a_i$ is a regular value of $f$, $\forall i=0,\cdots, n$.

\item  Every interval $[a_{i-1},a_i ]$, $0\leq i\leq n$ contains
at most   one critical value of $f$.

\item  If the interval $[a_{i-1},a_i ]$ contains     one critical
value of $f$ then this critical value must be the midpoint
\[
c_i=\frac{a_i+a_{i-1}}{2}
\]
Moreover, the interval $[a_i,a_{i-1}]$ is  very short, i.e.,
$(a_i-a_{i-1})<\si(c_i)$.
\end{itemize}

 Fix an $f$-slicing $a_0< \cdots <a_n$, and a tame set $A\subset M$. Now define
 \[
 M_i:=f^{-1}([a_{i-1}, a_{i}]), \;\;\pa_- M_i= \{f=a_{i-1}\}.
 \]
 Then
 \[
 A^\xi= \bigcup_i (A\cap M_i)^\xi.
 \]
 We will prove by induction over $i$ that for every tame set $B\subset M_i$ the set $B^\xi$ is also tame.  For $i=1$, the interval $[a_0,a_1]$ must contain a  critical value, the absolute minimum and we conclude using  Step 1   since
 \[
 B\subset M_1\Longrightarrow B^\xi \subset M_1 \Longrightarrow B^\xi=B^\xi\cap M_1.
 \]
  Consider now a tame set $B\subset M_{i+1}$.  Then
  \[
  B^\xi = B^\xi \cap M_{i+1} \cup ( B^\xi\cap \pa_- M_{i+1})^\xi.
  \]
  $B^\xi \cap M_{i+1} $ is  tame by Step 1, if the interval $[a_i,a_{i+1}]$ contains a critical value, or   by Step 2, if the interval $[a_i,a_{i+1}]$ contains no critical value.

  Now observe that
  \[
B^\xi\cap \pa_- M_{i+1}= (B^\xi\cap M_{i+1})\cap \pa_- M_{i+1}
\]
so that $B^\xi\cap \pa_- M_{i+1}$ is a tame subset of
$\pa_-M_{i+1} \subset M_i$.  The induction hypothesis now implies
that $( B^\xi\cap \pa_- M_{i+1})^\xi$ is tame.

\medskip

\noindent {\bf Step 4.} \emph{Conclusion.} Suppose $(\xi,f)$ is a
tame Morse pair on $M$.    We construct a new tame Morse pair
$(\hat{\xi},\hat{f})$ on $S^1\times M\times M$ defined by
\[
\hat{f}(\theta, x,y)= h_0(\theta)-f(x)+f(y),\;\;\forall (p,
x,y)\in S^1\times M\times M,
\]
where $h_0:S^1\ra \bR$ is the height function  we  considered in
Example \ref{ex: height}.     Similarly
\[
\hat{\xi}(\theta,x,y)=\xi_0(\theta)\oplus -\xi(x)\oplus \xi(y),
\]
where $\xi_0$ is the  gradient of $-h_0$. Denote by $\theta_0$ the
point $(1,0)$ on the unit circle and let
\[
\Delta= \{ (\theta_0,  x,x)\in S^1\times X\times X\}.
\]
By Step 3 the set
\[
G=  \Delta^{\hat{\xi}}\cup \Delta^{-\hat{\xi}}= \Bigl\{  (\theta,
u,v)\in S^1\times M\times M;\;\; \exists t\in \bR,\;\;x\in
M:\;\;(\theta,u,v)=\Phi^{\hat{\xi}}_t(\theta_0,x,x)\,\Bigr\}
\]
is tame.  Since the  negative  gradient flow of $h_0$ in the open
half-circle $S=\{ x^2+y^2=1;\;\;x>0\}$   is tamely conjugate to
the  translation flow  on $\bR$ we  deduce from  Corollary
\ref{cor: tame} that $\Phi^\xi$ is a tame flow.   \qed

\begin{theorem}  Suppose $X$ is a compact, real analytic manifold and $f: X\ra \bR$ is a real analytic Morse function.  Then  for every real analytic metric $g_0$  on $X$ and every $\ve >0$ there exist a real analytic metric $g$ on $X$ with the following properties.

\smallskip

\noindent $\bullet$  $\|g_0-g\|_{C^2}\leq \ve$.

\noindent $\bullet$  $(f, -\nabla^g f)$  is a tame Morse pair.

\smallskip

In particular, the flow generated by $-\nabla^gf$ is a tame Morse
flow. \label{th: PS}
\end{theorem}

\proof  The proof is based on a simple strategy. We  show that
we can find  real analytic metrics $g$  arbitrarily $C^2$-close to
$g_0$ such that  the gradient   vector field $\nabla^g f_0$ can be
linearized by a real analytic change of coordinates localized in a
neighborhood  of  the critical set. The linearizing  change of
coordinates  is   obtained  by invoking the Poincar\'e-Siegel
theorem \cite[Chap. 5]{Ar2} on the normal forms  of real analytic
vector fields in a neighborhood of an isolated  stationary point.

We digress to  recall the Poincar\'{e}-Siegel theorem. Suppose
$\vec{Z} $ is a real analytic vector field defined in a
neighborhood $\eN$ of the origin  $0$ of the Euclidean vector
space $\bR^n$.  Assume that $0$ is an isolated    stationary point
of the  vector field  $\vec{Z}$. If we regard $\vec{Z}$ as a real
analytic map
 \[
 \vec{Z}: \eN\ra \bR^N
 \]
 then we obtain a Taylor expansion near $0$
 \[
 \vec{Z}(x)=    L\cdot x +  \mbox{higher order terms},\;\;x\in\eN,
 \]
 where $L: \bR^n\ra \bR^n$ is a  linear operator. We regard it as a linear operator $T_0\bR^n\ra T_0\bR^n$. As such,  the operator $L$ is independent of the choice of coordinates.   More precisely, for every   linear connection  $\nabla$ on the tangent bundle $T\bR^n$ we have
 \[
 L v= (\nabla_v \vec{Z})_0,\;\;\forall v\in T_0\bR^n.
 \]
The Poincar\'e-Siegel theorem describes    conditions on $L$ which
imply the existence of  real analytic coordinates $y=(y^1,\cdots,
y^n)$  near $0\in \bR^n$ such that in these new coordinates the
vector field $\vec{Z}$ is linear,
\[
\vec{Z}(y)= Ly= \sum_j y^j L\pa_{y^j}.
\]
We describe these conditions only in the case when  $L$ is
semisimple (diagonalizable) and all its eigenvalues  are real
since this is the only case of interest to us.

Denote the eigenvalues of $L$ by
\[
\mu_1\leq \mu_2\leq \cdots\leq \mu_n.
\]
We write
\[
\vec{\mu}= (\mu_1,\dotsc, \mu_n)\in\bR^n.
\]
We say that  $L$  satisfies the \emph{Siegel $(C,\nu)$}-condition
if, for any $k=1,\dotsc, n$, and any  $\vec{m}=(m_1,\cdots,
m_n)\in (\bZ_{\geq 0})^n$ such that
\[
|\vec{m}|:= m_1+\cdots+m_n\geq 2,
\]
we have
\[
|\mu_k - (\vec{m},\vec{\mu})| \geq \frac{C}{|\vec{m}|^\nu}.
\]
We denote by $\eS_{C,\nu}\subset \bR^n$  the set of vectors
$\vec{\mu}$ satisfying the  Siegel $(C,\nu)$-condition,   and we
set
\[
\eS_\nu:=\bigcup_{C>0}\eS_{C,\nu}.
\]
Then the  set   $\bR^n\setminus \eS_\nu$ has  zero Lebesgue
measure if $\nu >\frac{n-2}{2}$, \cite[\S 24.C]{Ar2}. In other
words, if we fix $\nu>\frac{n-2}{2}$ then  almost every vector
$\vec{\mu}\in \bR^n$ satisfies the Siegel $(C,\nu)$-condition  for
some $C>0$.  We can now state the Poincar\'{e}-Siegel theorem
whose   very delicate proof can be found in   \cite[Chap.5]{Ar2}
or \cite{Pliss}.

\begin{theorem}[Poincar\'{e}-Siegel] Suppose  that the eigenvalues $(\mu_1,\dotsc, \mu_n)$ satisfy the Siegel $(C,\nu)$ condition for some $C>0$ and $\nu >0$. Then there exist local, real analytic coordinates $y=(y^1,\dotsc, y^n)$ defined in a neighborhood of $0\in \bR^n$ such that, in these coordinates, the vector field $\vec{Z}$ is linear,
\[
\vec{Z}(y)= L(y).\proofend
\]
\end{theorem}

After this     digression  we return to our original problem.

According to \cite{GJ},  can find a real analytic   isometric
embedding of  $(X,g_0)$ in some Euclidean space $\bR^N$. For every
real analytic metric $g$ on $X$ we set
$\vec{Z}_g:=\nabla^{g}f_0\in\Vect(X)$.  For every  $p_0\in
\Cr_{f_0}$ we denote by $L_{g,p_0}: T_{p_0}X\ra T_{p_0}X$  the
linear operator defined by
 \[
 L_{g,p_0}v:=(\nabla_v  Z_g)_{p_0},\;\;\forall v\in T_{p_0}X,
 \]
 where $\nabla$ is  any connection on $TX$. Since   $Z_g(p_0)=0$, the operator $L_{g,p_0}$ is independent of the choice of the connection $\nabla$.

 The operator $L_{g,p_0}$ is symmetric (with respect to the metric $g$) and  thus  diagonalizable. More precisely, if we choose local analytic coordinates  on $X$ near $p_0\in\Cr_{f_0}$ such that $\pa_{x_i}$ form a $g$-orthonormal basis of $T_{p_0}X$  which diagonalizes the Hessian matrix $\bigl(\, \pa^2_{x^ix^j}f_0(p_0)\,\bigr)$ then
 \[
 Z_g=\sum_{i,j} g^{ij} (\pa_{x^j}f_0)\pa_{x^i}
 \]
 \[
 L_{g,p_0}\pa_{x^k}=\sum_{i,j} \Bigl(\,\bigl(\,\pa_{x^k}g^{ij}\,\bigr)(p_0)\cdot( \pa_{x^j}f_0)(p_0) + g^{ij}(p_0) \bigl(\,\pa^2_{x^kx^j}f_0\,\bigr)(p_0)\,\Bigr)\pa_{x^i}.
 \]
 \[
 = \sum_{i,j}\delta^{ij}\bigl(\,\pa^2_{x^kx^j}f_0\,\bigr)(p_0)\pa_{x^i}= (\pa^2_{x^kx^k}f_0)(p_0)\pa_{x^k}
 \]
We want to prove  that arbitrarily close  to any real analytic
metric $g$ we can find real analytic metrics $h$   such that for
every $p_0\in \Cr_{f_0}$ there exist  real analytic  coordinates
$y$ so that in these coordinates   the vector field $Z_{h,p_0}$
has the linear form
  \[
  Z_{h,p_0}(y)= L_{h,p_0}(y).
  \]
Here is  the strategy. Near each $p_0$ choose local analytic
coordinates $(x^i=x_{p_0}^i)$  as above,  meaning that
$(\pa_{x^i})$ is a $g$-orthonormal basis of $T_{p_0}X$ which
diagonalizes the Hessian matrix, i.e.
\[
\pa^2_{x^ix^j}f_0(p_0)=0,\;\;\forall i\neq j.
\]
If $h$ is another metric on $X$ then the above computation shows
that
\[
L_{h,p_0} \pa_{x^k}=
\sum_{i,j}h^{ij}\bigl(\,\pa^2_{x^kx^j}f_0\,\bigr)(p_0)\pa_{x^i}=\sum_i
h^{ik}\bigl(\,\pa^2_{x^kx^k}f_0\,\bigr)(p_0)\pa_{x^i}.
\]
Denote by $\Sym^+(n)$ the space of positive definite, symmetric
$n\times n$ matrices.  We will show that  for any  map
\[
A:\Cr_{f_0}\ra \Sym^+(n),\;\; p\mapsto A_p
\]
close to the identity  map
\[
\bsI: \Cr_{f_0}\ra \Sym^+(n),\;\;p\mapsto \bsI_n,
\]
there exists a  real analytic metric  $h$, close to $g$, such
that, for every $p\in \Cr_{f_0}$,  the matrix   describing  $h$ at
$p$  in the coordinates $(x^i_p)$ chosen above is equal to
$A_p^{-1}$.  In other words, we want to show that as $h$ runs
through a small neighborhood of $g$,   the collection of matrices
\[
\Cr_{f_0}\ni p\mapsto  \bigl(\,h^{ij}(p)\,\bigr)_{1\leq i,j\leq
n}\in\Sym^+(n)
\]
 spans  a small neighborhood of the identity map.    This is achieved via a genericity result.      We can then prescribe  $h$ so that   at every $p\in\Cr_{f_0}$ the linearization $L_{h,p}$ satisfies  the conditions of the  Poincar\'{e}-Siegel theorem.

To prove that we can prescribe the metric $h$ any way we please at
the points in $\Cr_{f_0}$  we will prove an elementary genericity
result.   To state it we need a bit of terminology.

Fix a finite dimensional Euclidean space $E$ and denote by
$\eP_d(\bR^N,E)$ the  vector space of polynomial maps $\bR^N\ra E$
of degree  $\leq d$.    For every   $E$-valued,   real analytic
function $f$ defined in the neighborhood of a point $x\in\bR^N$,
and every nonnegative integer $k$ we denote by $j_k(f,x)\in
\eP_k(\bR^n, E)$ the $k$-th jet of $f$ at $x$.

\begin{lemma}  Let $B\subset \bR^N$ be an open ball and  $S\subset B$ a finite subset.        For every  integers $d>k>0$ define the  \emph{linear map}
\[
\ev_{d,k}: \eP_d(\bR^N,\bR^\ell)\ra  \prod_{s\in S}
\eP_k(\bR^N,\bR^\ell),\;\;   f\longmapsto   ( j_k(f,s)\,)_{s\in
S}.
\]
Then   $\ev_{d,k}$  is  onto if $d\geq 2(k+1)|S|-2$. \label{lemma:
poly-one}
\end{lemma}

\proof\footnote{The idea of this proof arose in conversations with
my colleague  R. Hind.}    It suffices to show that for every
$s_0\in S$, and every $P_0\in \eP_k(\bR^N,\bR^\ell)$ there exists
$f\in \eP_d(\bR^n, \bR^\ell)$ such that
\[
j_k(f,s_0)= j_k(P_0,s_0)\;\;j_k(f,s)=0,\;\;\forall s\neq s_0.
\]
Clearly it suffices to prove this only in the case $\ell=1$. For
every $s\in S$ define
\[
\rho_s(x)=|x-s|^2\in\eP_2(\bR^N,\bR).
\]
Observe that
\[
j_k(\rho_s^{k+1}, s)=0,\;\;\;\forall s\in S,\;\;\forall k\geq 0.
\]
Now define
\[
Q_{s_0}=\prod_{s\neq s_0} \rho_s^{k+1} ,\;\;\deg Q_{s_0}=
2(k+1)(|S|-1).
\]
Observe that for every  polynomial function $p$ we have
\[
j_k(pQ_{s_0},s)=0,\;\;\forall s\neq  s_0.
\]
The function  $1/Q_{s_0}$ is real analytic in a neighborhood of
$s_0$, and we denote by $R_{s_0}\in \eP_k(\bR^N,\bR)$ the $k$-th
jet of $1/Q_{s_0}$ at $s_0$. Then
\[
j_k(R_{s_0} Q_{s_0}, s_0)=1.
\]
Now define
\[
f= P_0 R_{s_0} Q_{s_0},\;\;\deg f= k + k +
(2k+2)(|S|-1)=2(k+1)|S|-2.
\]
Then
\[
j_k(f,s_0)=  j_k\bigl(\, j_k(P_0,s_0)\cdot j_k(R_{s_0} Q_{s_0},
s_0), s_0\,\bigr)= j_k(P_0,s_0),
\]
and
\[
j_k(f,s)=0,\;\;\forall s\neq s_0. \proofend
\]

Suppose now that the set $S$ lies on the  compact real analytic
submanifold $X\subset \bR^N$.   By choosing   real analytic
coordinates  on $X$ near each point $s\in S$ we obtain locally
defined  real analytic embeddings
\[
 i_s: U_s\subset \bR^n\ra \bR^N,\;\;i_s(0)=s,\;\;i_s(U_s)\subset X.
 \]
 Here, for every $s\in S$, we denoted by  $U_s$  a small,  open  ball centered at  $0\in \bR^n$, $n=\dim X$.    In particular, for every Euclidean vector space $E$, and every  positive integer $k$ we  obtain  \emph{surjective} linear maps
 \[
 \pi_s: \eP_k(\bR^N, E)\ra \eP_k(\bR^n,E),\;\; \eP_k(\bR^N, E) \ni f \mapsto j_k(\, f\circ i_s, 0).
 \]
For $x\in X$ we denote by $J_k(X,x,E)$ the space of $k$-jets at
$x$ of $E$-valued  real analytic maps  defined in a neighborhood
of $x$.  If $f$ is such a map, then we denote by $j_k(f,x)\in
J_k(X,x,E)$ its $k$-th jet.  We topologize $\eP_d(\bR^N,E)$ by
setting
\[
 |f| =\|f\|_{C^2(X)}
 \]
 Lemma  \ref{lemma: poly-one} implies the following result.

\begin{lemma}  Suppose $S$ is a finite subset  of $X$.  Then for any finite dimensional  Euclidean space $E$, and any integer     $d \geq 2|S|-2$ the linear map
\[
\ev:\eP_d(\bR^N,E) \ra \prod_{s\in S} J_0(X,s,E),\;\;
f\longmapsto  \bigl(\, j_0(f, s)\,\bigr)_{s\in S}
\]
is onto.  In particular, for every  $\ve >0$, the image of an
$\ve$-neighborhood of $0\in \eP_d(\bR^N,\bR)$  is an open
neighborhood of $0$ in $\prod_{s\in S} J_0(X,s,E)$.\qed
\label{lemma: poly-two}
\end{lemma}

We now specialize $E$ to  be the space $\Sym(N)$  of symmetric
bilinear forms $\bR^N\times \bR^N\ra \bR$, and thus      the space
of functions
\[
\bR^N\ra \Sym(N)
\]
can be viewed as the  space of deformations  of   Riemann metrics
on $\bR^N$.   The metric $g_0$ on $X$ is induced from the
Euclidean metric $\delta$ on $\bR^N$. If we deform   $\delta$
\[
\delta \ra \delta +   h,\;\; h\in \eP_d(\bR^N, E),\;\;d>2|S|,
\]
and $|h|$ is sufficiently small, then $\delta+h$ will still be a
metric on a neighborhood  of $X$ in $\bR^n$.

Fix $s\in X$.  Choose affine Euclidean coordinates
$(y^\alpha)_{1\leq \alpha\leq N}$ on $\bR^N$  such that
$y^\alpha(s)=0$, $\forall \alpha$.  Choose   local  real analytic
coordinates  $(x^1,\cdots, x^n)$ on $X$ in a neighborhood $U_s$ of
$s$ such that $x^j(s)=0$, $\forall j$. Along $X$ near $s$  the
vector   field $\pa_{x^i}$ is   represented by the vector field
\[
\sum_\alpha\frac{\pa y^\alpha}{\pa x^i}\pa_{y^\alpha}.
\]
If
\[
g(\pa_{y^\alpha},\pa_{y^\beta})=(\delta +
h)(\pa_{y^\alpha},\pa_{y^\beta})= \delta_{\alpha\beta}+
h_{\alpha\beta}
\]
then
\[
g(\pa_{x^i},\pa_{x^j})=\sum_{\alpha,\beta}
\bigl(\,\delta_{\alpha\beta}+ h_{\alpha\beta}\,\bigr)\frac{\pa
y^\alpha}{\pa x^i}\frac{\pa y^\beta}{\pa
x^j}=g_0(\pa_{x^i},\pa_{x^j})+\sum_{\alpha,\beta}
h_{\alpha\beta}\frac{\pa y^\alpha}{\pa x^i}\frac{\pa y^\beta}{\pa
x^j}.
\]
We think of $g_{ij}(x)$  as a real analytic map  from $U_s$ to the
space of symmetric  $n\times n$ matrices,    of $\frac{\pa
y^\alpha}{\pa x^i}$ a  real analytic map $ Y$ from $U_s$ to the
space of $N\times n$ matrices, and we think of $h$ as a real
analytic map from $U_s$ to the space of symmetric $N\times N$
matrices. Then
\[
g= g_0+ Y^t h|_{U_s} Y.
\]
Along $U_s$ we write
\[
Y= Y(0)+ O(1),\;\;  h= h(0)+ O(1)
\]
so that
\[
g = g_0 +  Y(0)^th(0)Y(0)+  O(1).
\]
$Y(0): \bR^n\ra \bR^N$ is an injective map   as it describes the
canonical  injection  $T_sX\hra T_s\bR^N$.      The correspondence
\[\
\Sym(N)\ni h\longmapsto  \eY(h):= Y^t(0) h Y(0)\in \Sym(n)
\]
is a linear map.  Intrinsically,  $\eY$ is the restriction map,
i.e.
\[
\eY(h)\,(u,v)=h(u,v),\;\;\forall u,v\in T_s X\subset T_s\bR^n.
\]
This shows that $\eY$ is onto, because any symmetric bilinear map
on $T_sX$ can be extended (in many different ways) to a symmetric
bilinear map on $\bR^N$. This concludes the proof of Theorem
\ref{th: PS}. \qed

We can refine the above existence result some more.

\begin{theorem} Suppose $M$ is a compact, real analytic manifold, $\dim M=m$,  $f:M\ra \bR$ is a real analytic  tame Morse function. For every critical point  $p$ of $f$ we denote by $\lambda(p)$  the Morse index of  $f$ at $p$.     For every $p\in\Cr_f$, we choose a vector
\[
\vec{a}(p)=\bigl(\, a_1(p),\dotsc, a_m(p)\,\bigr)\in \bR^m
\]
such that
\[
a_1\leq \cdots \leq a_{\lambda(p)} <0 < a_{\lambda(p)+1} \leq
\cdots \leq a_m.
\]
Then, for every $\ve >0$, we can find a  real analytic metric  $g$
on $M$,   such that  for every critical point  $p$ of $f$ there
exist  real analytic coordinates $(x^i)$ near $p$, and a vector
$\vec{b}=\vec{b}(p)\in\bR^m$ with the following properties.

\smallskip

\noindent (a) $x^i(p)=0$, $\forall i=1,\dotsc, m$.

\noindent (b) $|\vec{b}(p)-\vec{a}(p)|<\ve$.

\noindent (c) In the  coordinates $(x^i)$ the vector field
$\nabla^gf$ is  described by,
\[
\nabla^g f=\sum_{i=1}^m b_ix^i\pa_{x^i}.
\]
\label{th: spec}
\end{theorem}

\proof  From the Morse lemma we deduce that  we can find a
\emph{smooth}  metric $g_0$  on $M$ with the property  that for
every critical point $p$  there exist \emph{smooth} coordinates
$(y^i)$  near $p$ with the property that
\[
y^i(p)=0,\;\;\forall i=1,2,\dotsc, m
\]
and
\[
\nabla^{g_0} f =\sum_{i=1}^m a_i y^i\pa_{y^i}.
\]
Now choose a  real analytic metric $g_1$, sufficiently close to
$g_0$ such that the linearization  of $\nabla^{g_1} f$ at $p$ is
given by a diagonalizable operator $L_p: T_pM\ra T_pM$ with
eigenvalues
\[
\ell_1(p)\leq \cdots \leq \ell_m(p)
\]
satisfying
\[
|\vec{\ell}(p)-\vec{a}(p)|
<\frac{\ve}{2},\;\;\vec{\ell}(p)=(\ell_1(p),\dotsc, \ell_m(p)).
\]
Using Theorem \ref{th: PS} we can  find a real analytic metric
$g$ on $M$  such that the gradient vector field  $\nabla^gf$ can
be linearized  by an analytic    change of coordinates in a
neighborhood of every critical point, and  for  every critical
point $p$ the linearization of $\nabla^g f$ at $p$ is a
diagonalizable linear operator $B_p:T_pM\ra T_pM$  whose
eigenvalues
\[
b_1(p)\leq \cdots \leq b_m(p)
\]
satisfy
\[
|\vec{b}(p)-\vec{\ell}(p)| <\frac{\ve}{2}.
\]
This completes the proof of Theorem \ref{th: spec}. \qed

%% file: tameflow5.tex
Suppose $(f,\xi)$ is a tame Morse pair on the compact real
analytic manifold $M$, $\dim M=m$, such that the flow $\Phi^\xi$
generated by $\xi$ is tame.  Then,   for every critical point  $p$
of the Morse function $f$  we denote by $W^+(p,\xi)$ (respectively
$W^-(p,\xi)$)  the stable (respectively the unstable) variety of
$p$ with respect to the flow  $\Phi^\xi$.  Then $W^-(p,\xi)$  is a
$C^2$-submanifold of $M$ homeomorphic to $\bR^{\lambda(p)}$,
where $\lambda(p)$ denotes the Morse index of $f$ at $p$.
Similarly, $W^+(p,\xi)$ is a $C^2$-submanifold of $M$ homeomorphic
to $\bR^{m-\lambda(p)}$.

We say that $\Phi^\xi$ satisfies the \emph{Morse-Smale condition}
if, for every  pair of critical points $p,q$ such that
$f(p)>f(q)$,  the unstable manifold of $p$ intersects
transversally the stable manifold of $q$.

\begin{theorem} Suppose $M$ is a compact, real analytic manifold of dimension $m$, and $(f,\xi)$ is a tame Morse pair such that both $f$ and $\xi$ are real analytic.  Denote by $\Phi^\xi$ the flow generated by $\xi$. Then there exists a  smooth vector field $\eta$, which coincides with $\xi$ in an open neighborhood of the critical set of $f$, such that the pair $(f, \eta)$ is a tame Morse pair,  the flow generated by $\eta$ is tame and satisfies  the Morse-Smale condition.
\label{th: MS}
\end{theorem}

\proof  We follow closely the  approach pioneered by S. Smale (see
e.g. \cite[Section 2.4]{N2}). For simplicity, we assume $f$ is
nonresonant, i.e.,  every critical level set of $f$ contains a
unique critical point. Suppose the critical points are
\[
p_0,\dotsc, p_\nu,\;\; f(p_0)<f(p_1)<\cdots <f(p_\nu).
\]
For simplicity, we set $c_k:=f(p_k)$. Define
\[
\hbar =\min_{k=1,\dotsc, \nu} (c_k-c_{k-1}).
\]
 We will prove by induction  that for every  $k=0,1,\dotsc, \nu$, and for every $0 <\ve <\frac{\hbar}{100}$,  there exists  a tame $C^\infty$ vector field $\eta_k$ on $M$  with the following properties.

\begin{itemize}

\item $\eta_0=\xi$.

\item  $\eta_k(x)=\eta_j(x)$, $\forall 0\leq j<k$,  $\forall x\in
M$ such that $f(x)\not \in (c_k-2\ve, c_k-\ve)$.

\item  The pair $\bigl(f, \, \eta_k,\bigr)$  is a tame  Morse
pair.

\item The flow generated by $\eta_k$ is  tame and
\[
W^-(p_j, \eta_k)\pitchfork W^+(p_i, \eta_k),\;\;\forall 0\leq i<j
\leq k.
\]
\end{itemize}

\begin{figure}[h]
\centerline{\epsfig{figure=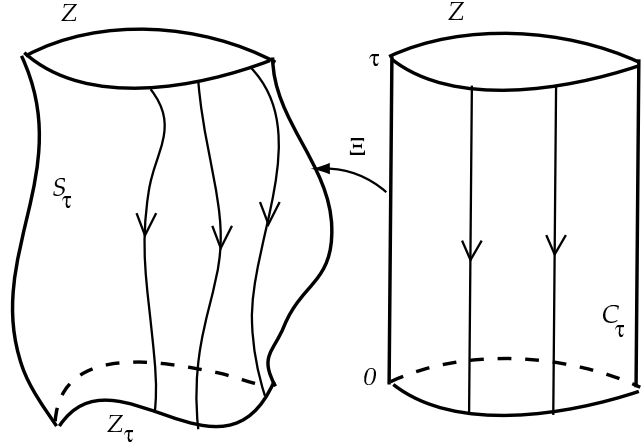,height=2.3in,width=3.3in}}
\caption{\sl Truncating a Morse flow. } \label{fig: 8}
\end{figure}

The statement is trivial for $k=0$ so we proceed directly to the
inductive step. Assume we have constructed $\eta_0,\dotsc,
\eta_k$, and we want to produce $\eta_{k+1}$.  Denote by
$\Phi^k_t$ the flow generated by $\eta_k$.  Set  $Z:=\{
f=c_{k+1}-\ve\}$.   Then there exists   $\tau>0$ such that
\[
\forall z\in Z,\;\;\forall t\in [0,t]:  f(\Phi^k_tz) > c_k-2\ve.
\]
Set (see Figure \ref{fig: 8})
\[
Z_\tau:=\Phi_\tau^k(Z_\ve),\;\;  S_\tau= S_\tau=
\bigcup_{t\in[0,\tau]}\Phi^k_t(Z),\;\;C_\tau=[0,\tau]\times Z,
\]

and define
\[
X= W^-(p_{k+1},\eta_k)\cap Z, \;\; Y= \bigcup_{j\leq k}
W^+(p_k,\eta_k)\cap Z.
\]
$Z$ is a real analytic  manifold, $X$ is a compact, real analytic
submanifold of  $Z_\ve$, while $Y$ is a smooth  submanifold of
$Z_\ve$. From to the classical transversality results of Whitney
(see \cite{Whit2} or \cite[Chap.3, 8]{Hir}) we deduce that there
exists a smooth  map
\[
h:[0,\tau]\times Z \ra Z,\;\; (t,z)\mapsto h_t(z)
\]
with the following properties.

\begin{itemize}

\item $h_0(z)= z$, $\forall z\in Z$.

\item $h_t$ is a diffeomorphism of $Z$, $\forall t\in [0,\tau]$.

\item $h_\tau(X)$ intersects  $Y$ transversally.

\end{itemize}

Using the approximation results in  \cite[Theorem 6]{Roy} we can
assume that $h$ is real analytic. Now choose a  smooth, increasing
tame function
\[
\alpha:[0,\tau]\ra [0,\tau]
\]
such that  $\alpha(t)=0$ for  all $t$ near zero, and
$\alpha(t)=\tau$,  for all $t$ near $\tau$. Define
\[
H: [0,\tau]\times Z\ra Z,\;\; H_t(z):= h_{\alpha(t)}(z).
\]
In other words, $H$ is a smooth, \emph{tame} isotopy between the
identity $\bsI_{Z}$ and $h_1$, which is independent of $t$ for $t$
near $0$ and $\tau$.

The tame  flow $\Phi^k$ defines a  smooth  tame diffeomorphism
(see Figure \ref{fig: 8}),
\[
\Xi: C_\tau=[0,\tau]\times Z \ra
S_\tau,\;\;\Xi_t(z)=\Phi^k_{t}(z).
\]
The diffeomorphism $\Xi$ maps $\eta_k$ to the vector field $\pa_t$
on $C_\tau$.

Using the isotopy $H$ we obtain a smooth tame diffeomorphism
\[
\hat{H}: C_\tau\ra C_\tau,\;\;\hat{H}(t,z)= H_t(z),
\]
such that
\[
\hat{H}_*\pa_t=\pa_t\;\; \mbox{near $\{0\}\times Z$ and
$\{1\}\times Z$}.
\]
The pushforward of the vector field $\eta_k|_{S_\tau}$ via the
diffeomorphism
\[
F=\Xi\circ\hat{H}\circ \Xi^{-1}: S_\tau\ra S_\tau
\]
is a smooth vector field   which coincides with $\eta_k$  in a
neighborhood of $Z$ and in a neighborhood of $Z_\tau$. Now define
the \emph{smooth} vector field $\eta_{k+1}$ on $M$ by
\[
\eta_{k+1}(x):=\begin{cases}
\eta_k(x) & x\in M\setminus S_\tau\\
F_*\eta_k(x) & x\in S_\tau.
\end{cases}
\]
$\eta_{k+1}$ is a smooth vector field,  and we denote by
$\Phi^{k+1}$ the flow on $M$ it generates.  Observe that
$\eta_{k+1}$ coincides with the original  vector field $\xi$ in an
open neighborhood of the critical set of $f$, and $f$ decreases
strictly on the nonconstant trajectories of $\eta_{k+1}$.  By
construction, we have
\[
W^-(p_j,\eta_{k+1})\pitchfork W^+(p_i,\eta_{k+1}),\;\;\forall
0\leq i < j\leq k+1.
\]
We want to prove that it is a tame flow. We will prove that the
maps
\[
\Phi^{k+1}:[0,\infty)\times M\ra M,\;\;
\Phi^{k+1}:(-\infty,0]\times M\ra M
\]
are  definable. We discuss only the first one, since the proof for
the second map is completely similar.

Observe first that,  $\hat{H}$ extends to a tame diffeomorphism
\[
\bR\times Z\ra \bR\times Z.
\]
We denote by $\Psi$  the tame  flow  $\bR\times Z$ obtained by
conjugating the translation flow with $\hat{H}$, i.e.,
\[
\Psi_t(s,z)= \hat{H}_{t+s} H_s^{-1} (z) .
\]
We divide  $M$ into three definable parts
\[
S_\tau,\;\; M_+:=\{ f> c_{k+1}-\ve\},\;\; M_-:= M\setminus
(S_\tau\cup M_+).
\]
We now have  definable functions
\[
T_+: M_+\ra (0,\infty],\;\; T_0,s:S_\tau\ra [0,\tau],
\]
\[
T_+(x):=\mbox{the moment of time when the trajectory of
$\Phi^{k+1}$ originating at $x$ intersects $Z$.}
\]
\[
T_0(x):=\mbox{the moment of time when the trajectory of
$\Phi^{k+1}$ originating at $x$ intersects $Z_\tau$,}
\]
and
\[
s(x)= \tau-T_0(x).
\]
We distinguish three cases.

\begin{itemize}

\item If $x\in M_-$  then $\Phi^{k+1}_t(x)=\Phi^k(x)$, $\forall
t\geq 0$.

\item If $x\in S_\tau$ then
\[
\Phi^{k+1}_t(x)=\begin{cases}
\Xi\circ \Psi_t\circ\Xi^{-1}(x) & t\leq T_0(x)\\
&\\
\Phi^k_{t-T_0(x)} \Xi\circ \Psi_{T_0(x)}\Xi^{-1}(x) & t>T_0(x)
\end{cases}
\]
\item If $x\in M_+$ then
\[
\Phi^{k+1}_t(x)=\begin{cases}
\Phi^k_t(x) &  t <T_-(x)\\
&\\
\Xi\circ \Psi_{t-T_-(x)}\circ\Xi^{-1}\circ \Phi^{k}_{T_-(x)}(x) & t\in (T_-(x), T_-(x)+\tau]\\
&\\
\Phi^k_{t-\tau-T_-(x)}\circ\Xi\circ \Psi_{\tau}\circ\Xi^{-1}\circ
\Phi^{k}_{T_-(x)}(x) & t>\tau+T_-(x).
\end{cases}
\]
\end{itemize}
This shows that $\Phi^{k+1}:[0,\infty)\times M\ra M$ is definable.
\qed

%% file: tameflow6.tex
In this section we  collect a few facts about the \emph{gap}
between two vector subspaces, \cite[IV.\S 2]{Kato}.

Suppose $E$ is a finite dimensional Euclidean space.  We denote by
$(\bullet,\bullet)$ the inner product on $E$, and by $|\bullet|$
the associated Euclidean norm. We define as usual the norm of a
linear operator $A:E\ra E$ by the equality
\[
\|A\|:=\sup\bigl\{ |Ax|;\;\;x\in E,\;\;|x|=1\,\bigr\}.
\]
The finite dimensional vector space $\End(E)$ is equipped with an
inner product
\[
\lan A,B\ran :=\tr(AB^*),
\]
and we set
\[
|A|:=\sqrt{\lan A,A\ran}=\sqrt{\tr(AA^*)}=\sqrt{\tr(A^*A)}.
\]
Since $E$ is finite dimensional, there exists a constant $C>1$,
depending only on the dimension of $E$, such that
\begin{equation}
\frac{1}{C}|A|\leq \|A\|\leq C|A|. \label{eq: norm-equiv}
\end{equation}
If $U$ and $V$ are two subspaces of $E$, then we define the
\emph{gap} between $U$ and $V$ to be the real number
\[
\de(U,V):= \sup\bigl\{\, \dist(u, V);\;\;u\in U, |u|=1\,\bigr\}
\]
\[
=\sup_u\inf_v \{\, |u-v|;\;\;u\in U, |u|=1,\;\;v\in V,\bigr\}.
\]

If  we denote by $P_{V^\perp}$ the orthogonal projection onto
$V^\perp$, then  we deduce
\begin{equation}
\de(U,V)=\sup_{|u|=1} |P_{V^\perp} u|=
\|P_{V^\perp}P_U\|=\|P_U-P_VP_U\|=\|P_U-P_UP_V\|. \label{eq:
gap-norm}
\end{equation}
Note that
\begin{equation}
\de(V^\perp, U^\perp)=\de(U,V). \label{eq: gap-dual}
\end{equation}
Indeed,
\[
\de(V^\perp, U^\perp)=\|P_{V^\perp}-P_{U^\perp}P_{V^\perp}\|=
\|\one -P_V -(\one-P_U)(\one-PV)\|
\]
\[
= \|P_U-P_UP_V\|=\de(U,V).
\]
We  deduce that
\[
0\leq \de(U,V)\leq   1,\;\;\forall U, V.
\]
Let us  point out that
\[
\de(U,V) < 1\Longleftrightarrow\;\;\dim U \leq  \dim V,\;\;U\cap
V^\perp =0.
\]
Note that  this implies that  the gap is \emph{asymmetric}  in its
variables, i.e. we cannot expect
\[
\de(U,V)=\de(V,U).
\]
Set
\[
\hat{\de}(U,V)= \de(U,V)+\de(V,U).
\]
\begin{proposition} (a) For any  vector subspaces $U,V\subset E$ we have
\[
\|P_U-P_V\|\leq \hat{\de}(U,V)\leq 2\|P_U-P_V\|.
\]
(b) For any vector subspaces $U,V,W$ such that $V\subset W$ we
have
\[
\de(U,V)\geq \de(U,W),\;\; \de(V,U)\leq \de(W,U).
\]
In other words, the function $(U,V)\mapsto \de(U,V)$ is
increasing in the first variable, and decreasing in the second
variable. \label{prop: element}
\end{proposition}

 \proof  (a) We have
 \[
 \hat{\de}(U,V)= \|P_U-P_UP_V\|+ \|P_V-P_VP_U\|
 \]
 \[
 = \|P_U(P_U-P_V)\|+ \|P_V(P_V-P_U)\|\leq 2\|P_U-P_V\|
 \]
 and
 \[
 \|P_U-P_V\|\leq \|P_U-P_UPV\|+\|P_UP_V-P_V\|
 \]
 \[
 = \|P_U-P_UP_V\|+ \|P_V-P_VP_U\|=\hat{\de}(U,V).
 \]
 (b)  Observe that for all $u\in U$, $|u|=1$ we have
 \[
\dist(u,V)\geq \dist(u,W)\Longrightarrow  \de(U,V)\geq \de(U,W).
\]
Since $V\subset W$ we deduce
\[
\sup_{v\in V\setminus 0} \frac{1}{|v|}\dist (v, U)\leq \sup_{w\in
W\setminus 0} \frac{1}{|w|}\dist (w, U).\proofend
\]
We denote by $\Gr_k(E)$ the Grassmannian of $k$ dimensional
subspaces of $E$ equipped with the metric
\[
\dist(U,V) =\|P_U-P_V\|.
\]
$\Gr_k(E)$ is a compact, tame  subset of $\End(E)$. We set
\[
\Gr(E):=\bigcup_{k=0}^{\dim E}\Gr_k(E).
\]
Let  $\Gr^k(E)$ denote the Grassmannian of  codimension $k$
subspaces. For  any subspace $U\subset E$ we set
\[
\Gr(E)_U:=\bigl\{ V\in \Gr(E);\;\;V\supset
U\,\bigr\},\;\;\Gr(E)^U:=\bigl\{ V\in \Gr(E);\;\;V\subset
U\,\bigr\}.
\]
Note that we have a metric preserving involution
\[
\Gr(E)\ni U\longmapsto U^\perp \in \Gr(E),
\]
such that
\[
\Gr_k(E)_U\longleftrightarrow \Gr^k(E)^{U^\perp},\;\;
\Gr^k(E)_U\longleftrightarrow \Gr_k(E)^{U^\perp}.
\]
Using (\ref{eq: dist}) we deduce  that for any  $1\leq j\leq k$,
and any  $U\in \Gr_j(E)$, there exits a constant  $c>1$ such that,
for every $L\in \Gr_k(E)$ we have
\[
\frac{1}{c}\dist(L,\Gr_k(E)_U\,)^2 \leq |P_U-P_UP_L|^2 \leq c
\dist(L,\Gr_k(E)_U\,)^2.
\]
The constant $c$ depends  on $j,k, \dim E$, and a priori it could
also depend on $U$.  Since the  quantities entering into the above
inequality are invariant with respect to the action of the
orthogonal group $O(E)$, and the action of $O(E)$ on $\Gr_j(E)$ is
transitive, we deduce that the constant $c$ is independent on the
plane $U$.  The inequality  (\ref{eq: norm-equiv}) implies the
following result.

\begin{proposition} Let $1\leq j\leq  k \leq \dim E$. There exists a positive constant $c>1$ such that, for any $U\in\Gr_j(E)$, $V\in \Gr_k(E)$ we have
\[
\frac{1}{c} \dist\bigl(\, V, \Gr_k(E)_U)\,\bigr)\leq \de(U, V)
\leq  c \dist\bigl(\, V, \Gr_k(E)_U)\,\bigr).\proofend
\]
\label{prop: de}
\end{proposition}

\begin{corollary}  For every $1\leq k\leq \dim E$ there exists a constant $c>1$ such that, for any $U,V\in \Gr_k(E)$ we have
\[
\frac{1}{c}\dist(U,V)\leq \de(U,V)\leq c \dist(U,V).\proofend
\]
\label{cor: de}
\end{corollary}

\proof In Proposition \ref{prop: de} we make $j=k$ and we observe
that $\Gr_k(E)_U=\{U\}$, $\forall U\in \Gr_k(E)$. \qed

 We would like to  describe a few simple geometric techniques for estimating the gap between two vector subspaces. Suppose  $U, V$ are two vector subspaces  of the Euclidean space $E$ such that
 \[
 \dim U\leq \dim V,\;\de(U,V) <1.
 \]
As remarked earlier, the condition $\de(U,V)<1$ can be rephrased
as $U\cap V^\perp =0$, or equivalently, $U^\perp +V=E$, i.e., the
subspace $V$ intersects  $U^\perp$ transversally. Hence
\[
U\cap\ker P_V=0.
\]
 Denote by $S$ the orthogonal projection of  $U$ on $V$.  We deduce that the restriction of $P_V$ to $U$ defines a bijection $U\ra S$.  Hence $\dim S=\dim U$, and we can find a linear map
\[
h: S\ra V^\perp
\]
whose graph is $U$, i.e.,
\[
U= \bigl\{ s+ h(s);\;\; s\in S,\bigr\}.
\]
Next, denote by $T$ the orthogonal complement of $S$ in $V$ (see
Figure \ref{fig: tame7}), $T:= S^\perp \cap V$, and  by $W$ the
subspace $W:= U+T$.
\begin{figure}[h]
\centerline{\epsfig{figure=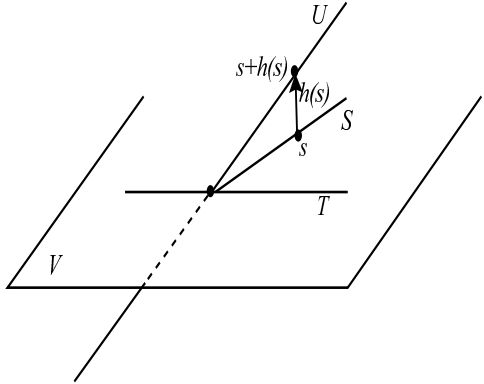,height=2in,width=2.5in}}
\caption{\sl Computing the gap between two subspaces.} \label{fig:
tame7}
\end{figure}

\begin{lemma}
\[
T= U^\perp \cap V.
\]
\end{lemma}

\proof Observe first that
\begin{equation}
(S+U)\subset T^\perp. \label{eq: perp}
\end{equation}
Indeed, let $t\in T$.       Any element in $S+U$ can be written as
a sum
\[
s+u= s+ s' +h(s'),\;\;s,s'\in S.
\]
Then  $(s+s')\perp t$ and $ h(s')\perp t$, because $h(s')\in
V^\perp$.  Hence  $T\subset U^\perp \cap S^\perp\subset U^\perp$.
On the other hand, $T\subset V$ so that
\[
T\subset U^\perp\cap V.
\]
Since $V$ intersects $U^\perp$ transversally we deduce
\[
\dim (U^\perp \cap V)= \dim U^\perp +\dim  V-\dim E= \dim V-\dim
U=\dim T.\proofend
\]

\begin{lemma}
\[
\de(W, V)= \de(U,V)=\de(U,S).
\]
\label{lemma: gap1}
\end{lemma}

\proof  The equality $\de(U,V)=\de(U,S)$ is  obvious. Let $w_0\in
W$ such that $|w_0|=1$ and
\[
\dist(w_0, V)=\de(W,V).
\]
To prove the lemma it suffices to show that $w_0\in U$. We write
\[
w_0= u_0+t_0,\;\;u_0\in U,\;\;t_0\in T,\;\;|u_0|^2+|t_0|^2=1.
\]
We have to prove that $t_0=0$. We can refine the above
decomposition of $w_0$ some more by writing
\[
u_0=s_0+h(s_0),\;\; s_0\in S.
\]
Then
\[
P_Vw_0= s_0+t_0.
\]
We know that for any $u\in U$, $t\in T$  such that $|u|^2+|t^2|$
we have
\[
|u_0^2-P_Vu_0|^2 = |w_0-P_V w_0|^2 \geq
|(u+t)-P_V(u+t)|=|u-P_Vu|^2.
\]
If in the above inequality we choose $t=0$ and $u=\frac{1}{|u|_0}$
we deduce
\[
|u_0^2-P_Vu_0|^2\geq \frac{1}{|u_0|^2} |u_0^2-P_Vu_0|^2.
\]
Hence $|u_0|\geq 1$ and since $|u_0|^2+|t_0|^2=1$ we deduce
$t_0=0$. \qed

The next result summarizes  the   above observations.

\begin{proposition} Suppose $U$ and $V$ are two subspaces of the Euclidean space $E$ such that $\dim U\leq \dim V$ and $V$ intersects $U^\perp$ transversally.  Set
\[
T:= V\cap U^\perp, \;\;W:= U+T,
\]
  and denote by $S$ the orthogonal projection of  $U$ on $V$. Then
\[
S= T^\perp\cap V,
\]
\[
\dim U=\dim S,\;\;\dim W=\dim V,
\]
and
\[
\de(W,V)=\de(U, V)=\de(U,S).\proofend
\]
\label{prop: gap}
\end{proposition}

\begin{proposition} Suppose  $E$ is an Euclidean vector space.  There exists a constant $C>1$, depending only on the dimension of $E$, such that, for any subspaces $U\subset  E$, and any linear operator  $S:U\ra U^\perp$, we have
\begin{equation}
\de(\Gamma_S, U)=\|S\|\bigl(\,1+\|S\|^2\,\bigr)^{-1/2}, \label{eq:
dist-graph0}
\end{equation}
and
\begin{equation}
\frac{1}{C}\|S\|\bigl(\,1+\|S\|^2\,\bigr)^{-1/2}\leq \de(U,
\Gamma_S) \leq C \|S\|\bigl(\,1+\|S\|^2\,\bigr)^{-1/2}, \label{eq:
dist-graph}
\end{equation}
where $\Gamma_S\subset U+U^\perp=E$ is the graph of $S$  defined
by
\[
\Gamma_S:=\bigl\{ \,u+Su\in E;\;\;u\in U\,\bigr\}.
\]
\label{prop: dist-graph}
\end{proposition}

\proof Observe that
\[
\de(\Gamma_S, U)^2=\sup_{u\in U\setminus
0}\frac{|Su|^2}{|u|^2+|Su|^2}=\sup_{u\in U\setminus
0}\frac{(S^*Su,u)}{|x|^2+(S^*Su,u)}.
\]
Choose an orthonormal basis $e_1,\dotsc, e_k$ of $U$ consisting of
eigenvectors of $S^*S$,
\[
S^*Se_i=\lambda_i e_i,\;\;0\leq \lambda_1\leq \cdots\leq
\lambda_k.
\]
Observe that
\[
\|S^*S\|=\lambda_k.
\]
We deduce
\[
\de(\Gamma_S, U)^2=\sup\bigl\{ \,\sum_i \lambda_i
u_i^2;\;\;\sum_i(1+\lambda_i)u_i^2=1\,\bigr\}
\]
\[
=\sup\bigl\{ 1-\sum_i
u_i^2;\;\;\sum_i(1+\lambda_i)u_i^2=1\,\bigr\}
\]
\[
=1-\inf\bigl\{\,\sum_i
u_i^2;\;\;\sum_i(1+\lambda_i)u_i^2=1\,\bigr\}=1-\frac{1}{1+\lambda_k}=\frac{\|S^*S\|}{1+\|S^*S\|}=
\frac{\|S\|^2}{1+\|S\|^2}.
\]
This proves (\ref{eq: dist-graph0}).  The inequality (\ref{eq:
dist-graph}) follows from (\ref{eq: dist-graph0}) combined with
Corollary \ref{cor: de}. \qed

Set
\[
\eP(E)=\bigl\{\, (U,V)\in \Gr(E)\times \Gr(E);\;\; \dim U\leq \dim
V,\;\; V\pitchfork U^\perp\,\bigr\}
\]
For every pair $(U,V)\in\eP(E)$ we denote by $\eS_V(U)$ the
\emph{shadow} of $U$ on $V$, i.e. the orthogonal projection  of
$U$ on $V$.   Let us observe that
\[
U^\perp\cap \eS_V(U)=0.
\]
Indeed, we have
\[
U^\perp \cap \eS_V(U)\subset  T:=U^\perp \cap V \Longrightarrow
U^\perp\cap \eS_V(U)\subset \eS_V(U)\cap T,
\]
and Proposition \ref{prop: gap} shows that $\eS(U, V)$ is the
orthogonal complement of $T$ in $V$.   Since $\dim U=\dim
\eS_V(U)$, we deduce that $\eS_V(U)$ can be represented as the
graph of a linear operator
\[
\eM_V(U): U\ra U^\perp
\]
which we will call the \emph{slope} of the pair $(U,V)$. From
Proposition \ref{prop: gap} we deduce
\[
\de(\eS_V(U),
U)=\frac{\|\eM_V(U)\|}{\bigl(\,1+\|\eM_V(U)\|^2\,\bigr)^{1/2}}\Longleftrightarrow
\|\eM_V(U)\|=\frac{\de(\,\eS_V(U), U)}{\bigl(\,
1-\de(\,\eS_V(U),U)^2\,\bigr)^{1/2}}.
\]
\begin{corollary} There exists a constant $C>1$, which depends only on the dimension of $E$ such that, for every pair $(U,V)\in \eP(E)$ we have
\[
\frac{1}{C}
\|\eM_V(U)\|\bigl(\,1+\|\eM_V(U)\|^2\,\bigr)^{-1/2}\leq
\de(U,V)\leq C
\|\eM_V(U)\|\bigl(\,1+\|\eM_V(U)\|^2\,\bigr)^{-1/2}.
\]
\label{cor: gap}
\end{corollary}

\proof  Use the equality $\de(U,V)= \de\bigl( \,
U,\eS_V(U)\,\bigr)$ and Proposition \ref{prop: dist-graph}. \qed

\begin{proposition} Suppose $A: E\ra E$ is an invertible   symmetric  operator,  and $U$ is  the  subspace of $E$ spanned by  the positive eigenvectors  $A$.  We denote by $m_+(A)$ the smallest positive eigenvalue of $A$, and by $m_-(A)$ the smallest positive eigenvalue of $-A$.    Then, for every subspace $V\subset E$, such that $(U,V)\in \eP(E)$, we have
\[
\de(U, e^{tA}V)\leq e^{-(\, m_+(A)+m_-(A)\,)t}
\|\eM_V(U)\|=e^{-(\, m_+(A)+m_-(A)\,)t}\frac{\de(\,\eS_V(U),
U)}{\bigl(\, 1-\de(\,\eS_V(U),U)^2\,\bigr)^{1/2}}.
\]
\label{prop: smale-flow}
\end{proposition}

\proof  Denote by $L$ the intersection of $V$ with $U^\perp$. Then
we have an orthogonal decomposition
\[
V=L +\eS_V(U),
\]
and if we write $\eM:=\eM_V(U): U\ra U^\perp$ we  obtain
\[
V= \bigl\{ \ell +u +\eM u;\;\;\ell \in L,\;\;u\in U\,\bigr\}.
\]
Using the orthogonal  decomposition $E=U+U^\perp$ we can
describe $A$ in the block form
\[
A=\left[
\begin{array}{cc}
A_+ & 0\\
0 &A_-
\end{array}
\right],
\]
where $A_+$ denotes the restriction of $A$ to $U$, and $A_-$
denotes the restriction of $A$ to $U^\perp$.

Set $V_t:=e^{tA}V$,  $L_t:=V_t\cap U^\perp$. Since $U^\perp$ is
$A$-invariant, we deduce that $L_t=e^{tA_-} L$, so that
\[
V_t= \bigl\{\, e^{tA_-} \ell + e^{tA_+}u +e^{tA_-}\eM u;\;\;\ell
\in L,\;\;u\in U\,\bigr\}
\]
\[
= \bigl\{\, e^{tA_-} \ell + u +e^{tA_-}\eM e^{-tA_+}u;\;\;\ell \in
L,\;\;u\in U\,\bigr\}.
\]
We deduce that for every  $u\in U$ the vector $u +e^{tA_-}\eM
e^{-tA_+}u$ belongs to $V_t$. Hence
\[
\de(U,V_t) \leq \sup_{|u|=1}|e^{tA_-}\eM e^{-tA_+}u|=\|e^{tA_-}\eM
e^{-tA_+}\|\leq e^{-(\, m_+(A)+m_-(A)\,)t}\|\eM\|.\proofend
\]

\begin{corollary} Let $A$ and $U$ as above.   Fix $\ell >\dim U$ and consider  a compact subset $K\subset \Gr_\ell(E)$ such that any $V\in K$ intersects $U^\perp$ transversally.  Then there exists a positive constant, depending only on $K$ and $\dim E$ such that
\[
\de(U, e^{tA}V)\leq C e^{-(\, m_+(A)+m_-(A)\,)t} ,\;\;\forall V\in
K.\proofend
\]
\label{cor: smale-flow}
\end{corollary}

%% file: tameflow7.tex
For any subset  $S$ of a topological space $X$ we will denote by  $\cl(S)$ its closure.

\begin{definition}   Suppose $X, Y$ are two $C^2$-submanifolds  of the Euclidean space $E$ such that $X\subset \cl(Y)\setminus Y$.

\noindent (a) We say that $(X,Y)$ satisfies   \emph{Verdier  regularity condition} $V$ at $x_0\in X$ if there exists an open neighborhood $U$ of $x_0$ in $E$  and a positive constant $C$ such that
\[
\gamma(T_{x}X, T_yY) \leq C|x-y|,\;\;\forall x\in U\cap X,\;\;y\in U\cap Y.
\]
\noindent (b) We say that $(X,Y)$ satisfies the \emph{Verdier  regularity condition} $V$ \emph{along $X$} if it satisfies the   condition $V$ at any point $x\in X$.\qed

\end{definition}

Note that   if  $X$ and $Y$ are connected  and if  $(X,Y)$ satisfies $V$ along $X$, then
\[
\dim X \leq \dim  Y.
\]
As explained in \cite{Verd}, the Verdier condition is invariant under $C^2$-diffeomorphisms.

\begin{remark}  The Verdier regularity condition is equivalent to   the microlocal  regularity condition $\mu$ of Kashiwara and Schapira, \cite[\S 8.3]{KaSch}. For a proof of this fact we refer to \cite{Trot2}.\qed
\end{remark}

The regularity condition $V$ is intimately related to  Whitney's regularity condition.

\begin{definition} Suppose $X, Y$ are two $C^1$-submanifolds  of the Euclidean space $E$ such that $X\subset \bar{Y}\setminus Y$.

\noindent (a) We say that the pair $(X,Y)$ satisfies the  \emph{Whitney regularity condition }   $(a)$  at $x_0\in X$ if, for any  sequence $y_n\in Y$ such that

\begin{itemize}

\item $x_n, y_n\ra x_0$,

\item the sequence of  tangent spaces $T_{y_n} Y$ converges to  the subspace $T_\infty$,

\end{itemize}

\noindent we have $T_{x_0}X\subset T_\infty$.

\noindent(b)  We say that the pair $(X,Y)$ satisfies the  \emph{Whitney regularity condition }   $(b)$  at $x_0\in X$ if, for any  sequence $(x_n,y_n)\in X\times Y$ such that

\begin{itemize}

\item $x_n, y_n\ra x_0$,

\item  the one dimensional subspaces $\ell_n=\bR(y_n-x_n)$ converge to the line $\ell_\infty$,

\item the sequence of tangent spaces $T_{y_n} Y$ converges to  the subspace $T_\infty$,

\end{itemize}

\noindent we have  $\ell_\infty\subset T_\infty$, that is, $\de(\ell_\infty, T_\infty)=0$.

\noindent (c) The pair  $(X,Y)$ is said to satisfy \emph{the regularity condition  $(a)$ or $(b)$ along $X$}, if it satisfies this condition at any $x\in X$.\qed
\end{definition}

The Whitney condition $(a)$ is weaker in the sense that $(b)\Longrightarrow (a)$  and it is fairly easy to construct instances when (a) is satisfied while (b) is violated. The Whitney regularity condition (b) is equivalent with the following geometric condition, \cite{Trot1}.

\begin{proposition}[Trotman]    Suppose  $(X,Y)$  is a pair of $C^1$ submanifolds  of  the $\bR^N$, $\dim X= m$.   Assume $X\subset\bar{Y}\setminus Y$.  Then the pair $(X,Y)$ satisfies the Whitney regularity   condition (b) along $X$ if and only  if, for any  open  set  $U\subset E$, and any    $C^1$-diffeomorphism  $\Psi :  U\ra V$, where $V$ is an open subset of $\bR^N$, such that
\[
\Psi(U\cap X) \subset \bR^m\oplus 0 \subset \bR^N,
\]
the map
\[
 \Psi(Y\cap U)\Lra \bR^m\times (0,\infty),\;\;y\longmapsto \bigl( \,{\rm proj}\,(y)\,, \,\dist(y,\bR^m)^2\,\bigr),
 \]
 is a submersion, where ${\rm proj}:\bR^N\ra \bR^m$ denotes the canonical orthogonal projection.\qed
 \end{proposition}

 For tame objects  the Verdier conditions implies the Whitney condition. More precisely, we have the following result.

\begin{proposition}[Verdier] Suppose $(X,Y)$ is a pair  of  $C^2$, tame submanifold of  the Euclidean space $E$ such that  $X\subset \bar{Y}\setminus Y$.  If $(X,Y)$ satisfies  the regularity condition $V$, then  it also satisfies the regularity  condition $W$.
\label{prop: vw}
\end{proposition}

We include  a detailed  proof of this result because  it   will give the reader a taste of both the significance of the Verdier condition, and  of the    versatility of the tame techniques.

\proof    We follow the approach of \cite{Loi}, which is  an $o$-minimal translation of the  original argument of  Verdier, \cite{Verd}.

First, without any loss of  generality, we  can assume the following.

\begin{itemize}

\item  $X$ is an open  neighborhood of $0$ in some subspace  $F\subset E$.

\item $Y$ is a  tame   $C^2$ submanifold, $Y\subset E\setminus F$,  $\dim Y=n$.

\end{itemize}

Define
\[
W_X(Y):=\bigl\{\, x\in X\cap\cl(Y);\;\;  \mbox{  $(X,Y)$ satisfies the $W$ condition at $x$}\,\bigr\},
\]
\[
V_X(Y)=\bigl\{ x\in X\cap \cl(Y);\;\;\mbox{  $(X,Y)$ satisfies the $V$ condition at $x$}\,\bigr\},
\]
and
\[
I_X(Y):=\bigl\{ (x, y, \ell, T_yY)\in X\times Y\times \Gr_1(E)\times \Gr_n(E);\;\; \ell=\bR\lan y-x\ran\,\bigr\},
\]
where $\bR(y-x)$ denotes the one dimensional subspace  generated by the vector $(y-x)$.

Denote by $\Delta_X\subset X\times X$ the diagonal subset, and  define
\[
J_X=\{ (x,x,\ell, S)\in \Delta_X\times \Gr_1(E)\times \Gr_n(E);\;\; \ell\subset S\}.
\]
If  $\pi_X:\Delta_X \times \Gr_1(E)\times \Gr_n(E)\ra X$ is the natural projection,  then
\[
x\in W_X(Y)\Longleftrightarrow\pi_X^{-1}(x)\cap  \cl\bigl(\,I_X(Y)\,\bigr)\subset J_X(Y).
\]
 We will prove that
\[
V_X(Y) \subset  W_X(Y).
\]
We argue by contradiction. Suppose that there  $\exists x_0\in V_X(Y)\setminus W_X(Y)$. For simplicity, we can assume that $x_0=0\in F$.

\begin{figure}[h]
\centerline{\epsfig{figure=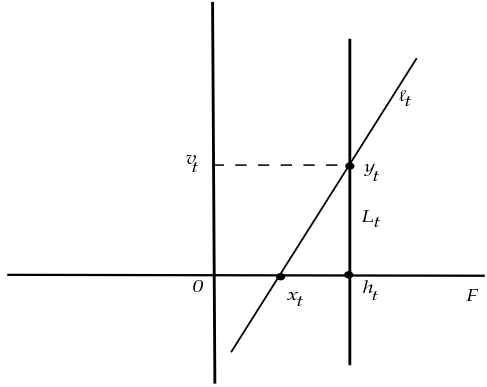,height=2in,width=2.5in}}
\caption{\sl Whitney condition.}
\label{fig: tame6}
\end{figure}

Since $0\not\in W_X(Y)$ we deduce that there exists $(\ell,S)\in \Gr_1(E)\times \Gr_n(E)$ such that
\[
(0,0,\ell,S)\in \cl\bigl(\, I_X(Y)\,\bigr)\;\; {\rm and}\;\; \de(\ell, S) >0.
\]
From the curve selection theorem we deduce that there exits  a   tame $C^1$  path
\[
\alpha: (0,r)\ra I_X(Y),\;\; t\mapsto (x_t,y_t, \bR\lan y_t-x_t\ran, T_{y_t} Y)
\]
such that $x_t\in X$, $y_t\in Y$,
\[
\lim_{t\ra 0^+} x_t=\lim_{t\ra 0^+} y_t=0,\;\;\lim_{t\ra 0^+}  \bR\lan y_t-x_t\ran= \ell,\;\; \lim_{t\ra 0^+}T_{y_t}Y= S.
\]
For simplicity, we write
\[
\ell_t= \bR\lan y_t-x_t\ran,\;\; S_t= T_{y_t} Y.
\]
Using the orthogonal decomposition
\[
E=F\oplus F^\perp
\]
we can decompose $y_t$ into a ``horizontal'' component $h_t\in F$, and a ``vertical'' component, $v_t=y_t-h_t\in F^\perp$ (see Figure \ref{fig: tame6}). Since $y_t \in Y\setminus F$ we deduce $v_t\neq 0$ and thus
\[
\dot{v}_t:=\frac{dv}{dt}(t)\neq 0,\;\;\mbox{for all sufficiently small $t$.}
\]
Moreover, we can choose the parametrization of $\alpha$ so that
\begin{equation}
|\dot{h}_t|\leq 1,\;\;\forall t
\label{eq: bound-h}
\end{equation}

Consider the line $L_t=\bR\lan y_t-h_t\ran$. The tame continuous path $t\mapsto L_t\in \Gr_1(E)$ has a limit as $t\ra 0^+$ which we denote by $L_0$. The limit line $\ell$ is the graph of a linear map
 \[
 f:L_0\ra F\Longrightarrow\ell\subset L_0\oplus F.
 \]
  Since $0\in V_X(Y)$ we deduce that
 \[
 \de(F,S)=0\Longrightarrow F\subset S.
 \]
Hence
\[
L_0\subset S\Longrightarrow S\Longrightarrow \ell \subset S.
\]
Hence $L_0\not\subset S$ so that
\begin{equation}
\de(L_0,S) >0.
\label{eq: l0s}
\end{equation}
To proceed further we need the following auxiliary result.
\begin{lemma} Suppose $u:[0,1)\ra \bR^m$ is a  continuous tame path which is $C^1$ on the open interval $(0,1)$, such that $u(0)=0$ and  $\frac{du}{dt}(t)=\dot{u}(t)\neq 0$ for all sufficiently small  $t$. Then
\begin{equation}
\lim_{t\ra 0^+}\frac{|u(t)|}{|\dot{u}(t)|}=0,
\label{eq: quot1}
\end{equation}
\begin{equation}
\lim_{t\ra 0} t|\dot{u}(t)|=0,
\label{eq: quot2}
\end{equation}
and
\begin{equation}
\lim_{t\ra 0^+}\de(\,\bR u(t),\bR \dot{u}(t)\,)=0.
\label{eq: chord}
\end{equation}
\end{lemma}

\proof Let us first prove (\ref{eq: quot1})  and (\ref{eq: quot2}) in the case $m=1$, i.e., when $u$ is a real valued function.   There exists a small interval $(0,a)\subset (0,1)$ such that $u(t)$  is strictly monotone on $(0,a)$.  For simplicity,  we assume $u(t)$ is strictly increasing, or else we replace $u$ by $-u$.
\[
u(t)>u(0),\;\;\forall t\in (0,a).
\]
By the mean value theorem there exists a   definable function $\xi: (0,a)\ra (0,a)$, $t\mapsto \xi_t$, such that $\xi_t\in (0,t)$ and $\dot{u}(t)=t\dot{u}(\xi_t)$. Hence
\[
0\leq \frac{u(\xi_t)}{\dot{u}(\xi_t)}=t\frac{u(\xi_t)}{u(t)}\leq t\;\; \mbox{and}\;\; 0\leq \xi_t\dot{u}(\xi_t)\leq t \dot{u}(\xi_t)=f(t),
\]
and consequently
\[
\lim_{t\ra 0^+}\frac{u(\xi_t)}{\dot{u}(\xi_t)}=\lim_{t\ra 0^+}\xi_t\dot{u}(\xi_t)=0.
\]
 We can find an even smaller interval $(0,a')\subset (0,a)$ such that $\xi$ is  strictly increasing and $C^1$ on $(0,a')$.  $\xi$ extends to a continuous, strictly increasing function $\xi:[0,a')\ra [0,a')$ such that $\xi_0=0$. Hence $\xi$ is a definable homeomorphism $[0,a)\ra [0,\xi_{a'})$. This   implies (\ref{eq: quot1}) and (\ref{eq: quot2}).

For $m\geq 1$ we write $u(t)=(u_1(t),\dotsc, u_m(t)\,)$, and define
\[
J=\bigl\{ j=1,2,\dotsc ,m;\;\; u_j(t)\;\;\mbox{is not constant on a neighborhood of $0$.}\,\bigr\}.
\]
For $t$ sufficiently small we have
\[
0\leq \frac{|u(t)|^2}{|\dot{u}(t)|^2}=\frac{\sum_{j\in J} u_j(t)^2}{\sum_{j\in J}\dot{u}_j(t)^2} \leq \sum_{j\in J}\frac{u_j(t)^2}{\dot{u}_j(t)^2} \stackrel{t\ra 0^+}{\Lra} 0.
\]
This proves (\ref{eq: quot1}) in the general case. The inequality (\ref{eq: quot2}) is proved in a similar fashion.

To prove  (\ref{eq: chord}) consider the continuous definable functions
\[
s:[0,1)\ra [0,\infty),\;\; s(t)=|u(t)|.
\]
Since $\dot{u}(t)\neq 0$, for $t$ sufficiently small, there exists $b>0$ such that the restriction of $s$ to $(0,b)$ is  strictly increasing and $C^1$.  Hence $s$  defines a $C^1$ diffeomorphism $s: (0,b)\ra (0,s(b))$.
\[
\dot{u}= u' \dot{s},\;u':=\frac{du}{ds},\;\;\dot{s}\neq 0.
\]
Hence $\bR\dot{u} =\bR u'$. We now write
\[
u(s)= sv(s),\;\;|v(s)|=1.
\]
We set $v_0=\lim_{s\ra 0} v(s)$. Then $\lim_{s\ra 0}\bR u(s) =\bR v_0$. On the other hand,
\[
\frac{du}{ds}= v(s) + s \frac{dv}{ds}.
\]
Using (\ref{eq: quot2}) we deduce
\[
\lim_{s\ra 0^+} s\frac{dv}{ds} =0\Longrightarrow \lim_{s\ra 0^+} \frac{du}{ds}=v_0\Longrightarrow\lim_{s\ra 0^+} \bR u'=\bR v_0.\proofend
\]

We have
\[
\de(\bR \dot{v}_t, T_{y_t}Y)=\frac{1}{|\dot{v}_t|} \dist(\dot{v}_t, T_{y_t}Y).
\]
From the equality
\[
\dot{v}_t+\dot{h}_t=\dot{y}_t\in T_{y_t}Y
\]
we deduce
\[
\dist(\dot{v}_t, T_{y_t}Y)=\dist(\dot{h}_t, T_{y_t}Y).
\]
Now observe that $\dot{h}_t\subset F$ and
\[
\dist(\dot{h}_t, T_{y_t}Y)= |\dot{h}_t|\de(\bR\dot{h}_t, T_{y_t}Y)\leq  |\dot{h}_t|\de(F, T_{y_t}Y)=|\dot{h}_t|\de(T_{h_t}X, T_{y_t}Y).
\]
From the Verdier regularity condition we deduce that
\[
|\dot{h}_t|\de(T_{h_t}X, T_{y_t}Y)\leq C |v_t|\,|\dot{h}_t|.
\]
Hence
\[
\de(\bR \dot{v}_t, T_{y_t}Y)=\frac{1}{|\dot{v}_t|} \dist(\dot{v}_t, T_{y_t}Y)\leq C|\dot{h}_t| \frac{|{v}_t|}{|\dot{v}_t|}\stackrel{(\ref{eq: bound-h})}{\leq} C\frac{|{v}_t|}{|\dot{v}_t|}.
\]
Using (\ref{eq: quot1}) we deduce
\[
\lim_{t\ra 0^+}\de(\bR \dot{v}_t, T_{y_t}Y)=0.
\]
On the other hand,    using (\ref{eq: chord})  we deduce
\[
\lim_{t\ra 0^+}\bR \dot{v}_t=\lim_{t\ra 0^+} \bR v_t=L_0
\]
from which we conclude that
\[
0=\de(L_0, \lim_{t\ra 0^+}T_{y_t}Y)=\de(L_0, S).
\]
This contradicts (\ref{eq: l0s}),  and completes the proof of Proposition \ref{prop: vw}.\qed

\begin{definition} Suppose $X$ is a  subset of an Euclidean space $E$.   A \emph{Verdier   stratification}  (respectively \emph{Whitney stratification}) of $X$ is an increasing, finite   filtration
\[
F_{-1}=\emptyset \subset F_0\subset F_1 \subset \cdots\subset F_m=X
\]
satisfying the following  properties.

\smallskip

\noindent (a) $F_k$ is closed in $X$, $\forall k$.

\noindent (b) For every $k=1,\dotsc, m$ the  set $X_k=F_k\setminus F_{k-1}$ is a $C^2$ manifold  of dimension $k$. Its connected components are called  the \emph{strata} of the stratification.

\noindent (c) (The frontier condition) For every $k=1,\dotsc, m$ we have
\[
\cl(X_k)\setminus X_k\subset F_{k-1}.
\]
\noindent (d) For every $0\leq j <k\leq m$ the pair $(X_j, X_k)$ is Verdier regular (respectively Whitney regular) along $X_j$.

 If $X$ is a tame set, then a Verdier (Whitney) stratification is called  tame if the sets $F_k$ are tame. \qed
\end{definition}

We have the following result due to essentially to Verdier  \cite{Verd} (in the subanalytic case) and Loi \cite{Loi} in the tame context.

\begin{theorem} Suppose $S_1,\dotsc, S_n$ are tame subsets of the Euclidean  space $E$. Then there exists a  tame  Verdier stratification of $E$ such that each of the sets $S_k$ is a  union of strata.\qed
\end{theorem}

\begin{remark} According to the results of Lion and Speissegger \cite{LS}, the strata in the above Verdier stratification can be chosen to be \emph{real analytic} submanifolds of $E$.\qed
\end{remark}

 A Whitney stratified space  $X$  has   a rather  restricted local structure.   More precisely, we have the following fundamental result whose  intricate proof can be found in \cite[Chap,II,\S 5]{GWPL}.

 \begin{theorem} Suppose  $X$ is a  subset of a   smooth manifold $M$ of dimension $m$, and
 \[
 F_0\subset F_1\subset \cdots \subset F_k= M
 \]
 is a Whitney stratification  of $X$. Then for   every stratum $S$  of dimension $j$ there exists

 \begin{itemize}

 \item  a closed tubular neighborhood $N$ of $S$ in $M$ with projection $\pi:N\ra S$,

 \item  a Whitney stratified  subset  $L_S$ of the sphere $S^{m-j-1}$

 \end{itemize}

\noindent  such that $\pi: \pa N\ra X$ is a  locally trivial  fibration     with fiber  homeomorphic to $L_S$, and  $N\cap X$ is homeomorphic with the mapping cylinder of the projection $\pi:\pa N\ra S$. The space $L_S$ is called the normal link of $S$ in $X$. \qed
\label{th: local-tr}
 \end{theorem}

%% file: tameflow8.tex
Suppose $M$ is a compact, connected  real analytic manifold of dimension $M$, and $(f,\xi)$ is a Morse pair on  $M$, \emph{not necessarily tame}.  Denote by $\Phi^\xi$ the flow generated by $\xi$,  by $W^-_p(\xi)$ (respectively $W^+_p(\xi)$) the   unstable (respectively stable) manifold of the critical point $p$, and    set
\[
M_k(\xi): =\bigcup_{p\in \Cr_f,\;\;\lambda(p)\leq k} W^-_p(\xi),\;\;W^-_k(\xi)=M_k(\xi)\setminus M_{k-1}(\xi).
\]
We say that the flow $\Phi^\xi$  satisfies the \emph{Morse-Verdier}   condition if the increasing filtration
\[
M_0(\xi)\subset M_1(\xi)\subset\cdots \subset M_m(\xi)
\]
is a Verdier regular stratification. In the sequel,  when no confusion is possible, we will write $W^\pm_p$ instead of $W_p^\pm(\xi)$.

\begin{theorem} If the Morse flow $\Phi^\xi$  satisfies the Morse-Verdier condition, then it also satisfies the Morse-Smale condition.

\label{th: vis}
\end{theorem}

\proof   Let $p,q\in \Cr_f$ such that $p\neq q$ and   $W^-_p\cap W^+_q\neq \emptyset$. Let $k$ denote the Morse index of $q$, and $\ell$ the Morse index of $q$ so that $\ell> k$. We want to prove that this intersection is  transverse.

Let $x\in W^-_p\cap W^+_q$ and set
\[
x_t :=\Phi_t^\xi(x).
\]
Observe that
\[
T_{x}W^+_q\pitchfork  T_xW^-_p\Longleftrightarrow \exists t \geq 0:\;\;T_{x_t}W^+_q\pitchfork  T_{x_t}W^-_p.
\]
We will prove  that $T_{x_t}W^+_q\pitchfork  T_{x_t}W^-_p$  if $t$ is sufficiently large.

 Since $(f,\xi)$ is a Morse pair, we can find coordinates $(u^i)$ in a neighborhood $U$ of $q$, and   real numbers
\[
\mu_1,\dotsc, \mu_m>0
\]
such that
\[
u^i(q)=0,\;\;\forall i,
\]
\[
\xi=  \sum_{i=1}^k \mu_i u^i\pa_{u_i} -\sum_{\alpha>k}\mu_\alpha u^\alpha\pa_{u_\alpha}.
\]
Denote by $A$ the diagonal matrix
\[
A=\diag(\mu_1,\dotsc,\mu_k,-\mu_{k+1},\dotsc,-\mu_m).
\]
Without any loss of generality, we can assume  that the point $x$ lies in the coordinate neighborhood $U$.  Denote by $E$ the Euclidean space with Euclidean coordinates $(u^i)$. Then  the path
\[
t\mapsto T_{x_t}W^-_p\in \Gr_\ell(E)
\]
is given by
\[
T_{x_t}W^-_p=e^{tA}T_{x}W^-_p,
\]
and in particular it has a limit
\[
\lim_{t\ra \infty} T_{x_t}W^-_p=T_\infty\in \Gr_\ell(E).
\]
 Since the pair $\bigl(\, W^-_q, W^-_p\,\bigr)$ satisfies the Verdier regularity condition along $W^-_q$, and $x_t\ra q$, as $t\ra \infty$, we deduce
\[
T_\infty\supset T_q W^-_q,\Longrightarrow T_\infty \pitchfork T_q W^+_q.
\]
Thus, for $t$ sufficiently large
\[
T_{x_t}W^-_p\pitchfork T_{x_t}W^+_q.\proofend
\]
Suppose $(f,\xi)$ is a Morse pair on the compact, real analytic manifold $M$.   Then for  every critical point  $p$ of $f$ of index $k$  we can find  local $C^2$-coordinates $(u^i)$ defined in an open neighborhood $U_p$, and positive real numbers $\mu_i$  such that
\[
u^i(p) =0,\;\;\forall i,
\]
and
\[
\xi=\sum_{i\leq k} \mu_i u^i\pa_{u_i} -\sum_{\alpha>k}\mu_\alpha u^\alpha\pa_{u_\alpha}.
\]
If $p$ is a hyperbolic point, i.e., $0<k<m$, we set,
\[
\gamma_u(p)=\gamma_u(\xi,p):= \min_{i\leq k}\mu_i,\;\gamma_s(p)=\gamma_s(\xi,p):=\min_{\alpha>k}\mu_\alpha,\;\; \Gamma_s(p)=\Gamma_s(\xi,p):=\max_{\alpha>k}\mu_\alpha,
\]
\[
g_s(p)=g_s(\xi,p):=\Gamma_s(p)-\gamma_s(p).
\]
Observe that $g_s(p)$ is  the length of the smallest interval containing all the negative (or stable) eigenvalues of the linearization of $\xi$ at $p$, while  $\gamma_u(p)$ is the smallest positive (or unstable) eigenvalue  of the linearization  of $\xi$ at $p$.

\begin{theorem} Suppose $(f,\xi)$ is a Morse pair on the real analytic manifold $M$ of dimension $m$ such that the flow  $\Phi^\xi$ satisfies the Morse-Smale condition. Define
\begin{equation}
\nu:=\min\bigl\{ \frac{\gamma_u(p)+\gamma_s(p)}{\Gamma_s(p)};\;\;p\in\Cr_f,\;\; 0<\lambda(p)<\dim M\,\bigr\}.
\label{eq: gap-nu}
\end{equation}
Assume $\xi$ is at least $(\lfloor\nu\rfloor+1)$-times differentiable.  Then the following hold.

\noindent (a)(Frontier property) $\cl\bigl(\,M_k(\xi)\,\bigr)\setminus M_k(\xi)\subset M_{k-1}(\xi)$, $\forall k$.

\noindent (b) For every pair of critical points $p,q$, and every  $z\in W^-_q\cap \cl( W^-_p)$, there exists an open neighborhood  $U$ of $z\in M$, and a positive constant $C$ such that
\begin{equation*}
\de\bigl( T_xW^-_q, T_yW^-_p\,\bigr)\leq C \dist(x,y)^\nu,;\;\forall x\in U\cap W^-_q, \;\;\forall y\in U\cap W^-_p.
\tag{$V_\nu$}
\label{tag: verd-nu}
\end{equation*}

\label{th: siv}
\end{theorem}

\begin{remark} (a) It is perhaps useful to visualize the condition (\ref{eq: gap-nu})  in which $\nu\geq 1$, as a  \emph{spectral clustering condition}.

Suppose  $p$ is an unstable   critical point of $f$. Denote  $H_p$ is the Hessian of $f$ at $p$. Using the metric $g$ we can identify $H_p$ with a symmetric operator. We denote by $\Sigma_p^\pm$ the collection of positive/negative eigenvalues of this operator. Then $\gamma_s(\xi,p)$ is the positive spectral gap,
\[
\gamma_s(\xi,p)=\min\Sigma^+_p=\dist(\Sigma^+_p,0),
\]
$\gamma_u(\xi,p)$ is the negative spectral gap
\[
\gamma_u(\xi,p)=\dist(\Sigma^-_p,0),
\]
and $\Gamma_s(\xi,p)$ is  the largest positive eigenvalue of $H_p$.  The condition
\[
\frac{\gamma_u(\xi,p)+\gamma_s(\xi,p)}{\Gamma_s(\xi,p)}\geq 1,
\]
than says that the length of smallest interval containing the positive spectrum  is  smaller than the length  largest interval containing $0$, and disjoint from the spectrum. In particular, if the positive eigenvalues   cluster in an tiny interval situated  far away from the origin, this condition is automatically satisfied.

\begin{figure}[h]
\centerline{\epsfig{figure=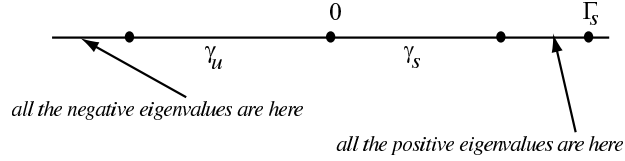,height=0.8in,width=3.2in}}
\caption{\sl Spectral gaps.}
\label{fig: tame26}
\end{figure}

(b)  If the    unstable manifolds of a Morse flow on a  compact smooth manifold  $M$  form a Whitney stratification,  then  Theorem \ref{th: local-tr} shows that the closure of any unstable manifold  is a  \emph{submanifold with conical singularities} in the sense of  \cite{Lauden}.\qed
\end{remark}

\proof  To prove part (a) it suffices to show that if
\[
W^-_q\cap \cl\bigl(\,W^-_p\,\bigr)\Longrightarrow \dim W^-_q <\dim W^-_p.
\]
Observe that the set $W^-_q\cap \cl\bigl(\,W^-_p\,\bigr)$ is flow  invariant, and its intersection with any compact subset of $W^-(p,\xi)$ is closed. We deduce that $p\in W^-_q\cap \cl\bigl(\,W^-_p\,\bigr)$.

Fix  a small neighborhood $U$ of $p$ in $W^-_p$.  Then there exists a  sequence $x_n\in \pa U$, and a sequence   $t_n\in[0,\infty)$, such that
\[
\lim_{n\ra \infty} t_n =\infty,\;\;\lim_{n\ra \infty} \Phi^\xi_{t_n}x_n= q.
\]
In particular, we deduce that $f(p)>f(q)$.

For every $n$  define
\[
C_n=\cl\bigl(\, \bigl\{ \Phi_t^\xi x_n;\;\;t\in (-\infty,t_n]\bigr\}\,\bigr).
\]
Denote by $\Cr^p_q$ the set of critical points $p'$ such that $f(q)<f(p')<f(p)$. For every $p'\in\Cr^p_q$ we denote  by $d_n(p')$ the distance from $p'$ to $C_n$.  We can find  a set $S\subset\Cr^p_q$ and a subsequence  of the sequence $(C_n)$,  which we continue to denote by $(C_n)$, such that
\[
\lim_{n\ra \infty} d_n(p')=0,\;\;\forall p'\in S\;\;\mbox{and}\;\;\inf_n d_n(p')>0,\;\;\forall p'\in \Cr^p_q\setminus S.
\]
Label the  points in $S$ by $s_1,\dotsc, s_k$,  so that
\[
f(s_1)>\cdots >f(s_k).
\]
Set $s_0=p$, $s_{k+1}=q$.  The critical points in $S$  are hyperbolic,  and we conclude  that  there exist trajectories $\gamma_0,\dotsc, \gamma_k$  of $\Phi$,  such  that
\[
\lim_{t\ra-\infty}\gamma_i(t)=s_i,\;\;\lim_{t\ra \infty}\gamma_i(t)=s_{i+1},\;\;\forall i=0,\dotsc,k,
\]
and
\[
\liminf_{n\ra \infty}\dist(C_n, \Gamma_0\cup\cdots \cup\Gamma_k)=0,
\]
where $\Gamma_i=\cl\bigl(\,\gamma_i(\bR)\,\bigr)$, and $\dist$ denotes the  Hausdorff distance.     We deduce
\[
W^-_{s_i}\cap W^+_{s_{i+1}}\neq \emptyset,\;\;\forall i=0,\dotsc,k.
\]
Since the flow $\Phi^\xi$ satisfies the Morse-Smale condition  we deduce from the above that
\[
\dim W^-_{s_i}> \dim W^-_{s_{i+1}},\;\;\forall i=0,\dotsc,k.
\]
In particular,
\[
\dim W^-_p>\dim W^-_q.
\]

To prove (b), observe first  that since the map $x\mapsto\Phi_t(x)$ is $(\lfloor\nu\rfloor+1)$-times differentiable for every $t$, the set of points $z \in W^-_p\cap \cl( W^-_q)$ satisfying (\ref{tag: verd-nu}) is open  in $W^-_q$ and flow invariant.  Since $q\in  \cl(W^-_p)\cap \cl( W^-_q)$ it suffices to prove (b) in the special case $z=q$. We will achieve this using an inductive argument.

For every $0\leq k\leq m=\dim M$ we denote by $\Cr^k_f$ the set of index $k$ critical points of $f$.  We will prove by decreasing induction over $k$  the following statement.

\begin{description}
\item[$\boldsymbol{S}(k)$]\emph{For every  $q\in \Cr^k_f$, and every $p\in \Cr_f$ such that $q\in \cl\bigl(\, W^-_p\,\bigr)$ there exists a neighborhood $U$ of $q\in M$, and a constant $C>0$ such that (\ref{tag: verd-nu}) holds.}
\end{description}

The statement is vacuously true when $k=m$.  We fix $k$, we assume that $\boldsymbol{S}(k')$  is true for any $k'>k$, and we will prove that the statement  its is true for $k$ as well.  If $k=0$ the statement is trivially true because  the distance between the trivial subspace and  any other subspace of a vector space is always zero. Therefore, we can assume $k>0$.

\begin{figure}[h]
\centerline{\epsfig{figure=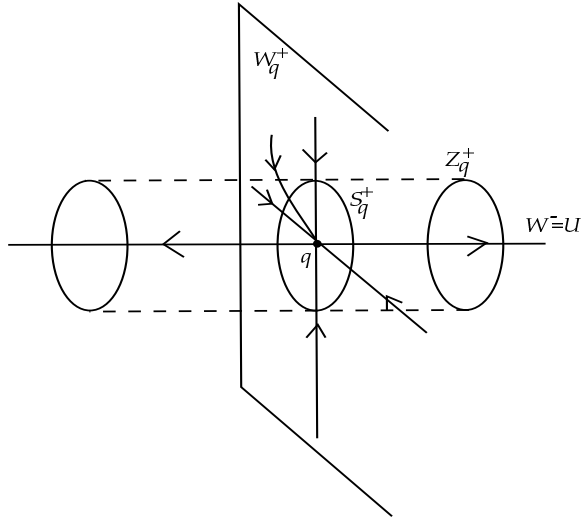,height=2.7in,width=3.0in}}
\caption{\sl The dynamics in a neighborhood of a hyperbolic point.}
\label{fig: tame9}
\end{figure}

Fix $q\in \Cr^k_f$, $\ell>k$ and $p\in \Cr^\ell_f$.   Fix adapted coordinates $(u^i)$  defined  in a  neighborhood  of $\eN$ of $q$ such that, there exist positive real numbers $R,\mu_1,\dotsc, \mu_m$ with the property
\[
\xi= -\sum_{i\leq k}\mu_i u^i\pa_{u^i}  +\sum_{\alpha >k}\mu_\alpha u^\alpha\pa_{u_\alpha},
\]
and
\[
\bigl\{\, (\,u^1(x),\dotsc, u^m(x)\,)\in\bR^m;\;\;x\in\eN\,\bigr\} \supset [-R,R]^m.
\]
For every $r\leq R$ we set
\[
\eN_r:= \bigl\{\, x\in \eN;\;\;|u^j(x)|\leq r,\;\;\forall j=1,\dotsc, m\,\bigr\},
\]
For every $x\in \eN_R$ we define, its horizontal and vertical components,
\[
\bh(x)= (u^1(x),\cdots, u^k(x))\in \bR^k,\;\;  \bv(x)=(u^{k+1}(x),\dotsc,u^m(x))\in\bR^{m-k}.
\]
Define (see Figure \ref{fig: tame9})
\[
S^+_q(r):=\bigl\{\,x\in W^+_q\cap \eN_r;\;\;|\bv(x)|=r\,\bigr\},\;\; Z^+_q(r)=\bigl\{\,x\in \eN_r;\;\;|\bv(x)|=r\,\bigr\}.
\]
The set $Z^+_q(r)$ is the boundary of a ``tube'' of radius $r$  around the unstable manifold $W^-_q$.

We denote by $U$  the vector subspace of $\bR^m$   given by $\{\bv(u)=0\}$, and by $U^\perp$ its orthogonal complement. Observe that for every $x\in W^-_q\cap\eN_R$ we have $T_xW^-_p=U$.

Finally, for $k'>k$ we denote by $\eT_{k'}(U^\perp)\subset \Gr_{k'}(\bR^m)$  the set of $k'$-dimensional subspaces of $\bR^m$ which intersect $U^\perp$ transversally.

 From part (a) we deduce that there exists  $r\leq R$
\begin{equation}
\eN_r\cap \cl\bigl(\, W^-_{q'}\,\bigr)=\emptyset,\;\;\forall j\leq k,\;\;\forall  q'\in\Cr^j_f,\;\;q'\neq q.
\label{eq: disj}
\end{equation}
For every critical point $p'$  we set
\[
C(p',q)_r:= C(p',q)\cap S^+_q(r).
\]
Now consider the set
\[
X_r(q):= C(p,q)_r\cup \bigcup_{k<\lambda(p')<\ell} C(p',q)_r.
\]
For any positive number $\hbar$   we set
\[
\eG_{r,\hbar}:=\cl\bigl(\,\bigl\{  T_xW^-_p;\;\;x\in Z^+_q(r);\;\;|\bh(x)|\leq \hbar\,\bigr\}\,\bigl)\subset \Gr_\ell(\bR^m).
\]
\begin{lemma} There exists a positive $\hbar\leq r$ such that
\[
\eG_{r,\hbar}\subset \eT_\ell(U^\perp).
\]
\end{lemma}

\proof  We argue by    contradiction. Assume that there exists  sequences $\hbar_n\ra 0$ and $x_n \in \eN_r$ such that
\[
|\bv(x_n)|=r,\;\;|\bh(x_n)|\leq \hbar_n,\;\;\de(U, T_{x_n}W^-_p)\geq 1-\frac{1}{n}.
\]
By extracting subsequences we can assume that $x_n\ra x\in S^+_q(r)$  and $T_{x_n}W^-_p\ra T_\infty$ so that
\begin{equation}
\de(U, T_\infty)=1 \Longleftrightarrow\mbox{ $T_\infty$ does not intersect $U^\perp$ transversaly.}
\label{eq: con-tra}
\end{equation}
From the frontier condition and (\ref{eq: disj}) we deduce $x\in X_r(q)$. If $x\in C(p,q)_r$ then   $x\in W^-_p\cap S^+_q(r)$, and we deduce  $T_\infty=T_xW^-_p$. On the other hand,   the Morse-Smale condition shows that  $T_xW^-_p$ intersects  transversally $T_xW^+_q=U^\perp$ which contradicts (\ref{eq: con-tra}).

Thus $x\in C(p',q)$ with $\lambda(p')=k'$, $k<k'<\ell$.    Since we assume that the statement $\boldsymbol{S}(k')$ is true,  we deduce $\de(T_xW^-_{p'},T_\infty)=0$, i.e.,
\[
T_\infty \supset T_xW^-_{p'}.
\]
From the Morse-Smale condition we deduce   that $T_xW^-_{p'}$ intersects  $T_xW^+_q=U^\perp$ transversally, and a  fortiori, $T_\infty$ will  intersect $U^\perp$ transversally. This again  contradicts (\ref{eq: con-tra}).\qed

Fix  $\hbar\in(0,r]$ such that  the   compact set
\[
\eG_{r,\hbar}= \bigl\{\, T_x W^-_p; x\in W^-_p\cap Z^+_q(r),\;\; |\bh(x)|\leq \hbar\,\bigr\}\subset \Gr_{\ell}(\bR^m)
\]
is a subset of $\eT_{\ell}(U^\perp)$. Consider the block
\[
\eB_{r,\hbar}=\bigl\{\, x\in \eN_r;\;|\bv(x)|\leq r,\;\;|\bh(x)|\leq \hbar\,\bigr\}.
\]
The set $\eB_{r,\hbar}$ is a compact neighborhood of $q$. Define
\[
A_u:\bR^k\ra \bR^k,\;\;A_u=\diag(\mu_1,\dotsc,\mu_k),
\]
\[
A_s:\bR^{m-k}\ra \bR^{m-k},\;\;A_s=\diag(\mu_{k+1},\dotsc, \mu_m),
\]
\[
A:\bR^m\ra \bR^m,\;\;A=\diag(A_u, -A_s).
\]
For every $x\in \eB_{r,\hbar}\setminus W^-_q$  we denote by $I_x$ the connected component of
\[
\{ t\leq 0;\;\;\Phi_t^\xi x\in \eB_{r,\hbar}\}
\]
which contains $0$. The set $I_x$ is a closed interval
\[
I_x:=[-T(x), 0],\;\;T(x)\in  [0,\infty].
\]
If $x\in \eB_{r,\hbar}\setminus W^-_q$ then $ T(x) <\infty$. We set
\[
z(x):= \Phi_{-T(x)}^\xi x,\;\;y(x):= v(z(x)).
\]
Then
\[
y(x)= e^{T(x) A_s} v(x),\;\; |y(x)|=r.
\]
We deduce
\[
|v(x)| =|e^{-T(x) A_s}y(x)| \geq   e^{-\Gamma_s(q) T(x)}|y(x)|= e^{-\Gamma_s(q)T(x)}r.
\]
Hence
\begin{equation}
e^{-\Gamma_s(q) T(x)}\leq \frac{1}{r} |v(x)|.
\label{eq: dist-st}
\end{equation}

Let $x\in \eB_{r,\hbar}\cap W^-_p$. Then
\[
T_xW^-_p= e^{T(x)A} T_{z(x)}W^-_p,\;\;  T_{z(x)}W^-_p\in \eG_{r,\hbar}
\]
and,  we deduce
\[
\de(U, T_xW^-_p)=\de(U, e^{T(x)A} T_{z(x)}W^-_p),\;\;U=T_qW^-_q.
\]
Using Corollary \ref{cor: smale-flow} we deduce that there exists a constant $C>0$ such that for every $V\in  \eG_{r,\hbar}$, and every  $t\geq 0$ we have
\[
\de(U, e^{tA}V)\leq  Ce^{-(\,\gamma_s(p)+\gamma_u(p)\,)t}
\]
Hence
\[
\forall x \in \eB_{r,\hbar}\cap W^-_p:\;\;\de(U, T_xW^-p)\leq C e^{-(\,\gamma_s(q)+\gamma_u(q)\,)T(x)}.
\]
Now observe that
\[
-(\gamma_s(q)+\gamma_u(q))\leq -\nu \Gamma_s(q)
\]
so that
\[
e^{-(\,\gamma_s(q)+\gamma_u(q)\,)T(x)}\leq  e^{-\nu\Gamma_s(q)T(x)} \stackrel{(\ref{eq: dist-st})}{\leq}\frac{1}{r^\nu} |\bv(x)|^\nu.
\]
We conclude that
\[
\forall x \in \eB_{r,\hbar}\cap W^-_p:\;\;\de(U, T_xW^-_p)\leq C \frac{1}{r^\nu} |\bv(x)|^\nu=\frac{C}{r^\nu} \dist(x, W^-_q)^\nu.
\]
Since for every $w\in \eB_{r,\hbar}\cap W^-_q$ we have $U=T_wW^-_q$, the last inequality proves $\boldsymbol{S}(k)$.\qed

\begin{corollary} Suppose $(f, \xi)$ is a  smooth Morse pair on the real analytic  manifold $M$ such that the flow $\Phi^\xi$ generated by $\xi$ satisfies the Morse-Smale condition,  and for every   hyperbolic critical point  $p$ we have
\[
\gamma_u(p)+\gamma_s(p) \geq \Gamma_s(p)\Longleftrightarrow  \gamma_u(p)\geq  \Gamma_s(p)-\gamma_s(p).
\]
Then  the filtration
\[
M_0(\xi)\subset M_1(\xi)\subset \cdots \subset M,\;\;M_k(\xi):=\bigcup_{\lambda(p)\leq k} W^-_p(\xi)
\]
is a Verdier stratification. In particular, if the flow $\Phi^\xi$ is also tame, then the above stratification  satisfies the Whitney regularity conditions as well.
\end{corollary}

From Theorem \ref{th: spec}  and  Theorem \ref{th: MS} we obtain  the following result.

\begin{corollary}  Suppose $M$ is a compact real analytic manifold of dimension $m$, $f: M\ra \bR$ is a real analytic Morse function, and $\nu$ is a positive  real number.   Then  there exist

\begin{itemize}

\item a real analytic metric $g$ on $M$,

\item a  smooth vector field $\xi$ on $M$,

\end{itemize}

\noindent such that
\begin{itemize}

\item $(f,\xi)$ is a Morse pair,

\item $\xi$  coincides with $-\nabla^g f$  in an neighborhood of the  critical set,

\item  the flow $\Phi^\xi$ generated by  $\xi$ is tame and satisfies  the Morse-Smale condition,

\item  for every  hyperbolic critical point $p$ of $f$ we have
\begin{equation*}
\frac{\gamma_u(\xi,p)+\gamma_s(\xi,p)}{\Gamma_s(\xi,p)}\geq \nu
\tag{$\gamma_\nu$}
\label{tag: gnu}
\end{equation*}
\end{itemize}

In particular, if $\nu\geq 1$, then  the stratification of $M$ by the unstable  manifolds of the flow $\Phi^\xi$ is   both Verdier and Whitney regular. \qed

\label{cor: vw}
\end{corollary}

\begin{remark} (a) Theorem \ref{th: siv} is  not optimal.  To see this consider the projective space $\bRP^n=\Gr_1(\bR^{n+1})$. We regard it as a  submanifold in the Euclidean space of symmetric operators $\bR^{n+1}\ra \bR^{n+1}$.

Any symmetric operator $A:\bR^{n+1}\ra \bR^{n+1}$ defines  a  function
\[
f_A: \bRP^n\ra \bR,\;\;L\mapsto \tr AP_L.
\]
Suppose
\[
A=\diag(\lambda_0,\dotsc, \lambda_n),\;\; \lambda_0<\lambda_1<\cdots <\lambda_n.
\]
Using the projective coordinates   $[x_0,\dotsc,x_n]$ on $\bRP^n$, we    can describe the critical points of $f_A$ as
\[
\Cr_A= \bigl\{\, p_0,\dotsc, p_n;\;\; p_i= [\delta_{i0},\delta_{i1},\dotsc,\delta_{in}]\,\bigr\}
\]
where $\delta_{ij}$ is the Kronecker symbol.

The eigenvalues of the Hessian of $f$ at $p_i$ are
\[
\mu_j=\lambda_j-\lambda_i,\;\; j\neq i.
\]
The hyperbolic critical points are $p_1,\dotsc, p_{n-1}$. The spectral clustering condition ($\gamma_{\nu=1}$) at $p_i$ reads
\[
\lambda_{i+1}-\lambda_{i-1} \geq \lambda_n-\lambda_i\Longleftrightarrow\lambda_i-\lambda_{i-1} \geq \lambda_n-\lambda_{i+1}.
\]
This condition is  satisfied if for example we choose $\lambda_i$ such that
\[
(\lambda_{i+1}-\lambda_i)\ll (\lambda_i-\lambda_{i-1}),\;\;\mbox{e.g.,}\;\; (\lambda_{i+1}-\lambda_i)< \frac{1}{i+1} (\lambda_i-\lambda_{i-1}).
\]
but fails in the case $\lambda_i=i$.

However, the unstable manifolds of the critical points are  \emph{independent} of the choice of $\lambda_i$.   In fact, these  unstable  varieties are the Schubert cells.
\[
W_i=\bigl\{\,[x_0,\dotsc, x_{i-1}, 1,0,\dotsc, 0];\;\; x_j\in \bR\,\bigr\}.
\]
By choosing $\lambda_i$ so that the clustering condition is satisfied, we deduce that the  unstable manifolds  satisfy the Verdier regularity condition, and they do so even when the spectral clustering condition is violated.

(b)  Although the clustering condition  is  not optimal,  it is in some sense necessary.  To understand this, suppose   we are on a compact, real analytic $3$-manifold $M$,  $f: M\ra \bR$ is a tame $C^2$-Morse function, and $\xi$ is a vector field on $M$ which generates a tame flow $\Phi_t$, and $\xi\cdot f<0$ on $M\setminus \Cr_f$.

 Suppose $q_0\in M$ is a  critical point of $f$ of index $1$, the Hessian of $f$ at $q_0$ has eigenvalues $-1,1,3$, and in a neighborhood of $q_0$  we can find real analytic coordinates  $(x,y,z)$ such that
\[
x(q_0)=y(q_0)=z(q_0)=0,\;\;\xi= x\pa_x-y\pa_y-3z\pa_z.
\]
Observe that the spectral clustering condition is violated since
\[
\gamma_s(q_0)+\gamma_u(q_0)= 2< \Gamma_s(p_0)= 3.
\]

Suppose the point $q=(0,0,1)\in W^+_{q_0}$ also  lies on the unstable  variety $W^-_p$ of a critical point $p$ of index $2$,  $q\in W^+_{q_0}\pitchfork W^-_p$.  Set $q_t= \Phi_t(q)$.  Then $q_t=(0,0,e^{-3t})\in W^-_p$ so that
\[
-3\pa_z=\dot{q}_0\in T_qW^-_p.
\]
Since  $W^-_p$ intersects  $W^+_{q_0}$  transversally at $q$ we deduce
\[
T_q W^-_p = {\rm span}\{ \pa_z, \pa_x+a\pa_y\}.
\]
Assume $a\neq 0$. Then
\[
T_0 W^-_{q_0}={\rm span}\,\{\pa_x\},\;\;T_{q_t}W^-_p= {\rm span}\,\bigl\{ \, e^{-3t} \pa_z,  e^{t}\pa_x+ e^{-t}a\pa_y \,\bigr\}=  {\rm span}\,\bigl\{ \, \pa_z,  \pa_x+ e^{-2t}a \pa_y\,\bigr\}
\]
We deduce that
\[
\delta\bigl(\,T_0W^-_{q_0},\,T_{q_t}W^-_p\,\bigr)\sim |a| e^{-2t},\;\;\mbox{as}\;\;t\ra \infty,
\]
so that
\[
\lim_{t\ra \infty} \frac{\delta\bigl(\,T_0W^-_{q_0},\,T_{q_t}W^-_p\,\bigr)}{{\rm dist}\,(0,q_t)}= \lim_{t\ra\infty}e^t= \infty.
\]
(c) We are inclined to believe that the Morse-Smale condition is  in fact equivalent with the Morse-Whitney condition, i.e., with the fact  that the  unstable manifolds define a Whitney stratification, but we have not attempted a proof. \qed
\end{remark}

%% file: tameflow9.tex
In this section we want to investigate  the Conley indices  of  the stationary points of  certain tame flows.  We begin with a  fast introduction to  Conley theory. For more details we refer to \cite{Con, Sal}.

Suppose $X$ is a compact metric space, and $\Phi: \bR\times X\ra X$, $(t,x)\mapsto \Phi_t(x)$  is a continuous flow.    Thinking of this flow as an action of $\bR$ on $X$, we will  denote $\Phi_t(x)$ by $t\cdot x$.     For  any set $W\subset X$ we define
\[
I^\pm(W):=\bigl\{\, x\in N;\;\; t\cdot x\in W,\;\;\forall t,\;\;\pm t\geq 0\,\bigr\},\;\;I(W):=I^+(W)\cap I^-(W).
\]
 An \emph{isolated} invariant set of the flow is a  closed   a closed,  flow invariant subset $S\subset X$ such that there exists  a compact neighborhood  $W$ of $S$ in $X$ with the property that $S=I(W)$.   $W$ is called an \emph{isolating neighborhood} of $S$.

Suppose $W\subset X$ is compact.  Then  the subset $A\subset W$ is  said to be \emph{positively invariant} with respect to $W$ if
\[
x\in A, \;\;t\geq 0,\;\; [0,t]\cdot x\subset W\Longrightarrow [0,t]\cdot x\subset A.
\]
Suppose $W$ is an isolating neighborhood of $S$. An \emph{index pair in $W$} for  the isolated invariant set $S$ is  a pair of compact sets $(N, N^-)$, $N^-\subset N\subset W$,  with the following properties.

\begin{enumerate}

\item[($I_0$)]  $N$ is positively invariant in $W$.

\item[($I_1$)] $N\setminus N^-$ is a neighborhood of $S$, and $S =I(\,\cl(N\setminus N^-)\,)$.

\item[($I_2$)] $N^-$ is positively invariant in $N$.

\item[($I_3$)]  If  $x\in N$, and $[0,\infty)\cdot x\not\subset N$, then there exists $t\geq 0$  such that $[0,t]\cdot x\subset N$, and $t\cdot x\in N^-$.
\end{enumerate}

A pair of compact sets $(N,N^-)$ satisfying the conditions ${\bf I}_1$, ${\bf I}_2$ and ${\bf I}_3$ will be called an \emph{index pair} of $S$.

\begin{theorem}[Existence of index pairs]  Suppose $S$ is an isolated invariant set of the  flow  $\Phi$, $W$ is an isolated neighborhood of $S$ and $U$ is a neighborhood of $S$. Then there exists an index pair $(N_U, N_U^-)$ of $S$ in $W$ such that
\[
\cl(N_U\setminus N_U^-)\subset U.\proofend
\]
\end{theorem}

 Suppose  $(N,N^-)$ is an index pair for the isolated  invariant set  $S$. For every $t\geq 0$ we define
 \[
 {}^tN:=\bigl\{ \,\ x\in N;\;\; [0,t]\cdot x\subset N\,\bigr\},
 \]
 \[
  {}^{-t}N^-:=\bigl\{ x\in  N;\;\;\exists  t'\in [0,t]: [0,t']\cdot x\subset N,\;\; t'\cdot x\in N^-\,\bigr\}
 \]
Then  ${}^tN\subset N$,  ${}^{-t}N^-\supset N^-$.  The inclusion induced map
\[
i_t:{}^tN/\,{}^tN\cap N^-\ra N/N^-
\]
is a homotopy equivalence with  homotopy inverse $f_t: N/N^-\ra {}^tN/\, {}^tN\cap N^-$ given by (see \cite[III.4.2]{Con})
\[
f_t([x])= \begin{cases}
[t\cdot x] & [0,t]\cdot x\subset N\\
[{}^tN\cap N^-] & {\rm otherwise}.
\end{cases}
\]
Similarly, the inclusion induced map
\[
j_t: N/N^-\ra N/{}^{-t}N^-
\]
is a homotopy equivalence with homotopy inverse   $g_t$ given as the composition $i_t\circ h_t$, where $h_t: N/{}^{-t}N^- \ra {}^tN/{}^tN\cap N^-$ is the \emph{homeomorphism} given by
\[
h_t([x]) =\begin{cases}
[t\cdot x] & [0,t]\subset N\setminus N^-\\
[{}^tN\cap N^-] & \mbox{otherwise}.
\end{cases}
\]
Suppose $(N_0, N_0^-)$ and $(N_1, N_1^-)$ are two index pairs  in $W$ for $S$.  Then   there exists $T=T(N_0,N_1)>0$  such that for any  $t>T$  we have
\[
({}^tN_0,{}^tN_0\cap N_0^-) \subset (N_1, {}^{-t} N_1^-),\;\; ({}^tN_1,{}^tN_1\cap N_1^-) \subset (N_0,  {}^{-t}N_0^-).
\]
Fix  $t>T(N_1,N_0)$, denote by $\alpha_t$ the inclusion induced map
\[
\alpha_t:{}^tN_0/{}^tN_0\cap N_0^-\ra N_1/{}^{-t} N_1^-,
\]
and  by $\beta_t$ the inclusion induced map
\[
\beta_t: {}^tN_1/{}^tN_1\cap N_1^- \ra  N_0/{}^{-t}N_0^-.
\]
Define $\eC^t_{N_1,N_0}: N_0/N_0^-\ra N_1/N_1^-$ as the composition
\[
N_0/N_0^-\stackrel{f_t^0}{\Lra} {}^tN_0/ {}^tN_0\cap N_0^- \stackrel{\alpha_t}{\Lra}N_1/{}^{-t}N_1^-\stackrel{g^1_t}{\Lra} N_1/N_1^-.
\]
For any $t,t'>T(N_1,T_0)$, the maps $\eC^t_{N_1,N_0}$ and $\eC^{t'}_{N_1,N_0}$ are homotopic.   We denote by
\[
\eC_{N_1,N_0}\in [N_0/N_0^-, N_1/N_1^-],
\]
 the homotopy class  determined by  this family of maps,  and we will refer to it as the connector  from $N_0$ to $N_1$

If $(N_0,N_0^-)$, $(N_1, N_1^-)$ and $(N_2,N_2^-)$ are three index pairs, and
\[
t>\max\bigl\{\,\; T(N_2,N_1),\; T(N_1,N_0),\;\; T(N_2,N_0)\,\bigr\},
\]
then we have a homotopy
\[
\eC^t_{N_2, N_0}\simeq \eC^t_{N_2,N_1}\circ \eC^t_{N_1,N_0}.
\]
In particular, if  $N_2=N_0$ we deduce
\[
\eC^t_{N_0,N_1}\circ \eC^t_{N_1,N_0} \simeq \eC^t_{N_0,N_0} \simeq \one,
\]
so that all the connectors are homotopy equivalences.

The homotopy type of the pointed space $[N/N^-]$  is therefore independent of the    index pair $(N, N^-)$ of $S$. It is called the \emph{Conley index} of $S$ and  it is denote by $\bh(S)$, or $\bh(S,\Phi)$.

Consider now a compact tame set $X$ embedded in some Euclidean space $E$. Denote  by $|\bullet|$ the Euclidean norm on $E$. Suppose   $\Phi$ is a tame flow on $X$.

\begin{definition} A stationary point $p$ of $\Phi$ is called \emph{Morse like}   if  there exists    tame continuous function $f:X\ra \bR$ with the  following properties.

\begin{itemize}

\item $f(p)=0$.

\item The function $f$ decreases, along the trajectories of the flow, not necessarily strictly.

\item There exists $c_0>0$ such that
\[
\Cr_\Phi\cap \{0<|f|<c_0\}=\emptyset.
\]
\item  The set $\Cr^0:=\Cr_\Phi\cap\{f=0\}$ is finite.

\item  The function $f$ decreases \emph{strictly} along any  portion of nonconstant trajectory situated in $\{|f|<c_0\}$

\end{itemize}

The function $f$  is called   a \emph{local Lyapunov} function adapted to the stationary point $p$.\qed

\end{definition}

Suppose $p\in X$ is a Morse like stationary point of the flow, and $f$ is a local Lyapunov function adapted to  $p$. For every $c\in \bR$ we denote by $X_c$ the level set $\{f=c\}$. Denote by $W^+_p$  and respectively $W^-_p$ the stable and respectively unstable varieties  of the point $p$, and set
\[
\eL^-_p(\ve):= W^-_p\cap X_{-\ve},\;\; \eL^+_p(\ve):= W^+_p\cap X_\ve.
\]

\begin{lemma} Suppose  $\ve \in (0,c_0)$. Then the following hold.

\noindent (a) The link $\eL^\pm_p(\ve)$ is a compact subset of  $X_{\pm \ve}$.

\noindent (b) The tame set $W^\pm_p(\ve)=W^\pm_p\cap \{|f|\leq \ve\}$ is tamely homeomorphic to  a cone on $\eL^+_p(\ve)$.
\label{lemma: link}
\end{lemma}

\proof  (a) We prove only the case $\eL^+_p(\ve)$  since the other case is obtained from this by time reversal.     We argue by contradiction. Suppose
\[
x_0\in \cl(\,\eL^+_p(\ve)\,)\setminus \eL^+_p(\ve)
\]
Then there exists a tame continuous path $(0,1]\ni s\mapsto x_s\in\eL^+_p(\ve)$ such that
\[
\lim_{s\ra 0^+}x_s=x_0.
\]
Since $f(t\cdot x_s)\in [0,\ve]$, $\forall s,t>0$ we deduce $f(t\cdot x_0)\in[0,\ve]$, $\forall t\geq 0$. If we set
\[
q=\lim_{t\ra \infty} t\cdot x_0
\]
we deduce that $q$ is a stationary point of $\Phi$ such that $f(q)\in [0,\ve]$. Since $\ve<c_0$ we deduce $q\in \Cr^0$, and since $x_0\not\in W^+_p$ we deduce $q\neq p$.

\begin{figure}[h]
\centerline{\epsfig{figure=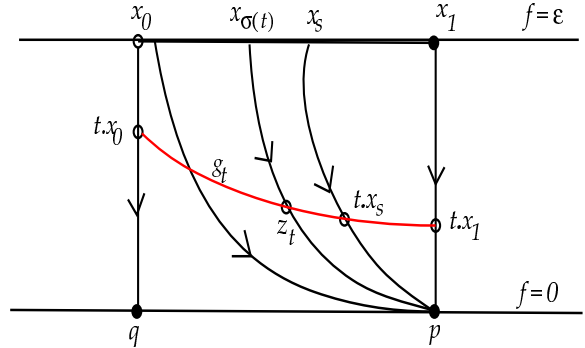,height=1.8in,width=3in}}
\caption{\sl The stable variety of $p$ is arbitrarily close to that of $q$.}
\label{fig: tame10}
\end{figure}
Consider the family of paths (see Figure \ref{fig: tame10})
\[
g_t: [0,1]\ra X,\;\;  g_t(s)= tx_s.
\]
Let
\[
\delta:=\min\{ |q'-q''|;\;\;q',q''\in \Cr^0,\;\;q'\neq q''\},\;\;d_t:= |t\cdot x_0-t\cdot x_1|,
\]
and   consider the definable family of closed subsets of the unit interval
\[
I_t:=\bigl\{\, s\in [0,1];\;\; |t\cdot x_s- t\cdot x_0|=\frac{1}{2}\min(\delta, d_t)\,\bigr\}.
\]
Note that $I_t\neq \emptyset$,   $\forall t>0$.    We can then find a definable function
\[
\si: [0,\infty)\ra [0,1]
\]
such that $\si(t)\in I_t$,  $\forall t>0$.  Set $z_t:= t\cdot x_{\si(t)}$  so that
\[
|z_t-t\cdot x_0|=\frac{1}{2}\min\{\delta, d_t\},\;\;\forall t>0.
\]
 The function $\si$ is continuous for $t$ sufficiently large  and   the limit
\[
\si_\infty:=\lim_{t\ra \infty}\si(t)
\]
exists and it is finite. Observe that the definable path
\[
t\mapsto t\cdot x_{\si(t)}\in \{0\leq f\leq \ve\},
\]
has a limit as $t\ra \infty$ which we denote by $z_\infty$. Since $d_t\ra |q-p|\geq \delta$ we deduce
\[
|z_\infty-q|=\frac{1}{2}\delta.
\]
In particular, we deduce that $z_\infty$ is not a stationary point of the flow.

Consider now the function
\[
e: X\ra (-\infty,0] ,\;\;e(x)=f(x)-f(\Phi_1(x) ),
\]
where $\Phi_1$ denotes the time-$1$ map  determined by the flow $\Phi$. Since $z_\infty$ is not a  stationary  point  we deduce
\[
e(z_\infty) <0.
\]
Because the  time-$1$ map $\Phi_1$  is continuous, we deduce that, for every positive $\hbar$  such that $\hbar\leq |e(z_\infty)|$, there exists an open neighborhood  $U_\hbar$ of $z_\infty$ in $X$ such that
\[
e(x)<\hbar,\;\;\forall z\in U_\hbar.
\]
In particular, for sufficiently large  $t$, we have $z_t\in U_\hbar$, and thus
\[
0\leq f(\,\Phi_1(z_t)\,) < f(z_t)-\hbar.
\]
If we let $t\ra \infty$ we deduce
\[
0\leq f(\,\Phi_1(z_\infty)\,) \leq  f(z_\infty)-\hbar=-\hbar.
\]
This contradiction proves the compactness of $\eL^+_p(\ve)$.

\noindent (b) From part (a) we deduce easily that  $W^+_p(\ve)$ is compact. Consider the tame homeomorphism
\[
[0,1)\ni t\mapsto  t(s)=\frac{s}{1-s}\in[0,\infty).
\]
Now consider the map
\[
[0,1]\times \eL^+_p(\ve)\ra W^+_p(\ve),\;\; (s,x)\mapsto t(s)\cdot x.
\]
This maps the slice $\{1\}\times \eL^+_p(\ve)$ to $p$ and it induces a   tame continuous  bijection  from the cone on $\eL^+_p(\ve)$ to $W^+_p(\ve)$. Since $W^+_p(\ve)$ is compact we deduce that this map is a homeomorphism. \qed

The (tame) topological type of $\eL^+_p(\ve)$ and respectively $\eL^-_p(\ve)$ is independent of $\ve$ if $\ve$ is sufficiently small because the tame continuous map
\[
f:W^\pm_p(\ve)\setminus \{p\}\ra \bR
\]
is locally trivial for $\ve>0$. We will refer  to this tame  homeomorphism class  as the \emph{stable} and respectively \emph{unstable link} of $p$, and we will denote it by $\eL^\pm_p$.

Observe that for $\ve>0$ sufficiently small  the tame set  $W^\pm_p\cap \{|f|\leq \ve\}$ is tamely homeomorphic to the cone on $\eL^\pm_p$, and that the links  $\eL^\pm_q(\ve)$, $q\in \Cr^0$ are mutually disjoint compact subsets of $X_{\pm\ve}$.

\begin{proposition}  Let $\ve \in (0, c_0)$ and  let $K\subset U_\ve$ be a tame compact  neighborhood  of $\eL^-_p(\ve)$ in the level set $X_{-\ve}$  such that
\[
K\cap W^-_q=\emptyset,\;\;\forall q\in\Cr^0,\;\;q\neq p,
\]
and set
\[
N= N_{\ve, K}:= \bigl(\, W^-_p\cup W^+_p \cup (-\infty, 0]\cdot K\,\bigr)\cap \{|f\leq \ve\}.
\]
Then the pair $(N,K)$ is an index pair for $p$.
\label{prop: ipair}
\end{proposition}

\proof   The conditions $I_2$ and $I_3$ in the definition of an index pair are clearly satisfied  due to the existence of the Lyapunov function $f$, so it suffices to show that $N$ is a \emph{compact, isolating, neighborhood} of $p$.   In the proof we will need   several auxiliary results.

\begin{lemma} Suppose
\[
(0,1]\ni s\mapsto x_s \in X_{-\ve},\;\; (0,1)\ni s\mapsto  t_s\in (0,\infty)
\]
are tame continuous  paths   such that
\[
\lim_{s\ra 0^+} t_s=\infty\;\;\mbox{and}\;\;f(\,(-t_s)\cdot x_s)\leq 0,\;\;\forall s\in(0,1).
\]
Then there exists  $q\in\Cr^0$  such that $x_0\in \eL^-_q(\ve)$ and $\lim_{s\ra 0^+} (-t_s)\cdot x_s=q$.
\label{lemma: compN}
\end{lemma}

\proof    Observe that
\[
(-T)\cdot x_0\in \{-\ve \leq f\leq 0\},\;\;\forall T>0
\]
so that there exists $q\in\Cr^0$ such that $x_0\in \eL^-_q(\ve)$. Set $z_s=(-t_s)\cdot x_s$.

The definable path $s\mapsto z_s$ has a limit   $z_0=\lim_{s\ra 0+} z_s$.  Since
\[
T\cdot z_0\in \{-\ve\leq  f\leq 0\},\;\;\forall T>0,
\]
the point  $z_0$ must be a stationary point.  We claim $z_0=q$. We argue by    contradiction, so we assume  $z_0\neq q$.

  Set $y_s:= (-t(s))\cdot x_0$. For every $s\in (0,1]$ consider the definable continuous path
\[
g_s:[0,1]\ra X,\;\; g_s(\lambda)= (-t(s))\cdot x_{\lambda\cdot s}.
\]
Observe that $g_s(0)= y_s$ and  $g_s(1)=  z_s$.  Arguing as in  the proof of Lemma \ref{lemma: link} we can find a definable function
\[
(0,1)\ni s\mapsto \lambda_s\in [0,1]
\]
such that
\[
\dist\bigl(\, g_s(\lambda_s), q) =\frac{1}{2}\min \{ \delta, |z_s-y_s|\,\},\;\;\delta:=\min\{\, |q'-q''|;\;\;q',q''\in\Cr^0,\;\;q'\neq q''\,\}.
\]
We set
\[
\gamma_s:=g_s(\lambda_s) =(-t(s)) \cdot x_{\lambda_s s}.
\]
 Then,    as $s\searrow 0$, the point $\gamma_s$ converges to a point $\gamma_0$ such that
\[
\gamma_0\in \{ -\ve \leq f\leq 0\},\;\;\dist(\gamma_0, q)=\frac{1}{2}\min\{\delta, |z_0-q|\}=\frac{1}{2}\delta.
\]
Thus $\gamma_0$ \emph{is not a stationary point} of  $\Phi$.    We claim that
\begin{equation}
f(T\cdot \gamma_0) \geq -\ve, \;\; \forall T>0.
\label{eq: long1}
\end{equation}
Indeed, for every $T>0$, and for every $\hbar>0$ there exists a small neighborhood $U=U_{T,\hbar}$ of $\gamma_0$ such that for every $x\in U$ we have
\[
|f(T\cdot x)- f(T\cdot \gamma_0)|<\hbar.
\]
We can now find $s>0$ such that $\gamma_s\in U_{T,\hbar}$ and  $t(s)> T$,  from which we deduce
\[
f(T\cdot \gamma_0)\geq  f(T\cdot \gamma_s) -\hbar \geq  f(t(s)\cdot \gamma_s)-\hbar =f(x_{\lambda_s\cdot s})-\hbar=-\ve-\hbar.
\]
This proves the claim (\ref{eq: long1}) which in turn implies that  $\gamma_0$ \emph{has to be a stationary point}. This contradiction completes the proof of  Lemma \ref{lemma: compN}. \qed

Observe that for every $x\in X_{-\ve}$ we have
\[
\Phi_{-\infty} (x)\in \{ f\geq 0\}.
\]
Define  $T=T_{-\ve}: X_{-\ve}\ra [0,\infty]$ by setting   $T(x)=\infty$ if $\Phi_{-\infty} x\in\Cr^0$, and otherwise, we let $T(x)$ to be the unique positive real number  such that
\[
(-T(x))\cdot x\in X_0.
\]
 Using the definable  homeomorphism
\[
\si:[0,\infty)\ra [0,1),\;\; t\mapsto \si(t)=\frac{t}{1+t}.
\]
we obtain a compactification  $[0,\infty]$ of $[0,\infty)$ tamely homeomorphic to $[0,1]$.

\begin{lemma}[Deformation Lemma] (a) The  tame function
\[
X_{-\ve} \ni x\mapsto  T_{-\ve}(x)\in [0,\infty]
\]
is continuous.

\noindent (b) The tame function
\[
\eD^{-\ve}_\Phi: \bigl\{ (x,t)\in X_{-\ve}\times [0,\infty];\;\; t\leq T_{-\ve}(x) \,\bigr\} \ra \{-\ve\leq f\leq 0\},\;\;(x,t)\mapsto  (-t)\cdot x
\]
is continuous.
\label{lemma: cont}
\end{lemma}

\proof  For simplicity, during this proof, we will write $T(x)$ instead of $T_{-\ve}(x)$.

\noindent (a) By invoking  the closed graph theorem  it suffices to show  that for  any  continuous definable path
\[
(0,1)\ni s\mapsto (x_s, T(x_s))\in X_{-\ve}\times [0,\infty]
\]
such that
\[
x_s\ra x_0,\;\; T(x_s) \ra T_0\in [0,\infty],
\]
then $T_0= T(x_0)$. Observe that if $T(x_s)=\infty$, for all $s$ sufficiently small, then there exists  $q\in \Cr^0$ such that $x_s\in \eL^-_q((\ve)$, and since $\eL^-_q(\ve)$ is compact, we deduce $x_0\in \eL^-_q(\ve)$. Thus, we can assume that $T(x_s)<\infty$, for all $s$.

If $T_0 <\infty$, the conclusion follows from the  continuity  of the flow. Thus, we can assume $T_0=\infty$,  and $T(x_s)\nearrow \infty$ as $s\searrow 0$, and we have to prove that  there exists  $q\in  \Cr^0$ such that $x_0\in \eL^-_q(\ve)$. This follows immediately from the fact that
\[
(-T)\cdot x_0 \in \{ -\ve\leq f \leq 0\},\;\;\forall T>0,
\]
 so that   $x_0$ must belong to the unstable variety of a stationary point situated in the   region  $\{-\ve \leq f\leq 0\}$.

 \noindent (b)   Again we rely on the closed graph theorem. We have to show that for every  tame continuous paths
 \[
 (0,1)\ni s\mapsto (x_s, t_s)\in X_{-\ve}\times [0,\infty],
 \]
 such that
 \[
 0\leq t_s\leq T(x_s),\;\;\lim_{s\ra 0^+} x_s= x_0,\;\;\lim_{s\ra 0^+}t_s=t_0,\;\;  \lim_{s\ra 0^+} (-t_s)\cdot x_s =y_0\in X_0,
 \]
we have $y_0=(-t_0)\cdot x_0$.

Arguing as in (a),  we  see that the only nontrivial situation is when $t_s\nearrow\infty$ as $s\searrow 0$. In this case, we have to prove that $y_0\in \Cr^0$  and $x_0\in\eL^-_{y_0}(\ve)$. This follows from Lemma \ref{lemma: compN}. \qed

The Deformation Lemma  has many useful corollaries.

\begin{corollary} The continuous tame map
\[
\eT^{-\ve}_\Phi: X_{-\ve}\ra X_0,\;\; x \ra \eT^{-\ve}_\Phi(x) = (-T_{-\ve}(x))\cdot x\in X_0
\]
induces a tame homeomorphism
\[
X_\ve^*= X_{-\ve}\setminus \bigcup_{q\in\Cr^0}\eL^-_q(\ve)\ra X_0^*= X_0\setminus \Cr^0.
\]
\label{cor: homeo}
\end{corollary}

\proof The map  $\eT^{-\ve}_\Phi: X_{-\ve}^*\ra X_0^*$ is  continuous and bijective. Its inverse is continuous because its graph is closed. \qed

 Consider the strip
\[
\eS_{-\ve}: =\bigl\{ (x,s)\in X_{-\ve}\times [0,1];\;\;   s\leq \si_{-\ve}(x)=\frac{T_{-\ve}(x)}{1+T_{-\ve}(x)},\;\;\forall x\in X_{-\ve}\,\bigr\}.
\]
Observe that we have a tame homeomorphism
\[
\eA_{-\ve}: X_{-\ve}\times [0,1]\ra \eS_{-\ve},\;\; (x, \lambda)\mapsto(x, \si_{-\ve}(x)\cdot \lambda)
\]
and a tame homeomorphism
\[
S: \eS_{-\ve} \ra \bigl\{ (x,t)\in X_{-\ve}\times [0,\infty];\;\; t\leq T_{-\ve}(x) \,\bigr\} ,
\]
given by
\[
S(x,s)\mapsto (x, \frac{s}{1-s}).
\]
The composition
\[
\hat{\eD}^{-\ve}_\Phi:=\eD_\Phi^{-\ve}\circ S\circ\eA_{-\ve}:  X_{-\ve}\times [0,1]\ra   \{-\ve\leq f\leq 0\}
\]
is a tame continuous map, which  along  $X_{-\ve}\times \{1\}$  it coincides with the map $\eT^{-\ve}_\Phi:X_{-\ve}\ra X_0$.

The  natural deformation retraction of $X_{-\ve}\times [0,1]$ onto $X_{-\ve}\times \{1\}$  determines a deformation retraction of
\[
\eR_\Phi^{-\ve}: \{-\ve\leq f\leq 0\}\times [0,1]\ra \{-\ve\leq f\leq 0\}
\]
of $\{-\ve\leq f\leq 0\}$  onto $\{f=0\}$. The next result summarizes  the above observations.

\begin{corollary} The deformation $\hat{\eD}^{-\ve}_\Phi$ induces a homeomorphism between the mapping cylinder of $\eT^{-\ve}_\Phi:X_{-\ve}\ra X_0$ and the region $\{-\ve\leq f\leq 0\}$.\qed
\label{cor: map-cyl}
\end{corollary}

\begin{remark}  The maps  $T_{-\ve}$, $\eD^{-\ve}_\Phi$, $\eT^{-\ve}_\Phi$, $\eR^{-\ve}_\Phi$   have ``positive'' counterparts $T_\ve(x): X_\ve\ra [0,\infty)$,
\[
 \eD^\ve_\Phi: \bigl\{ (x,t)\in X_{\ve}\times [0,\infty];\;\; t\leq T_{\ve}(x) \,\bigr\} \ra \{\,0\leq f\leq -\ve\,\},
 \]
 \[
 \eT^\ve_\Phi: X_\ve\ra X_0,
 \]
 and
 \[
 \eR_\Phi^{\ve}: \{\,\ve\geq f\geq \ve\,\}\times [0,1]\ra \{\, \ve \geq f\geq 0\},
 \]
 and their similar properties follow by time reversal     from their ``negative'' counterparts. \qed

\label{rem: def}
\end{remark}

Now set
\[
K_0:=\eT^{-\ve}_\Phi(K)\subset X_0,\;\;K^+=(\eT^\ve_\Phi)^{-1}(K_0)\subset X_\ve.
\]
Then $K_0$ is a compact neighborhood  of $p$ in $X_0$, $K^+$ is a compact neighborhood of $\eL^+_p(\ve)$ in $X_\ve$,  and we have the equality
\begin{equation}
N= \underbrace{\hat{\eD}^{-\ve}_\Phi\bigl(\, K\times [0,1]\,\bigr)}_{N_{\leq 0}}\, \cup \, \underbrace{\hat{\eD}^\ve_\Phi\bigl(\, K^+\times [0,1]\,\bigr)}_{N_{\geq 0}}.
\label{eq: N}
\end{equation}
Now observe that  $N_{\leq 0}$ is a \emph{compact} neighborhood of $p$ in $\{f\leq 0\}$ and $N_{\geq 0}$ is a \emph{compact} neighborhood of $p$ in $\{f\geq 0\}$.

The fact that $N$ is an isolating neighborhood of $p$ follows from  (\ref{eq: N}). This completes the proof of Proposition \ref{prop: ipair}. \qed

\begin{theorem} Suppose $\Phi$ is a tame flow on a tame compact set $X$,  $p$ is a Morse like stationary point of $\Phi$, and $f$ is a local Lyapunov function adapted to $p$.   We denote  by $W^-_p$ the unstable variety of $p$, and for every $\ve>0$ we set
\[
W^-_p(\ve):= W^-_p\cap \bigl\{\, -\ve\leq f\leq 0\,\bigr\},\;\;\eL^-_p(\ve):=W^-_p\cap\bigl\{\, f= -\ve\,\bigr\}.
\]
Then  the Conley index $h_\Phi(p):=h(\{p\}, \Phi)$ is homotopy equivalent to the pointed space $W^-_p(\ve)/\eL^-_p(\ve)$, for all sufficiently small $\ve>0$.
\label{th: conley}
\end{theorem}

\proof   We continue to use the same notations as in the proof of Proposition \ref{prop: ipair}.   Fix $\ve>0$ sufficiently small so that   the only stationary points of $\Phi$ in $\bigl\{|f|\leq \ve\bigr\}$ lie on the level set $X_0$.

Because both $X_{-\ve}$ and $\eL^-_p(\ve)$ are tame  compact tame sets we can find a triangulation  of $X_{-\ve}$ so that $\eL^{-\ve}_p$ is a    subcomplex of the triangulation of $X_{-\ve}$.        From the classical results of J.H.C. Whitehead \cite{White} we deduce that  for any neighborhood $U$  of $\eL^-_p(\ve)$ we can find triangulations of the pair $(X_{-\ve},\eL^-_p(\ve))$ such that the simplicial neighborhood of $\eL^-_p(\ve)$ in $X_{-\ve}$ is contained in $U$ and collapses onto $\eL^-_p(\ve)$.  Fix such a simplicial  neighborhood $K$ which is disjoint from $W^-_q$,  $\forall q\in \Cr^0_\Phi$, $q\neq p$.  Because $K$ collapses onto $\eL^-_p(\ve)$ we can find a \emph{tame} deformation retraction onto   $\eL^-_p(\ve)$.

Form the index pair $(N, K)= (N_{\ve,K}, K)$.  Let us point out that both $N$ and $K$ are  compact sets, and in particular the inclusion $K\hra N$ is a cofibration.

 Using the  deformation retraction $\hat{\eD}^\ve_\Phi$ we see that the pair $(N,K)$ is homotopy equivalent to the pair $(N_{\leq 0},K)$, $N_{\leq 0}=N\cap\{f\leq 0\}$. Corollary \ref{cor: map-cyl} implies that $N_{\leq 0}$ is  homeomorphic to the mapping cylinder  of the  tame map
 \[
 \eT_\Phi^{-\ve}: K\ra  K_0=\eT^{-\ve}_\Phi(K)\subset X_0.
 \]
Corollary \ref{cor: homeo} shows that $\eT^{-\ve}_\Phi$ induces a tame homeomorphism
\[
K^*=K\setminus \eL^-_p(\ve)\ra K^*_0=K_0\setminus \{p\}.
\]
Now observe that $W^-_p(\ve)$ is also homeomorphic to the mapping cylinder of the map $\eT^{-\ve}_\Phi:\eL^-_p(\ve)\ra \{p\}$.  We deduce that   $N_{\leq 0}$ is homeomorphic to the mapping cylinder of the natural projection
\[
\pi: K \ra K/\eL^-_p(\ve).
\]
A tame deformation retraction  of $K$ onto $\eL^-_p(\ve)$ extends to  a deformation    of the mapping cylinder of $\pi$ to the mapping cylinder of $\pi|_{\eL^-_p(\ve)}$ which is homeomorphic to $W^-_p(\ve)$.  \qed

\begin{remark} The Conley index computation   in this section  bares a striking resemblance  with the  computation of  Morse data in the Goresky-MacPherson stratified Morse theory, \cite{GM}.   We believe  this resemblance goes beyond  the level of accidental coincidence, but we will pursue this line of thought elsewhere.  \qed
\end{remark}

Here is a simple application of the above result.  Suppose $X$ is a tame space, and $\Phi$ is a  \emph{Morse like} tame flow on $X$. This means  that $\Phi$ has finitely many stationary points, and admits a tame Lyapunov function $f$.  Observe that the local minima are stationary points  of $\Phi$. They are characterized by the condition $W^-_p=\{p\}$. We denote by $\Cr_\Phi\subset X$ the set of stationary points.

For every  compact  tame space $Y\neq \emptyset$ we denote by $\eP_Y(t)$ the  Poincar\'e polynomial of $Y$
\[
\eP_Y(t)=\sum_{k\geq 0} (\dim H_k(Y,\bR))t^k\in \bZ[t].
\]
If $A$ is a compact tame subset  of $A$  we denote by $\eP_{Y,A}(t)\in\bZ[t]$ the Poincar\'e polynomial of the pair $(Y,A)$ defined in a similar  fashion.  In particular, for  every  $p\in \Cr_\Phi$,  we denote by  $\bM_{\Phi,p}(t)$ the Poincar\'e polynomial of the Conley index of $p$, and we will refer to it as the \emph{Morse polynomial of the stationary point $p$}.
 As in \cite{N2}, we define an order relation $\succ$ on $\bZ[t]$ by declaring $A\succeq B$ if there exists a polynomial $Q\in\bZ[t]$ with nonnegative coefficients such that
\[
 A(t)= B(t)+(1+t)Q(t).
\]
\begin{corollary}[Morse inequalities] Let $\Phi$ be a Morse like flow on $X$ with Lyapunov function  $f$ be as above.
\[
\sum_{p\in\Cr_\Phi}  \bM_{\Phi,p}(t)\succeq \eP_X(t).
\]
\label{cor: morse-i}
\end{corollary}

\proof Define the discriminant set,
\[
\Delta_f: =f(\Cr_\Phi).
\]
$\Delta_\Phi$ is a finite set of real numbers
\[
\Delta_\Phi=\bigl\{ c_0<c_1<\cdots <c_n\,\bigr\}.
\]
For $k=0,\dotsc, n$  we set
\[
\Cr^k_\Phi:=\Cr_\Phi\cap \{f=c_k\}.
\]
Now choose $r_{0}=c_0< r_1< c_1\cdots <c_{n-1}< r_n<c_n=r_{n+1}$ and set $X^k:= \{f\leq r_k\}$.  For each   $k=0,1,\dotsc, n$ the pair $[X^{k+1}, X^k]$ is an  index pair for the isolated invariant set $\Cr^k_\Phi$.  We  deduce that
\[
h_\Phi(\Cr^k_\Phi) = \bigvee_{p\in \Cr^k_\Phi}  h_\Phi(p).
\]
Hence
\[
\eP_{X^{k+1}, X^k}(t)=\sum_{p\in \Cr^k_\Phi}  \bM_{\Phi,p}(t).
\]
Using \cite[Remark 2.16]{N2} we deduce
\[
\sum_k\eP_{X^{k+1}, X^k}(t) \succeq \eP_X(t)
\]
from which  the Morse inequality follow immediately. \qed

\begin{remark} Let  $p\in \Cr^k_\Phi$. If $\eL^-_p(\ve)=\emptyset$  then $\bM_{\Phi, p}(t)=1$.  Otherwise
\[
h_\Phi(p)\simeq (C\eL^-_p(\ve),\eL^-_p(\ve)\,),
\]
where $CA$ denotes the  cone on the topological space $A$.  From the long exact sequence of  the pair  $(C\eL^-_p(\ve),\eL^-_p(\ve)\,)$ we deduce that if $\eL^-_P(\ve)\neq \emptyset$ then
\[
\dim H_0(C\eL^-_p(\ve),\eL^-_p(\ve)\,)=0, \;\;\dim H_1(C\eL^-_p(\ve),\eL^-_p(\ve)\,)= \dim H_{0}(\eL^-_p(\ve)\,)-1
\]
\[
\dim H_{k+1}(C\eL^-_p(\ve),\eL^-_p(\ve)\,)=\dim H_{k}(\eL^-_p(\ve)\,),\;\;\forall  k> 0.
\]
If we denote by  $\tilde{\bM}_{\Phi,p}(t)$ the Poincar\'e polynomial of the reduced homology of $\eL^-_p(\ve)$ we deduce
\[
\bM_{\Phi,p}(t)= t \tilde{\bM}_{\Phi,p}(t).
\]
If we define for uniformity
\[
\tilde{\bM}_{\Phi,p}(t)=t^{-1},\;\;{\rm if}\;\;\eL^-_p(\ve)=\emptyset
\]
then the previous equality holds in all the cases. We can rephrase the  Morse inequalities as
\begin{equation}
\sum_{p\in\Cr_\Phi}  t\tilde{\bM}_{\Phi,p}(t)\succeq \eP_X(t).
\label{eq: morse-ineq}
\end{equation}
\qed
\end{remark}

%% file: tameflow10.tex
The results in the previous section allows us to give a more detailed picture of the  gradient like tame flows on compact tame sets.   For any compact tame space $X$ we denote by $\Flg(X)$  the set of gradient-like tame flows on $X$ with finitely many stationary points.

\begin{definition} (a) A \emph{tame blowdown} is a   continuous tame map $\beta:Y\ra X$,  such that $X$ and $Y$ are compact tame sets, and there exists  finite nonempty subset $L=L_\beta\subset X$ such that the induced map
\[
\beta: Y\setminus \beta^{-1}(L_\beta)\ra  X\setminus L_\beta
\]
is a homeomorphism. The compact set $L_\beta$ is called the blowup locus of $\beta$. The compact set $\beta^{-1}(L_\beta)$ is called the \emph{exceptional locus} of $\beta$ and it is denoted by $E_\beta$. We will also say that $Y$ is a \emph{tame blowup} of $X$.  A \emph{weight} for the blowdown map $\beta$ is a tame continuous function $w: X\ra [0,\infty)$ such that $w^{-1}(0)=L_\beta$.  We will refer to  a pair (blowdown, weight) as a \em{weighted blowdown}.

 (b) A  \emph{tame flop} is a diagram of the form
 \[
 \begin{diagram}
 \node{}\node{Y}\arrow{sw,t}{\beta_-}\arrow{se,t}{\beta_+}\node{}\\
 \node{X_-}\node{}\node{X_+}
 \end{diagram}
 \]
 where $\beta_\pm :Y\ra X_\pm$ are tame blowdowns.   The \emph{connector} associated  to the flop is obtained by gluing the mapping cylinder of $\beta_-$ to the mapping cylinder of $\beta_+$ along $Y$ using the identity map $\one_Y$. We will denote the connector by $(\stackrel{\beta_-}{\longleftarrow}Y\stackrel{\beta_+}{\Lra})$.

\noindent(c)  A \emph{tame flip}\footnote{The  ``o'' in flop  indicates  that the arrows  arrow coming \emph{out} of the middle of the diagram, while the ``i'' in flip  indicate that the arrows are coming \emph{into} the middle  of the diagram.} is a diagram
\[
\begin{diagram}
\node{Y_-}\arrow{se,b}{\beta_-}\node{}\node{Y_+}\arrow{sw,b}{\beta_+}\\
\node{}\node{X}\node{}
\end{diagram}
\]
where $\beta_\pm:Y_\pm X$ are blowdown maps. The  \emph{connector} of the flip is  the tame space obtained by gluing the mapping cylinder  of $\beta_-$ to the mapping cylinder of $\beta_+$ along $X$ using the  identity map $\one_X$. We will denote it by $(\stackrel{\beta_-}{\Lra} X\stackrel{\beta_+}{\longleftarrow})$.\qed
\end{definition}

\begin{remark} In the above definition  of a blowdown map  $\beta: Y\ra X$ we allow for the possibility that  the exceptional locus $E_\beta$ is empty. For example, the map
\[
\{0\}\ra \{0,1\},\;\;0\mapsto 0,
\]
is a blowdown map, with blowup locus $\{1\}$, and  empty exceptional locus.\qed
\end{remark}

Suppose that $\Phi$ is a gradient-like tame flow on  a compact  tame space $X$, and that the set $\Cr_\Phi$ of stationary points is finite. Fix a Lyapunov function $f:X\ra \bR$, and let
\[
\{ c_0<c_1<\cdots < c_\nu\}
\]
be the set $f(\Cr_\Phi)$.  For $i=1,\dotsc,\nu$ we set $ d_i:=\frac{c_{i-1}+c_i}{2}$, and we define
\[
Y_i=\{f=d_i\},
\]
and
\[
\Cr^j_\Phi:=\Bigl\{x\in \Cr_\Phi;\;\;f(x)=c_j\,\},\;\;X_j:=\{f=c_j\,\bigr\},\;\;\forall j=0,\dotsc, \nu.
\]
For every point $x\in X$  we  denote  by $\Phi(x)$ the trajectory of $\Phi$ through $x$,
\[
\Phi(x):=\bigl\{\Phi_t(x);\;\;t\in \bR\,\bigr\}.
\]
In the previous section we have proved that the flow defines   tame blowdowns
\[
\lambda_{i}=\lambda_{i}^\Phi: Y_{i+1}\ra X_{i},\;\;  \lambda_{i}(y) =\Phi(y)\cap X_{i},
\]
\[
\rho_i=\rho_i^\Phi: Y_i\ra X_i,,\;\;\rho_i(y)=\Phi(y)\cap X_i,
\]
and $\Cr^i_\Phi=L_{\rho_i}=L_{\lambda_i}$.  The space   $X$ is obtained via the attachments
\[
\Cyl_{\lambda_0}\cup_{Y_1} ( \stackrel{\rho_1}{\Lra} X_1\stackrel{\lambda_1}{\Lla}) \cup_{Y_2}\cdots \cup_{Y_{n-2}} (\stackrel{\rho_{n-1}}{\Lra} X_{n-1}\stackrel{\lambda_{n-1}}{\Lla})\cup_{Y_n} \Cyl_{\rho_n},
\]
where $\Cyl_g$ denotes the mapping cylinder of  a tame continuous map $g$.

The tame blowdowns $\lambda_i$ and $\rho_i$ carry natural weights. To define them we first need to define the tame maps
\[
T_i^+: X_i\ra (0,\infty],\;\; T_i^-: X_i\ra (0,\infty]
\]
where for every $x\in X_i$, we denote by $T_i^+(x)$ the moment of time when the   flow line  through $x$ intersects $Y_i$,  and by $T_i^-(x)$ the  moment of time when the  backwards  flow line trough $x$ intersects $Y_{i+1}$. Equivalently,
\[
T_i^+(x)=\sup\bigl\{\, t>0;\;\;  f(\Phi_t(x))\geq d_i\,\bigr\},\;\; T_i^-(x)=\sup\bigl\{\, t>0;\;\;f(\Phi_{-t}(x)\,)\leq d_{i+1}\,\bigr\}.
\]
Observe that
\[
T_i^\pm(x)=\infty \Longleftrightarrow x\in \Cr^i_\Phi.
\]
In Section \ref{s: 9} we have proved that the tame functions $T_i^\pm$ are continuous. Now define
\[
w_i^\pm:=\frac{1}{T_i^\pm}.
\]
The  functions $w_i^\pm$ are continuous, nonnegative  and
\[
w_i^\pm(x)=0 \Longleftrightarrow x\in \Cr^i_\Phi.
\]
In other words, $w_i^+$ is a weight for $\rho_i$, and $w_i^-$ is a weight for $\lambda_i$.

\begin{definition} A \emph{weighted chain} of tame flips is a sequences  $\BXi=\BXi(\lambda_i,\rho_i, w^\pm_i)$ of flips
\[
Y_{-1} \stackrel{\rho_0}{\Lra} X_0\stackrel{\lambda_0}{\Lla} Y_1\stackrel{\rho_1}{\Lra} X_1\stackrel{\lambda_1}{\Lla}\, \cdots\,   \stackrel{\rho_{n-1}}{\Lra} X_{n-1}\stackrel{\lambda_{n-1}}{\Lla}Y_n\stackrel{\rho_n}{\Lra} X_n\stackrel{\lambda_n}{\Lla} Y_{n+1},
\]
and weights $w_i^+$ for $\rho_i$, $w_i^-$ for $\lambda_i$ such that $X_0$ and $X_n$ are finite sets, $Y_{-1}=Y_{n+1}=\emptyset$, and $L_{\rho_i}=L_{\lambda_i}$, $\forall i$. The tame space associated to a weighted chain is defined as
\[
\Cyl_{\lambda_0}\cup_{Y_1}\Cyl_{\rho_1}\cup_{X_1}\cup\Cyl_{\lambda_1}\cup_{Y_2}\, \cdots\,\cup_{Y_n} \Cyl_{\rho_n} X_n.
\]
We denote by $\eC_w$ the set of weighted tame chains and, for every compact tame set $X$, we denote  by $\eC_w(X)$ the set of weighted chains  whose associated space is $X$. \qed
\end{definition}

The  discussion  preceding the above definition shows that we have a natural map
\[
\BXi:\Flg(X)\ra \eC_w(X),\;\; \Phi\longmapsto \BXi_\Phi.
\]
Under this map, the stationary points of $\Phi$  correspond bijectively to the  points in the blowup loci $L_{\rho_i}=L{\lambda_i}$.   The exceptional loci of $\rho_i$ correspond to unstable links of stationary points, while  the exceptional loci of $\lambda_i$ correspond to the stable links.
\begin{theorem} The map $\BXi:\Flg(X) \ra \eC_w(X)$ is surjective.
\end{theorem}

\proof  The strategy is simple.  Given a weighted chains of  flips $\Xi(\lambda_i, \rho_i, w_i^{\pm}, 0\leq i\leq n)\in \eC_w(X)$  we will construct local flows and Lyapunov functions  on the various mapping cylinders  associated to this chain, and then we concatenate them.  It suffices to do this for a single blowdown map $\beta: Y\ra X$, with weight $w$.

For us, a local  tame flow on a tame set $S$  will be a tame continuous map
\[
\Psi: R_S\ra S,
\]
where $R_S\subset \bR\times S$ is a tame subset  such that

\begin{figure}[h]
\centerline{\epsfig{figure=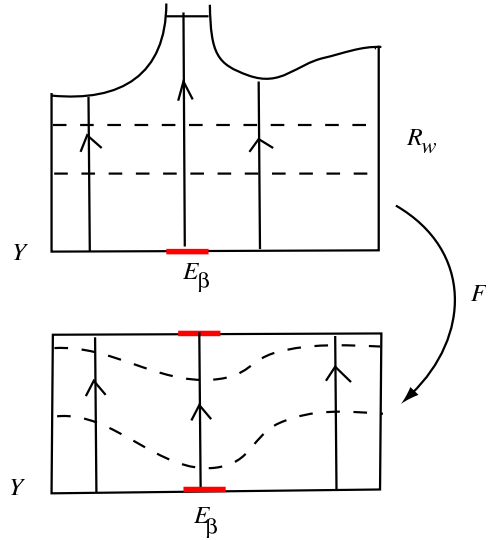,height=2.8in,width=2.5in}}
\caption{\sl Constructing  a gradient flow on the mapping cylinder of a tame map.}
\label{fig: tame24}
\end{figure}

\begin{itemize}
\item $\{0\}\times S\subset R_S$,

\item  for every $s\in S$, the set $I_s:=\{t\in \bR;\;\;(t,s)\in R_S\}$ is  an interval of positive length, and

\item  for every $s\in S$ and $t_0,t_1\in I_s$ such that $t_0+t_1\in I_s$ we have
\[
\Phi_{t_0+t_1}(s)=\Phi_{t_0}\bigl(\,\Phi_{t_1}(s)\,\bigr).
\]
\end{itemize}
Define
\[
T: Y\ra (0,\infty],\;\;T(y)=\begin{cases}
\frac{1}{w(y)} & w(y)\neq 0\\
\infty & w(y)=0
\end{cases},
\]
and set
\[
R_{w}:= \bigl\{ (y,t)\in Y\times [0,\infty];\;\; t\leq T(y)\,\bigr\}.
\]
Fix a tame homeomorphism $F:R_w\ra Y\times [0,1]$ such that, $F(y,0)=(y,0)$, and  the diagram below is commutative
\[
\begin{diagram}
\node{R_w}\arrow[2]{e,t}{F}\arrow{se,b}{}\node[2]{Y\times [0,1]}\arrow{sw,b}{}\\
\node[2]{Y}
\end{diagram}
\]
where the maps $R_w,Y\times [0,1]\ra Y$ are the natural projections (see Figure \ref{fig: tame24}).

Consider the translation flow on $Y\times [-\infty, \infty]$, whose stationary points are $(y,\pm\infty)$, $y\in Y$.

It restricts to a (local) flow on $R_w$ whose trajectories are the vertical lines $[0,T(y)]\ni t\mapsto(y,t)$.   Via $F$  we obtain  a tame  local flow on $Y\times [0,1]$, whose orbits are the vertical segments $\{0\}\times [0,1]$.  The bottom point $(y,0)$ will reach the top point $(y,1)$ in $T(y)$ units of time.

The natural map $Y\times [0,1]\ra [0,1]$ decreases strictly along the trajectories  of this local flow, and thus defines a Lyapunov function. The points in $E_\beta\times\{1\}$ are stationary points.

This  local flow descends to a local flow on
\[
\Cyl_\beta=Y\times [0,1]\cup_\beta X
\]
where the points in $L_\beta\subset X$ are stationary points.

\qed

To transform the above theorem into an useful  technique for producing   gradient like flows, we need to  explain how to construct  weighted blow-down/up maps.

Note that given a compact tame space $Y$ and $E\subset Y$ a compact tame subset, then $X/E$ is a compact tame space, and the natural projection $Y\ra Y/E$ is a blowdown map. We would like to investigate the opposite process.

Suppose   we are given a  compact tame space $X$, a point  $p_0\in X$, and a continuous tame function $w: X\ra [0,\infty)$ such that $w^{-1}(0)=\{x_0\}$.

We can then find $r_0>0$ such that the induced  map
\[
w:\{\, 0<w<r_0\}\ra (0,r_0]
\]
is a  (tamely) locally trivial fibration.   The level sets $\{w=\ve\}$, $\ve\in (0,r_0)$ are all (tamely) homeomorphic. We will refer to any one of them as the $w$-link\footnote{We do not know if the homeomorphism type of the $w$-link  depends on the weight $w$, or   that it is homeomorphic to the link of the point $x_0$ in $X$  as defined in Appendix \ref{s: b}.} of $p_0$, and we will denote it by $\eL_w(p_0)$.

Observe that the closed neighborhood $\{w\leq r_0\}$ of $p_0$ is  tamely homeomorphic to the cone on $\eL_w$, or equivalently the mapping cylinder of the constant map $\eL_w\ra \{p_0\}$.

Consider now an  arbitrary, tame continuous map $\mu: \eL_w\ra E$, where $E$ is a tame compact  set. Observe that   the canonical map from the mapping cylinder of $\mu$ to the mapping cylinder of the constant map  $\eL_w\ra \{p_0\}$  is a blowdown map $\Cyl_\mu\ra \{w\leq r_0\}$ with blowup locus $\{x_0\}$, and exceptional locus $E$.  We can now define the blowup space $\widehat{X}_{w,\mu}$  to be
\[
\widehat{X}_{w,\mu}= \{w\geq r_0\}\cup_{\eL_w} \Cyl_\mu.
\]
\begin{figure}[h]
\centerline{\epsfig{figure=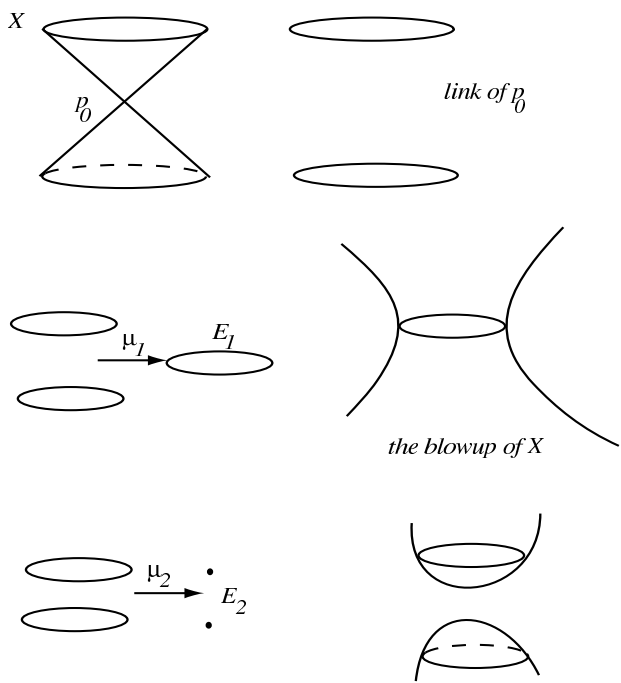,height=3.5in,width=3.2in}}
\caption{\sl Blowing up the vertex of a cone in two different ways.}
\label{fig: tame25}
\end{figure}
\begin{ex} (a) Suppose $X$ is the Euclidean space $\bR^n$,  $p_0$ is the origin, and $w$ denotes the Euclidean norm. The $w$-link of $p_0$ is the round sphere $S^{n-1}$. If we denote by $\mu:S^{n-1}\ra \bRP^{n-1}$ the   canonical double covering, then the blowup $\widehat{X}_{p_0,w}$ is the usual blowup in algebraic geometry.

(b) Suppose $X$ is the semialgebraic cone (see Figure \ref{fig: tame25})
\[
X=\bigl\{\, (x,y,z)\in \bR^3;\;\; z^2=x^2+y^2,\;\;|z|\leq 1\,\bigr\}
\]
and $p_0$  is the origin. Assume $w(x,y,z)=|z|$. Then the link of $p_0$ consists of two circles.

We can choose  $\mu$ in many different ways. For example, we can choose $\mu=\mu_1: S^1\sqcup S^1 \ra S^1$ to be the natural identification map, or we can choose $\mu=\mu_2: S^1\sqcup S^1\ra \{0,1\}$ to be the map which collapses  each of the two circles to a different point. The resulting blowup spaces $\widehat{X}_{w,\mu}$ are depicted in Figure \ref{fig: tame25}.

These types  of blowups appear in Morse theory, when we cross a level set of a $3$-dimensional  Morse function   containing   saddle point. \qed

\end{ex}

%% file: tameflow11.tex
In the sequel,  a  combinatorial simplicial complex (or CSC  for
brevity) is a \emph{finite} collection $\eK$ of nonempty  finite
sets  with the property  that
 \[
 A\in K,\;\;B\subset A \Longrightarrow  B\in \eK.
 \]
 The sets   in $K$ are called the \emph{open faces} of  $K$.   The union of all the sets in $K$  is called the vertex set of $K$ and will be denoted by  $\eV(\eK)$.  The dimension of an  open face $A\in \eK$ is the nonnegative integer
 \[
\dim A:= \#A -1.
\]
We set
\[
\dim K:=\max\bigl\{ \,\dim A;\;\;A\in \eK\,\bigr\}.
\]
A vertex  is a $0$-dimensional face.

For every  subset $ \eA\subset \eK$  we define its
\emph{combinatorial closure}  to be
\[
\cl_c(\eA)=\bigl\{\, B\in\eK;\;\;\exists A\in \eA:\;\;B\subset
A,\bigr\}.
\]
A \emph{subcomplex} of $\eK$ is a subset  $\eA\subset \eK$ such
that $A=\cl_c(\eA)$. The \emph{$\ell$-th skeleton} of $\eK$ is the
subcomplex
\[
\eK_\ell=\bigl\{\, A\in K;\;\;\dim A\leq \ell\,\bigr\}.
\]
For any  subset $S\subset \eV(\eK)$ we denote by   $\eF(S)$ the
subcomplex of $\eK$ spanned by the  faces with vertices in $S$,
\[
\eF(S):= \bigl\{ \,A\in\eK;\;\;A\subset S\,\bigr\}.
\]
For any vertex  $v$  of $\eK$,   we denote by $L(v)=L(v,\eK)$ the
set of vertices adjacent to $c$ in $\eK$, and we set
 \[
 S(v)=S(v,\eK):=\{p\}\cup L_p.
 \]
 The \emph{combinatorial star} of $v$ in $\eK$ is then the subcomplex
 \[
 \eS(v)=\eS(v,\eK):=\eF(S(v)),
 \]
  while the \emph{combinatorial link} of $v$ in $\eK$ is the subcomplex
 \[
 \eL(v)=\eL(v,\eK):=\eF(L(v)).
 \]
For every finite set $S$ we denote by $\bR^S$ the vector space of
functions $S\ra \bR$.   $\bR^S$ has a canonical basis  consisting
of the Dirac functions $(\delta_s)_{s\in S}$, where
\[
\delta_s(s')=\begin{cases}
1 & s'=s\\
0 &s'\neq s.
\end{cases}
\]
For any  subset $A\subset S$ we  denote by $\conv (A)$ the convex
hull of the set
\[
\bigl\{ \delta_a;\;\;a\in A\,\bigr\}\subset \bR^S.
\]
If $\eK$ is a CSC, then  the \emph{geometric realization} of $\eK$
is the closed subset
\[
|\eK|=\bigcup_{A\in \eK} \conv(A) \subset \bR^{\eV(\eK)}.
\]
We denote by $\St(v)$ the geometric realization of $\eS(v)$ and by
$\Lk(v)$ the geometric realization of $\eL(v)$.

\begin{ex} (a) Suppose $(P,\leq)$ is a finite poset (partially ordered set).   Then the  \emph{nerve} of $(P,\leq)$ is the CSC $\eN(P,\leq)$, with vertex set $P$, and open  faces given by the chains of  $P$, i.e., the linearly ordered subsets of $P$.  For any poset $P$ we will denote by $|P|$ the geometric realization of its nerve
\[
|P|:=|\eN(P)|.
\]
We say that two posets are homeomorphic or homotopic if the
geometric realizations of their nerves are such.

\noindent (b)  Suppose $\eK$ is a CSC. Then $\eK$ is a finite
poset,  where the order relation is given by inclusion.  The nerve
of $(\eK,\subset)$ is called the  \emph{first barycentric
subdivision of $\eK$}, and it is denoted by $D\eK$.   We define
inductively
\[
D^{n+1}\eK:=D(D^n\eK)
\]
We say that $D^n\eK$ is the $n$-th barycentric subdivision of
$\eK$.

\noindent (c) Suppose $\eK_1$ and $\eK_2$ are two  CSC's with
disjoint vertex sets $\eV_1,\eV_2$. We define the \emph{join} of
$\eK_1$ and $\eK_2$ to be the   CSC $\eK_1\ast\eK_2$ with vertex
set $\eV_1\cup \eV_2$, and  faces $F_1\cup F_2$, $F_i\in \eK_i$.
The  join of a CSC and a point  which is not a vertex of $\eK$ is
called the cone on $\eK$ and it is denoted be ${\rm Cone}\,(\eK)$.

\noindent (d) If $\eK$ is a $CSC$ then the suspension of $\eK$ is
the CSC  $\Sigma\eK$ defined as the join of $\eK$ with  the CSC
$S^0=\bigl\{ \{N\}, \{S\}\,\}$, where $N,S\not\in \eV(\eK)$.  The
$n$-th iterated suspension of $\eK$ is defined inductively as
\[
\Sigma^n\eK:=\Sigma(\Sigma^{n-1}\eK).\proofend
\]

\end{ex}

If $\eK_0$ and $\eK_1$ are two CSC-s, then a  CSC morphism from
$\eK_0$ to $\eK_1$ is a map
\[
f:\eV(\eK_0)\ra \eV(\eK_1),
\]
 such that
\[
A\in\eK_0\Longrightarrow f(A)\in\eK_1.
\]
A morphism $f:\eK_0\ra \eK_1$ induces  a morphism  $Df: S\eK_1\ra
S\eK_1$ between the  first barycentric subdivisions, and a
continuous, piecewise linear map $f_\sharp:|\eK_0|\ra |\eK_1|$.

\begin{definition} (a)  A  \emph{dynamical orientation} on the CSC $\eK$ is a  binary relation $\rightsquigarrow$ on $\eV(\eK)$  with the following properties.

 \begin{itemize}

\item  If $u\rightsquigarrow v$  then $\{u,v\}$ is a one
dimensional face of $\eK$.

 \item For any open face $A\in \eK$, the restriction of $\rightsquigarrow$ to $A$  is a linear order.

 \end{itemize}

\noindent (b)  A \emph{combinatorial flow} is  a pair $(\eK,
\rightsquigarrow)$, where $\eK$ is a CSC and $\rightsquigarrow$ is
a dynamical orientation  on $\eK$.\qed
\end{definition}

 If $(\eK,\rightsquigarrow)$ is a  combinatorial flow, and $p\in \eV(K)$ then we set
 \[
 L(p\rightsquigarrow):=\bigl\{ u\in \eV(\eK);\;\;p\rightsquigarrow u\,\bigr\},\;\; \eL(p,\rightsquigarrow):=\eF(\,L(p\rightsquigarrow)\,),
 \]
 \[
 W(p\rightsquigarrow):=L(p\rightsquigarrow)\cup \{p\},\;\;\eW(p\rightsquigarrow):=\eF(W(p\rightsquigarrow)\,).
 \]
The sets $W(\rightsquigarrow p)$, $L(\rightsquigarrow p)$ etc.,
are defined in a similar fashion. We will say that $\eL(p\rsq)$ is
the  \emph{ unstable  combinatorial link}.

Using the construction  in Example \ref{ex: simplicial} we can
associate to any combinatorial flow $(\eK, \rightsquigarrow)$ a
tame flow $\Phi=\Phi^\rightsquigarrow$ on the geometric
realization $|\eK|$. We will say that $\Phi^\rightsquigarrow$ is
the \emph{simplicial flow} determined by the dynamical orientation
$\rightsquigarrow$.

 \begin{theorem} Suppose $(\eK,\rightsquigarrow)$ is a combinatorial flow,   and $\Phi$ is the  simplicial  flow on $|\eK|$ associated to $\rightsquigarrow$.  Then the following hold.

\smallskip

\noindent  (a) The map
\[
 \eV(\eK)\ni v\mapsto \delta_v\in|\eK|
 \]
 is a bijection from the vertex set of $\eK$ to the set of stationary points of $\Phi$.

 \noindent (b) For every vertex $v$ of $\eK$, the Conley index of   $\delta_v\in |\eK|$ is homotopy equivalent to  the pointed space
 \[
 |{\rm Cone}\,(\eL(v\rsq)|/|\eL(v\rightsquigarrow)|.
 \]
 \label{th: simp-flow}
 \end{theorem}

 \proof  Part (a) is obvious.  To prove  (b) observe  that the  star $\St(v)$ is a compact, flow invariant neighborhood of  $\delta_v$. Thus, the Conley index of $\delta_v$ in $\Aff(\eK)$  is homotopy equivalent to the Conley index of $\delta_v$ in $\St(v)$.

 Observe that we have a partition
 \[
 S(v)=\{p\} \sqcup \eL(\rightsquigarrow v)\sqcup \eL(v \rightsquigarrow).
 \]
 Now define
 \[
 f: S(v)\ra \{-1,0,1\}
 \]
 by setting
 \[
 f(u) :=\begin{cases}
 0  & u=v\\
 1 &u\in \eL(\rightsquigarrow v)\\
 -1 & u\in\eL(v\rightsquigarrow).
 \end{cases}
 \]
 The function $f$ induces a piecewise linear function $\St(v)\ra [-1,1]$ which, for simplicity, we continue to denote by $f$.

 From the explicit description in Example \ref{ex: simplicial}  of the canonical tame flow on an affine simplex we deduce that $\delta_v$ is a Morse like stationary point of the flow $\Phi$ on $\St(v)$, and $f$ is a tame local Lyapunov function adapted to $\delta_v$. The result now follows from Theorem \ref{th: conley}. \qed

\begin{ex} The cheapest way of producing a  dynamical orientation on  a CSC  $\eK$ is to  choose an injection
\[
f:\eV(\eK)\ra \bR.
\]
Then we define
\[
x\rsqf y\Longleftrightarrow f(x)>f(y),\;\;\{x,y\}\in \eK.
\]
Then  $f$ defines a piecewise  linear function $f:|\eK|\ra \bR$
which is a tame Lyapunov function  for the  simplicial flow
determined by  $\rsqf$.

Alternatively, the restriction  of a generic linear map
$f: \bR^{\eV(K)}\ra \bR$ to the affine realization $|\eK|$ is injective on the
vertex  set.    This function is a stratified Morse function in the sense of Goresky-Macpherson,   and in this case, the Conley index computations   also follow from the  computaions in \cite{GM} of the local Morse data  of a stratified Morse function. \qed
\end{ex}

Let us present a few applications of this result to the homotopy
theory of  posets.  We need to introduce some terminology

Suppose $(P,\leq )$ is a finite poset. Recall that for any $x,y\in
P$  we define the order intervals
\[
[x,y]:=\bigl\{z \in P;\;\;x\leq x \leq y\,\bigr\},
\;\;(x,y)=\bigl\{ z\in P;\;\;x<z<y\,\bigr\}
\]
and  we say  that $y$ \emph{covers} $x$ if $[x,y]=\{x,y\}$. We
write this $y\gtrdot x$. We define
\[
P_{<x }:=\bigl\{ x\in P;\;\;x<y\,\bigr\}.
\]
An order ideal of $P$ is a subset $I\subset P$ such that
\[
x\in I \Longrightarrow P_{\leq x}\subset I.
\]
For every chain $x_0<x_1<\cdots  <x_k$ in $P$, we   will refer to
the  integer $k$ as the length of the    chain.    Given $x\leq
y$,  we define $\ell(x,y)$ the be the  maximal length of a chain
originating at $x$ and ending at $y$. Observe that
\[
x\lessdot y \Longleftrightarrow \ell(x,y)=1.
\]
 Finally, we will say that a poset is contractible, if it is homotopic to the poset consisting of a single point.

A function  $f: P\ra \bR$ on a poset $P$ is called
\emph{admissible} if
\[
f(x)=f(y) \Longrightarrow \;\;\mbox{$x$ and $y$ are not
comparable}.
\]
Suppose $f: P\ra \bR$ is admissible.  For every $x\in P$ we set
\[
V^+(x)=V^+(x,f):=\bigl\{\,
y>x;\;\;f(x)>f(y)\,\bigr\},\;\;\eS^+(x)=\eS^+(x,f):=\{x\}\cup
V^+(x),
\]
\[
 V^-(x)=V^-(x,f)=\bigl\{ z<x;\;\;f(z)>f(x)\,\bigr\},\;\;\eS^-(x)=\eS^-(x,f):=\{x\}\cup V^-(x).
\]
\begin{remark} Here is the intuition behind the  sets $V^\pm(x,f)$.       Note that these sets are empty for every $x\in P$ if and only if $f$ is a strictly decreasing function. In other words, the sets $V^\pm(x)$ collect the ``\emph{violations}''     at $x$ of the strictly decreasing condition.\qed
\end{remark}
The admissible function $f$ defines  a partial order
\[
x\precf y \Longleftrightarrow f(x)>f(y) \;\;{\rm and}\;\;  x<y,
\]
so that
\[
V^-(x)=\bigl\{ y\in P;\;\;y\prec_f x\,\bigr\},\;\; V^+(x)=\bigl\{
z\in P;\;\; x\prec_f z\,\bigr\}.
\]
If $f: P\ra \bR$ is a  we have a   simplicial flow $\Phi^f$ on the
nerve of $P$  given by the dynamical orientation
\[
x\rsqf y \Longleftrightarrow f(x)>f(y)   \mbox{ and $x$ and $y$
are comparable elements of $P$}.
\]
Every  point $x\in P$ is a stationary  point of this flow.  We
denote by $h_f(x)$ its Conley index.

The  unstable combinatorial  link $\eL(x\rsqf)$ of $x$ is the
nerve of the poset
\[
V^+(x) \cup (P_{<x}\setminus V^-(x)),
\]
which  is the  join
\[
\eN\bigl( V^+(x)\,\bigr)\ast \eN\bigl(\,P_{<x}\setminus
V^-(x)\,\bigr).
\]
\begin{definition} Suppose $f$ is a real valued  admissible function on the poset $P$. A point  $x\in P$ is called a regular point of $f$ if one of the posets  $V^+(x)$ or $P_{<x}\setminus V^-(x)$ is contractible.\qed
\end{definition}

\begin{corollary}  If $x$ is a regular point of  $f$ then its  Conley index is trivial.\qed
\end{corollary}

\begin{definition} Suppose $f:P\ra \bR$ is a real valued  admissible function on a poset $P$.

\noindent (a) We say that $f$ is \emph{coherent} if
\[
x\precf y \Longrightarrow \mbox{ $f$ is strictly increasing on the
interval $[x,y]$.}
\]
 The \emph{order} of a coherent function is the nonnegative integer
\[
\omega(f):=\max\bigl\{\ell(x,y);\;\; \mbox{$f$ is increasing on
$[x,y]$}\,\bigr\}.
\]
A \emph{Morse-Forman} function is a coherent function of order
$\leq 1$.

\noindent (b) We say that $f$ satisfies the condition $C_+$ if it
is coherent and there exists  a map
\[
C_+: P\ra P
\]
 such that
\[
x\leq C_+(x)\;\;{\rm and}\;\;\eS^+(x,f)= [x, C_+(x)],\;\;\forall
x\in P.
\]
We  say that $f$ satisfies the condition $C_-$ if it is coherent
and there exists a map $C_-:P\ra P$ such that
\[
x\geq C_-(x)\;\;{\rm and}\;\;\eS^-(x,f)= [C_-(x), x],\;\;\forall
x\in P.
\]
We say  that $f$ satisfies the condition $C$ if it satisfies  both
$C_+$ and $C_-$. \qed
\end{definition}

\begin{ex} (a) Any strictly decreasing function  on a finite poset $P$ is a coherent  function of order zero.

\noindent (b) Suppose $\eK$ is a  CSC with vertex set $V$, i.e. an
ideal of $2_*^V$. Then a discrete Morse function $f:\eK\ra \bR$ of
the type introduced  by R. Forman in \cite{For} is a Morse-Forman
function on the poset of faces  of a CSC.

\noindent (c) If $f:P\ra \bR$ satisfies $C_-$, and $I\subset P$ is
an ideal, then $f|_I$  satisfies $C_-$.

\noindent (d) In  Figure \ref{fig: tame12} we have depicted a
coherent  function of order two on the poset of faces of the  two
dimensional  simplex. The arrows   indicate the dynamical
orientation  determined by this function. This function also
satisfies condition $C$.

\begin{figure}[h]
\centerline{\epsfig{figure=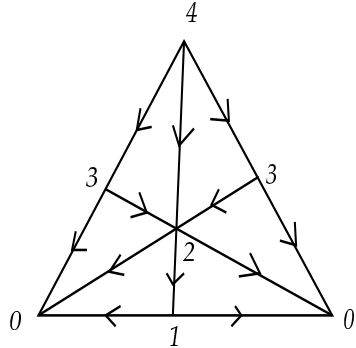,height=1.8in,width=1.8in}}
\caption{\sl A coherent function of order $2$.} \label{fig: tame12}
\end{figure}

\qed
\end{ex}

\begin{corollary}    If $f: P\ra \bR$ is a function satisfying condition $C_+$, then any point  $x\in P$ such that $x\neq C_+(x)$ is a regular point.
\end{corollary}

\proof  If $x\neq C_+(x)$, then the nerve of  $V_+(x)$ is a cone
with vertex $C_+(x)$, hence contractible. \qed

\smallskip

\emph{In the remainder of this section, we will assume that $P$ is
the poset  $\eF(X)$ of faces of a regular  $CW$-decomposition  of
a compact space $X$.}

\smallskip

Observe that  the intersection of two faces is either empty, or a
face of $X$, i.e., $\eF(X)$ is a meet semilattice.    By
\cite[Thm. III.1.7]{LW}, geometric realization of the nerve of
the poset $\eF(X)$ is PL homeomorphic to $X$. In particular, if
$F\in \eF(X)$   is a closed face, then $\eF(X)_{<F}$ is the union
of all the proper faces  of $F$ so that    the geometric
realization of  the nerve of $\eF(X)_{<F}$ is   PL homeomorphic to
the $PL$ space $\pa F\cong S^{\dim F-1}$. Similarly, the geometric
realization of $\eF(X)_{\leq F}$ is     PL homeomorphic to the
closed  ball $F\cong \bD^{\dim F}$ equipped with its the natural
PL structure.

Suppose $f: \eF(X)\ra \bR$ is a coherent function. For any face
$F$ we denote by $V^+_{max}(F)$ the maximal elements in $V^+(F)$.
Since $f$ is coherent, we deduce
\[
V^+(F)=\bigcup_{T\in V_{max}^+(F)}  (F, T].
\]
Observe that
\[
(F,T_1]\cap (F, T_2]=\begin{cases}
\emptyset & T_1\cap T_2=\emptyset\\
(F, T_1\cap T_2] &  T_1\cap T_2\neq \emptyset.
\end{cases}
\]
Denote by $M^+(F)$ the  CSC with vertex set $V^+_{max}(F)$ such
that $\{T_1,\dotsc, T_k\}$ is a face if and only if $T_1\cap\cdots
\cap T_k\neq \emptyset$.   In other words, $M^+(F)$ is the nerve
of the cover  $\bigcup_{T\in V_{max}^+(F)}  (F, T]$. The   order
intervals $(F,G]$ are contractible and we deduce from the Nerve
Theorem  \cite[Thm. 10.6]{Bjorn} that  the nerve of $V^+(F)$ and
$M^+(F)$ have the same homotopy type.  We obtain the following
consequence.

\begin{corollary}  Suppose $f:\eF(X)\ra \bR$ is a coherent function.  If $M^+(F)$ is a (non-empty) contractible CSC, then $F$ is a regular point of $f$.\qed
\end{corollary}

\begin{remark}  Observe that the  coherent function $f:\eF(X)\ra \bR$ satisfies condition $C_+$ if and only if, for every face  $F$, the simplicial complex $M^+(F)$ is either empty, or consists of a single point.\qed
\end{remark}

Suppose now that $f$ satisfies $C_-$. Fix a closed face $F$, and
set $F_-=C_-(F)$. In other words
\[
B\leq F,\;\; f(B)\geq f(F) \Longleftrightarrow B\in [F_-,
F]\Longleftrightarrow V^-(F)=[F^-, F).
\]
For any $B<F$ set
\[
\eC(B,F):=\eF_{<F}(X)\setminus [B,F).
\]

\begin{lemma} For any $B <F$, the geometric realization of the nerve of the poset  $\eC(B,F)$ is homeomorphic to  the ball $\bD^{\dim F}$.
\label{lemma: low-con}
\end{lemma}

\proof   Denote by $Y$ the  union of proper faces of $F$ which do
not contain $B$, i.e.
\[
Y=_{G\in \eC(B,F)} G.
\]
$Y$ is a PL space and the geometric realization of the  nerve of
$\eC(B,F)$ is PL homeomorphic to  $Y$.

 We set $n:=\dim F$, and we assume $F\subset \bR^n$. Choose a point $b_0$ a point  in the relative interior  of $B$,    and  for every $r>0$ denote by $L_r$ the intersection of $F$ with the sphere of radius $r$ in $\bR^n$ centered at $b_0$.  For  $r$ sufficiently small, $L_r$ is homeomorphic to  a closed ball of dimension $n-1$.    For every $x\in F\setminus \{b_0\}$ we denote by $\si_r(x)$ the intersection of  the line $[b_0,x]$ with the link $L_r$.      For $r>0$ sufficiently small, the map $\si_r$ defines a homeomorphism $Y\ra L_r$. \qed

Note that $\eC(B,F)$ consists of all closed faces of $\pa F$ which
do not contain  $B$.  Denote by $\pa F\setminus \St_B$ the union
of all the closed faces  $F'\in \eC(B,F)$.

\begin{theorem} Suppose $f:\eF(X)\ra \bR$ satisfies $C_-$.    If $F\neq  C^-(F)$ then $F$ is a regular point of $f$,  while if $F=C_-(F)$ then the Conley index of $F$ with respect to the simplicial flow defined by $f$ is
\[
h_f(F)= \bigl[\, |{\rm Cone}( S^{\dim F-1}\ast  M^+(F)\,)|, \,
|S^{\dim F-1}\ast  M^+(F)|\,\bigr]
\]
\[
\simeq \bigl[\, {\rm Cone}\,\Sigma^{\dim F}|M^+(F))|,\,
\Sigma^{\dim F}|M^+(F))|\,\bigr].
\]
\label{th: cminus}
\end{theorem}

\proof If $F\neq C^-(F)$ then Lemma \ref{lemma: low-con} shows
that  the poset  $\eF(X)_{<F}\setminus V^-(F)$ is contractible,
and  thus $F$ is a regular point. If $F=C^-(F)$,  then
$\eF(X){<F}\setminus V^-(F)=\eF(X)_{<F}$, and the poset
$\eF_{<X}$ is  PL homeomorphic to the sphere $S^{\dim F-1}$. The
poset $V^+(F)$ is  homotopic to  $|M^+(F)|$ and thus
\[
|\eL(F\rsqf)|\simeq  S^{\dim F -1}\ast |M^+(F)|\simeq
\Sigma^{\dim F}|M^+(F))|.
\]
The result now follows from  Theorem \ref{th: simp-flow}. \qed

For any face  $F\in \eF(X)$  we denote by
$\tilde{\eP}_{M^+(F)}(t)$ the  Poincar\'e  polynomial of the
reduced homology of $|M^+(F)|$, with the convention that
\[
\tilde{\eP}_{M^+(F)}(t)=t^{-1}\;\;{\rm if}\;\; M^+(F)=\emptyset.
\]
Denote by $\tilde{M}_{F, t}(t)$  the Poincar\'e  polynomial of the
reduced homology of $|\eL(F\rsqf)|$.  Then
\[
\tilde{M}_{F, t}(t)= \eP_{\Sigma^{\dim F}|M^+(F))|}(t)= t^{\dim
F}\tilde{\eP}_{M^+(F)}(t),
\]
and from (\ref{eq: morse-ineq}) we  deduce the Morse inequalities
\begin{equation}
\sum_{F=C_-(F)} t^{\dim F+1}\tilde{\eP}_{M^+(F)}(t)\succeq
\eP_X(t). \label{eq: morse-ineq1}
\end{equation}
Observe that if $f$ satisfies $C$, then for any face $F$ the
simplicial complex $M^+(F)$  is either empty, i.e. $C_+(F)=F$, or
consists of a single point, and $F\neq C_+(F)$. In this case, the
Morse inequalities       are very similar to the classical ones
\begin{equation}
\sum_{F=C_-(F)=C_+(F)} t^{\dim F}\succeq \eP_X(t). \label{eq:
morse-ineq2}
\end{equation}

\begin{ex} In the left-hand side of  Figure \ref{fig: tame14} we have depicted a  coherent function $f$ of order two on the  poset of faces of a $2$-dimensional    (affine) simplicial complex $X$. It satisfies the condition $C_-$, but it does not satisfy the condition $C_+$.  The simplicial flow determined by this function is depicted the right-hand side of the figure.

The vertices labelled by $4$, $1$ and $-1$  correspond to the
faces $F$ satisfying the condition
\[
F=C_-(F),
\]
so these are the only stationary points of the flow  which could
have  nontrivial Conley index, and thus   could potentially affect
the topology of  $X$.   For a vertex $v$ such that $f(v)=1$, the
simplicial complex $M^+(v)$ is contractible (it corresponds to the
barycenter of an edge labelled $0$) and thus  the Conley index is
trivial.

If  $v$ is a  vertex such that $f(v)=4$, then the simplicial
complex $M^+(v)$ consists of two points (labelled $A,B$ in the
figure) and we deduce that the Conley index of such a point is
$[S^1,\ast]$.  In this case we observe that the Morse inequalities
become equalities,  and  we see  that  we can use the flow to
collapse $X$ to a wedge of $3$ circles. \qed

\begin{figure}[h]
\centerline{\epsfig{figure=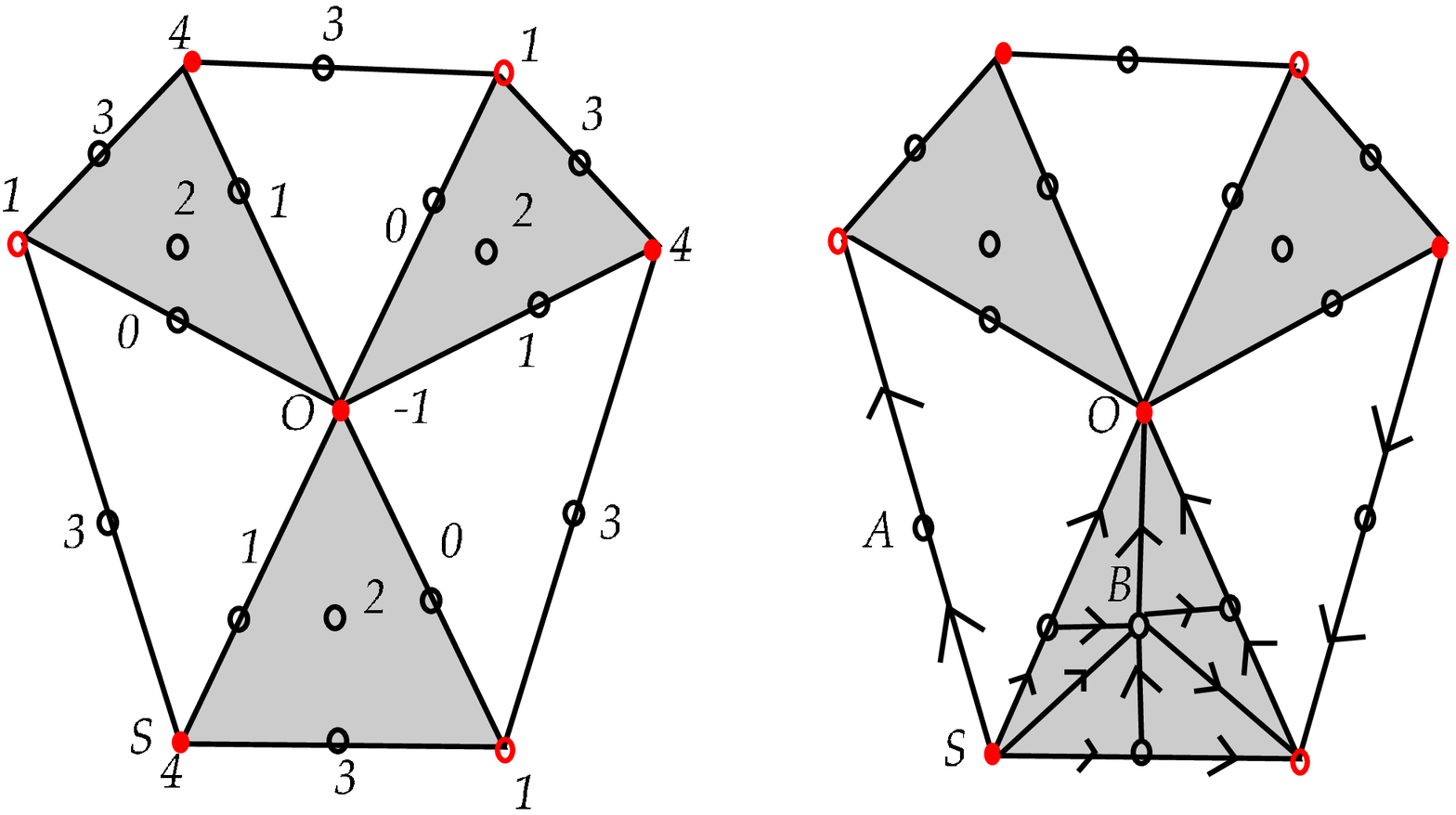,height=2.1in,width=3.8in}}
\caption{\sl A  $C_-$ function of order $2$ and its associated
``gradient'' flow.} \label{fig: tame14}
\end{figure}

\end{ex}

\begin{remark} (a) It is not hard to prove that if $f$ is a coherent function  on the poset $\eF(X)$ of faces of a polytopal decomposition of PL space $X$, then we can modify $f$ to a coherent \emph{injective} function $g:\eF(X)\ra \bR$  such that $\rsqf=\stackrel{g}{\rsq}$.

\noindent (b) If $f:\eF(X)\ra \bR$ satisfies the condition $C$,
then  we can use  the flow  determined by $f$ to    extract
information about the \emph{simple} homotopy type   of $X$.    If
the order of $f$ is $\leq 1$, then  $f$  is a discrete Morse
function of the type  introduced by R. Forman, and  many of the
results \cite{For} follow  from general properties  of  the Conley
index and tame flows.  We will not pursue   this point of view.

\noindent (c) One  of the drawbacks of this theory  is the lack
of a general technique for constructing $C_-$ functions.      Here
is one possible method of  addressing this issue.

 We observe that any  CSC  with vertex set $\eV$ is an ideal of the poset  $2_*^\eV$ of nonempty  subsets of $\eV$. Any function $f:2_*^\eV\ra \bR$   satisfying the condition $C_-$ restricts to a function with the same property on any ideal.  Thus  a good understanding of the $C_-$ functions  on the poset of faces of a simplex  could potentially lead to a large supply  of $C_-$ functions on the poset of faces of a CSC.

 A $C_-$ function on $2^\eV_*$ has a very combinatorial description. It is a function $f:2_*^\eV\ra \bR$ satisfying the conditions

\begin{itemize}
\item  If $A\varsubsetneq B$ then  $f(A)\neq f(B)$.

\item  If $A, B\varsubsetneq A,B\subset C$ and  $f(A), f(B) \geq
C$ then  $A\cap B\neq \emptyset$,  and for every $D$ such that
$A\cap B\subseteq D\subseteq C$ we have
\[
f(A\cap B)\geq f(D) \geq f(C).
\]
\end{itemize}

The group of permutations of $\eV$  acts on the right on the space
of functions $2_*^\eV\ra \bR$.  Also, the  group of increasing
bijections of $\bR$ acts on the left on the space of  such
functions.   It would be  interesting to have   an  estimate on
the  number of orbits of \emph{injective} $C_-$-functions
$2_*^\eV\ra \bR$. \qed
\end{remark}

%% file: tameflow12.tex
 In this   final section we will  describe a  natural tame   generalization of the subanalytic currents introduced by R. Hardt in \cite{Hardt1}.  Our  terminology    concerning currents closely  follows that of Federer \cite{Feder} (see also the more accessible \cite{Morg}). We will then use the finite volume flow technique of Harvey-Lawson  \cite{HL} for   certain tame flows on  compact  real analytic manifolds to produce interesting deformations of the DeRham complex.

Suppose $X$ is a $C^2$, oriented manifold of dimension $n$. We denote by $\Omega_k(X)$ the space of $k$-dimensional currents in $X$, i.e., the topological dual space of the space  $\Omega^k_{cpt}(X)$ of smooth, compactly supported $k$-forms on $M$.  We will denote by
\[
\lan\bullet,\bullet\ran: \Omega^k_{cpt}(X)\times \Omega_k(X)\ra \bR
\]
the natural pairing.  The boundary of a current $T\in \Omega_k(X)$ is  the $(k-1)$-current defined via the Stokes formula
\[
\lan \alpha, \pa T\ran :=\lan d\alpha, T\ran,\;\;\forall \alpha\in \Omega^{k-1}_{cpt}(X).
\]
For every  $\alpha\in \Omega^k (M)$,  $T\in \Omega_m(X)$,  $k\leq m$ define $\alpha \cap T\in \Omega_{m-k}(X)$ by
\[
\lan \alpha \cap T,  \beta \ran =\lan T,\alpha\wedge \beta \ran,\;\;\; \forall \beta\in \Omega^{n-m+k}_{cpt}(X).
\]
We have
\[
\lan \pa (\alpha \cap T),\beta\ran = \lan \,(\alpha\cap T), d\beta\ran = \lan T, \alpha\wedge d\beta\ran
\]
\[
=(-1)^k \lan T,d(\alpha\wedge \beta) -d\alpha\wedge \beta\ran = (-1)^k \lan \alpha \cap\pa T,\beta\ran +(-1)^{k+1} \lan d\alpha \cap T,\beta\ran
\]
which yields the \emph{homotopy formula}
\begin{equation}
\pa (\alpha\cap T)= (-1)^{\deg \alpha} \bigl(\, \alpha \cap \pa T-(d\alpha) \cap T\,\bigr).
\label{eq: homotop}
\end{equation}
The pair $(X, \ori_X)$, $\ori_X$ orientation on $X$, defines a current $[X,\ori_X]\in \Omega_n(X)$,  called the \emph{the current of integration along $X$}.  The current $[X,\ori_X]$ defines an inclusion
\[
\Omega^k(X)\ra \Omega_{n-k}(X),\;\;\alpha\mapsto \alpha\cap[X,\ori_X].
\]
If $X_0,X_1$ are oriented $C^2$-manifolds of dimensions $n_0$ and respectively $n_1$, and $f:X_0\ra X_1$ is a $C^2$-map, then  to every current $T\in \Omega_k(X_0)$ such that the restriction of $f$ to $\supp T$ is proper, we can associate a current $f_*T\in \Omega_{k-(n_1-n_0)}(X_1)$  defined   by
\[
\lan\beta , f_*T\ran =\lan f^*\beta, T\ran,\;\;\forall \beta \in \Omega^{k-(n_1-n_0)}_{cpt}(X_1).
\]

If $D\subset \bR^n$ is a tame $C^1$ submanifold of $\bR^n$ of dimension $k$.  Then any orientation $\ori_D$ on $D$ determines a  $k$-dimensional current $[D,\ori_D]$    via the equality
\[
\lan \alpha, [D, \ori_D]\ran:=\int_D \alpha,\;\;\forall \alpha\in \Omega^k_{cpt}(\bR^n).
\]
The integral in the right-hand side is well defined because any compact, $k$-dimensional tame set has finite $k$-dimensional Hausdorff  measure. We denote by $\eT_k(\bR^n)$ the    Abelian subgroup of $\Omega_k(\bR^n)$     generated  by currents of the form $[D,\ori_D]$ as above, and by $\eT_k^\bR(\bR^n)$ the vector space spanned by such currents. We will refer to the   currents in $\eT_k(\bR^n)$  as (integral) \emph{tame currents}.  The support of a  tame current is a tame closed set.

For every closed tame set $S\subset \bR^n$ we define
\[
\eC_k(S):=\bigl\{\, T\in \eT_k(\bR^n);\;\;\supp T,\;\;\supp\pa T\subset S\,\bigr\}.
\]
Observe that we obtain a chain complex $(\eC_\bullet(S),\pa)$
\[
\cdots \ra \eC_k(S)\stackrel{\pa}{\Lra} \eC_{k-1}(S)\ra\cdots.
\]
Suppose $C^1$-map $f: \bR^n\ra \bR^m$  whose restriction to the tame set $S\subset \bR^n$ is proper.  Then  $f$ induces  a   morphism  of chain complexes $f_\#: \eC_\bullet(S)\ra \eC_\bullet ( f(S) )$. Arguing as in the proof  of \cite[Lemma 4.3]{Hardt1} we obtain the following result.

\begin{lemma}[Lifting Lemma] Suppose $f$ is a tame  $C^1$-map  of an open neighborhood of a tame set $S$  such that the induced map $S\mapsto f(S)$ is a homeomorphism. Then the induced map $f_\#:\eC_\bullet(S)\ra \eC_\bullet(f(S))$ is an isomorphism of chain complexes. \qed
\end{lemma}

We can use the lifting lemma as in R. Hardt did in \cite{Hardt1} to show the following result.

\begin{proposition}  Suppose $S_i\in\bR^{n_i}$, $i=0,1$ are  tame sets. Then  any  proper, continuous tame map $f: S_0\ra S_1$ induces a morphism of chain complexes
\[
f_\#: \eC_\bullet(S)\ra \eC_\bullet(S_1).\proofend
\]
\label{prop: tame-chain}
\end{proposition}

We recall the construction  of this map. Denote by $\Gamma_f\subset \bR^{n_1}\ra \bR^{n_1}$ the graph of $f$. We obtain a ``roof''
\[
\begin{diagram}
\node[2]{\Gamma_f}\arrow{sw,t}{\ell}\arrow{se,t}{r}\node{}\\
\node{S_0}\node[2]{S_1}
\end{diagram}
\]
where the  left map $\ell$ and the right map $r$ are induced by the canonical projections $\bR^{n_0}\times \bR^{n_1}\ra \bR^{n_i}$.  Observe that  $\ell_0$ is a \emph{homeomorphism} and  $r$ is proper of $\Gamma_f$.  If $T\in \eC_k(S)$ we define using the Lifting Lemma
\[
f_\#T:= r_\# \ell_\#^{-1} T.
\]

We would like to explain how to  geometrically   describe  the  boundary of a tame current.   This would require  the notion of tame tube around a   tame submanifold of $\bR^n$.

Suppose $M\subset \bR^n$ is a $C^p$-manifold, $p\geq 2$.  We denote by $\eN(M)$ the normal bundle of $M$ in $\bR^n$, i.e.,
\[
\eN(M):=\bigl\{ (v,x)\in \bR^n\times M;\;\; v\perp T_xM\,\bigr\}.
\]
Observe  that if $M$ is tame, so is $\eN(M)$.   We let $\bp=\bp_M: \eN(M)\ra M$ denote the natural projection, and $\br=\br_M:\eN(M)\ra [0,\infty)$ denote the radial distance function defined by
\[
r(v,x)=|v|,
\]
where $|v|$ denotes the Euclidean length of $v$.  Observe that $\bp$ and $\bR$ are  tame if $M$ is tame.

We denote by $\exp: \eN(M)\ra \bR$ the \emph{exponential map}
\[
\exp(v,x)=x+v.
\]
Observe that if $M$ is tame, then so is  $\exp$.

A \emph{tube}  around  $M$ in $\bR^n$ is an open neighborhood $U$ of $M$ such that the exponential map   induces a $C^2$-diffeomorphism
\[
\exp: \exp^{-1}(U)\ra U.
\]
To each tube   we can associate a \emph{projection} $\pi=\pi_U: U\ra M$,   and a \emph{radial distance} function $\rho=\rho_U: U\ra [0,\infty)$  defined by the commutative diagrams
\[
\begin{diagram}
\node{\exp^{-1}(U)}\arrow[2]{e,t}{\exp}\arrow{se,b}{\bp}\node[2]{U}\arrow{sw,b}{\pi}\\
\node{}\node{M}\node{}
\end{diagram},\;\; \begin{diagram}
\node{\exp^{-1}(U)}\arrow[2]{e,t}{\exp}\arrow{se,b}{\br}\node[2]{U}\arrow{sw,b}{\pi}\\
\node{}\node{\bR_{\geq 0}}\node{}
\end{diagram}
\]
A tube  $U$ around $U$ is called \emph{tame} if $U$ is tame, and there exists a tame $C^p$ function
\[
\ve: M\ra (0,\infty)
\]
 such that
\[
\exp^{-1}(U)=\bigl\{ (v,x)\in \eN(M);\;\;|v|<\ve(x)\,\bigr\}.
\]
We will refer to $\ve$ as the \emph{width  function} of the tame tube. From \cite[Thm. 6.11, Lemma 6.12]{Co}  we obtain the following result.

\begin{theorem}[Abundance of tame tubes] Suppose $M$ is a tame $C^p$ submanifold of $\bR^n$, $p\geq 2$. Then any tame open  neighborhood $\eO$ of $M$ contains a tame tube with width function  strictly smaller than $<1$.\qed
\label{th: tube}
\end{theorem}

  Fix $p\geq 2$, and suppose $D$ is a  tame, connected,  orientable  $C^p$-submanifold  of  $\bR^n$ of dimension $m$. Fix an orientation $\ori=\ori_D$ on $D$, and a Verdier stratification of $D$ such that $D$ is a stratum.  Recall that this implies that the Whitney regularity condition is satisfied as well.   Denote by $(\dot{D}_w^i)_{1\leq i\leq \nu}$ the $(m-1)$-dimensional strata of this stratification.  Set
  \[
  \dot{D}:=\cl(D)-D,\;\;\dot{D}_w:=\bigcup_{i=1}^n \dot{D}_w^i.
  \]
   Then $\dot{D}\setminus\dot{D}_w$  is a  tame set of dimension $<m-1$.

 Choose a tube  $U_i$  (not necessarily tame) around  $\dot{D}_w^i$ in $\bR^n$ with projection $\pi_i$, and radial distance $\rho_i$  with the following properties.

 \begin{itemize}

 \item The map
\[
\pi_i\times \rho_i: U_i \cap D\ra \dot{D}_w^i\times (0,\infty),
\]
submersive.

\item  There exists a smooth function $d_i: \dot{D}^i_W\ra (0,\infty)$ such that the restriction of $\pi$ to the  set $D\cap\{\rho_i=d_i\}$ is a  locally trivial fibration,  and the set $D\cap \{ \rho_i\leq d_i\}$ is homeomorphic to the mapping cylinder of $\pi_i:D\cap\{\rho_i=d_i\}\ra \dot{D}_w^i$.

\end{itemize}

The existence of such a tube is  guaranteed  by the normal equisingularity of strata of a Whitney stratification (see \cite[Lemma II.2.3, Thm. II.5.4]{GWPL}.

Using  Theorem \ref{th: tube} we deduce that there exists a
\emph{tame} tube $W_i\subset U_i$  around $\dot{D}_w^i$.  Using
\cite[Lemm 6.12]{Co} we can even arrange that the width function
of $W_i$ satisfies $\ve_i(x)< \frac{1}{2}d_i(x)$, $\forall x\in
\dot{D}_w^i$. We will say that $W_i$ is a \emph{Whitney tube} of
$\dot{D}_w^i$ (relative to $D$).

Fix $x_i\in \dot{D}_w^i$, and set
\[
S_i:= \bigl\{y\in
D;\;\;\pi_i(y)=x_0,\;\;{\rho}_i(y)=\ve_i(\pi_iy)\,\bigr\}=(\pi_i\times {\rho}_i)^{-1}\bigl(x_i,\ve_i(x_i)\bigr).
\]
$S_i$ is a tame zero dimensional set so that it is finite.

The  restriction of $\pi_i$ to $L_i:=D\cap\{\rho_i(y)=\ve_i(\pi_iy)\,\}$ is a locally trivial fibration over $\dot{D}_w^i$ with fiber $S_i$, and  the set $D\cap\{\rho_i\leq \ve_i(\pi_iy)\}$ is homeomorphic to the mapping  cylinder of $\pi_i:L_i\ra \dot{D}_w^i$.

\begin{figure}[h]
\centerline{\epsfig{figure=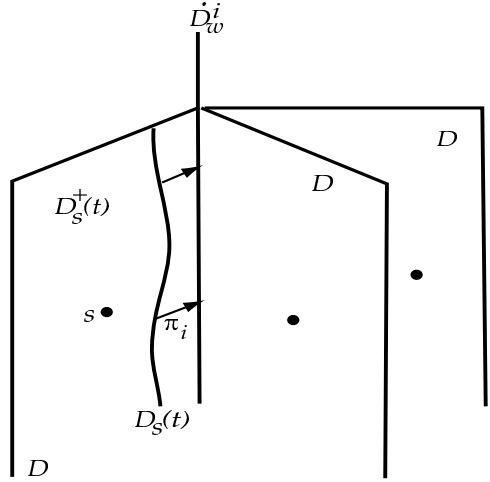,height=2.5in,width=2.5in}}
\caption{\sl   Normal equisingularity in codimension $1$.}
\label{fig: 16}
\end{figure}
For  $t\in (0,1)$, and $s\in S_i$,    we  denote by $D_s(t)^+$ the component of  $D\cap \bigl\{\, t\ve_i\leq {\rho}_i<\ve_i,\bigr\}$ containing the point $s$, and we denote by $D_s(t)$ its boundary (see Figure \ref{fig: 16}),
\[
D_s(t):= \bigl\{ y\in D_s^+(t);\;\;{\rho}_i(y)=t\ve_i(\pi y)\,\bigr\}.
\]
The orientation $\ori$ on $D$ induces orientations on the components $D^+_s(t)$, and in turn, these   define  orientations on their boundaries $D_s(t)$ via the outer-normal-first convention. The projection $\pi_i$ induce   diffeomorphisms $\pi_i: D_s(t)\ra \dot{D}_w^i$, and thus orientations $\ori_s$ on $M$. We have the following result.

\begin{theorem}[Generalized Stokes formula]
\[
\pa [D,\ori]=\sum_{i=1}^\nu\sum_{s\in S_i}[\dot{D}_w^i,\ori_s].
\]
\label{th: curr-bound}
\end{theorem}

\proof   We  choose a triangulation of $D$  such that all the open faces are tame $C^3$-manifolds and each one of them  is contained in a stratum of the Verdier stratification.  This reduces the problem to the following special case.

 Denote by $\Delta_m$ the \emph{standard    $m$-simplex}
\[
\Delta_k:= \bigl\{ (t_0,\dotsc,t_m)\in \bR_{\geq 0}^{m+1};\;\;\sum_{j=0}^mt_j=1\,\bigr\}
\]
We denote by $e_0,\dotsc, e_m$ the vertices of $\Delta_m$, and for every $I\subset \{0,\dotsc, m\}$ we denote by $\eO_I$ the  \emph{open} face spanned by the vertices $e_i$, $i\in I$.

We define a \emph{tame $m$-simplex} to be  a pair $D=(\Delta_m,f)$, where  $f: \Delta_m\ra \bR^n$ is a tame continuous map with the following properties.

\begin{itemize}

\item  The map $f$ is a homeomorphism onto its image.

\item  The images of the open faces are $C^3$-submanifolds.

\item   The collection of images of the open faces of $\Delta_m$ form a Verdier stratification of $f(\Delta_k)$.

\end{itemize}

 For a tame  $m$-simplex  $D= (\Delta_m, f)$  and $ I\subset\{0,\dotsc,m\} $  we will write
 \[
 D_f(I):= f\bigl(\,\eO_I\,\bigr).
 \]
 For simplicity, we will write
 \[
 D_f:=D_f(\{0,\dotsc, m\}),\;\; D^k_f= D(\{0,\dotsc, \hat{k},\dotsc, m\}),\;\;\bd(D):=D\setminus D_f.
 \]
 $f$  induces orientations $\ori_f$ on   $D_f$, and $\ori_k$ on $D^k_f$.   The theorem is then a consequence of the following equality
 \begin{equation}
 \pa [D_f,\ori_f]=\sum_{k=0}^m(-1)^k[D^k_f,\ori_k].
 \label{eq: pa-symp}
 \end{equation}
Denote by  $[\Delta_m]$ the tame  current defined by the standard  $m$-simplex  equipped with the orientation defined by the frame $(\be_1-\be_0,\dotsc, \be_m-\be_0)$, where $(\be_i)$ is the canonical basis  of $\bR^{m+1}$.   Using Proposition \ref{prop: tame-chain} (or rather its proof) we deduce
\[
[D_f,\ori_f]= f_\# [\Delta_m].
\]
The equality  (\ref{eq: pa-symp}) follows from the fact that $f_\#$ is a morphism of chain complexes
\[
\pa [D_f,\ori_f]=f_\#\pa[\Delta_m]. \proofend
\]

\begin{remark}(a)  To detect the boundary contributions $[\dot{D}_w^i,\ori_s]$ we do not need to know precisely a Whitney stratification of  $D$.   We look at the $(m-1)$-dimensional stratum $\dot{D}$, and orient it in some fashion using an orientation $\ori_\pa$.

Next, find a tube  $(T,\pi,\rho, \ve)$ around $\dot{D}$  and consider the   shrinking tubes
\[
T_s:=\bigl\{ z\in T;\;\;\rho(z)\leq s\ve(\pi z)\,\bigr\},\;\;s\in (0,1).
\]
We denote by $\pa T_s\pitchfork D$ the  subset of $\pa T_s$ where the intersection is transversal.

Suppose that, for all $s$, the set  $\pa T_s\pitchfork D$ projects properly via $\pi$ onto a dense  open subset $\dot{D}_{reg}$ of $\dot{D}$.  We denote by $\dot{D}^i_{reg}$ the components of $\dot{D}_{reg}$, and by $\pa T_s^i\cap D$ the  preimage  of $\dot{D}^i_{reg}$ in $\pa T_s \cap D$ via $\pi$.  Then  the components  of $\pa T_s^i\cap D$  are equipped  with  orientations as boundary components of $D\setminus T_s$,  and we denote by $n_i$ the degree of   the map
\[
\pi:\pa T_s^i\cap D\ra \dot{D}^i_{reg}.
\]
Then
\[
\pa [D]=\sum_i n_i[\dot{D}^i_{reg},\ori_\pa].
\]
(b) The  proof of the very  simple and natural  statement    of the Lifting Lemma        requires  quite sophisticated   results in geometric measure theory.  In Appendix A we present a proof of  (\ref{eq: pa-symp}) which does not use the Lifting Lemma so that the  reader  could appreciate the subtlety of this  result, and the strength of  tame geometric techniques. \qed
\label{rem: stokes}
\end{remark}

We want to apply the above facts concerning  tame currents to the study of  asymptotics of certain simple tame  flows.

Consider the standard simplex $\Delta_m$, with vertices $\be_0,\dotsc, \be_m\in \bR^{m+1}$.  This   labelling of the vertices defines a tame flow $\Phi_t$ on $\Delta_m$ and a  flow $\Phi^\pa_t$ on its boundary $\bd(\Delta_m)$. Consider the    tame, increasing homeomorphism
\[
\tau:\bR\ra (-1,1),\;\;t\mapsto \frac{t}{\sqrt{1+t^2}}.
\]
For every  tame subset of $S\subset \bR$ (i.e. a finite union of  open intervals and singletons) we define
\[
\Gamma^S:= \Bigl\{\, \bigl(\tau(t), x, \Phi_t x\,\bigr);\;\;x\in \Delta_m,\;\;t\in S\,\Bigr\}\subset[-1,1]\times \Delta_m\times \Delta_m.
\]
The projection
\[
[-1,1]\times \Delta_m\times [-1,1]\times \Delta_m\ra \Delta_m,\;\;(\tau, x,y,)\mapsto  (\tau,x)
\]
defines   homeomorphisms
\[
\Gamma^S \ra S\times \Delta_m.
\]
We orient $S$ in the canonical way as a tame subset of $\bR$.   \emph{We  fix  an orientation $\ori_m$ on  $\Int\Delta_m$.}   Using the above homeomorphism  and the orientation $\ori_m$ we obtain and orientation on the top dimensional part of $\Gamma^S$, and thus  a tame current $[\Gamma^S]$.

For simplicity  we will  set
\[
\Gamma^+=\Gamma^{[0,\infty)},\;\;\Gamma^-=\Gamma^{(-\infty,0]},\;\;\Gamma^t=\Gamma^{\{t\}}.
\]
The   boundary of $\Delta_m$ is $\Phi$ invariant, and we denote by $\Phi^\pa$ the  flow    induced by $\Phi$ on the boundary.  We orient the boundary using the orientation induced from $\ori_m$. Using the flow $\Phi_\pa$ and the orientation $\pa\ori_m$ we define  in a similar way the currents $[\Gamma^{S}_{\pa}]$, $S$ tame subset of $\bR$.

Every tame subset   $S\subset \bR$ canonically defines a tame current $[S]\subset \eT_\bullet(\bR)$.  To avoid notational overload we will continue to denote the current $[S]$ simply by $S$. and  we can extend by linearity
\[
S\mapsto [\Gamma^S],\;\;[\Gamma^S_\pa]
\]
to  the  maps
\[
\eT_\bullet(\bR) \ni [S]\mapsto[\Gamma^S],\;\;[\Gamma^S_\pa]\in \eT_\bullet([-1,1]\times \Delta_m\times \Delta_M.
\]
If $S\in \eT_\bullet(\bR)$ is a \emph{compactly supported} tame current, then
\[
\pa [\Gamma^S]= \Gamma^{\pa S}+(-1)^{\dim S} [\Gamma^S_\pa].
\]
In particular, we have
\begin{subequations}
\begin{equation}
\pa [\,\Gamma^{[0,T]}\,]=[\Gamma^T]-[\Gamma^0] -[\,\Gamma^{[0,T]}_{\pa}\,],
\label{eq: pa-poz}
\end{equation}
\begin{equation}
\pa [\,\Gamma^{[-T,0]}\,]= [\Gamma^0]-[\Gamma^T]-[\,\Gamma^{[-T,0]}_{\pa}\,].
\label{eq: pa-neg}
\end{equation}
\end{subequations}

We denote by $\eH^d$ the $d$-dimensional Haudorff measure.   If $S\subset\bR$ is a compact tame  set then both $\Gamma^S$ and $\Gamma_\pa^S$ have finite Hausdorff measures of dimensions $m+\dim S$ and $m-1+\dim S$ respectively. Arguing exactly as in the proof of Lemma \ref{lemma: hausdorf} we obtain  the following result.

\begin{lemma} (a)  As $T\ra\infty$  the current $[\,\Gamma^{[0,T]}\,]$ converges in the mass norm to $[\,\Gamma^+\,]$, and the current $[\,\Gamma^{[-T,0]}\,]$ converges in the mass norm to $[\,\Gamma^-\,]$, i.e.,
\[
\lim_{T\ra \infty} \eH^{m+1}\bigl(\,\Gamma^{[T,\infty)}\,\bigr)= 0=\lim_{T\ra \infty}\eH^{m+1}\bigl(\,\Gamma^{(-\infty,-T]}\,\bigr).
\]
(b) Similarly, as $T\ra \infty$  the current  $[\,\Gamma_\pa^{[0,T]}\,]$ converges in the mass norm to $[\Gamma^+_\pa]$ and the current $[\,\Gamma_\pa^{[-T,0]}\,]$ converges in the mass norm to $[\Gamma_\pa^-]$.\qed
\end{lemma}
If we let $T\ra \infty$ in the equalities  (\ref{eq: pa-poz}) and (\ref{eq: pa-neg}) we obtain

\begin{subequations}
\begin{equation}
\pa \Gamma^{+}=[\Gamma^\infty]-[\Gamma^0] -[\,\Gamma^{+}_{\pa}\,],
\label{eq: pa-poz1}
\end{equation}
\begin{equation}
\pa \Gamma^{-}= [\Gamma^0]-[\Gamma^{-\infty}]-[\,\Gamma^{-}_{\pa}\,],
\label{eq: pa-neg1}
\end{equation}
\end{subequations}
where $[\Gamma^\infty]$ is a tame current supported  in $\cl(\Gamma^+)\setminus \Gamma^+$ and $[\Gamma^{-\infty}]$ is a tame current   supported in $\cl(\Gamma^-)\setminus \Gamma^-$.  We will use the generalized Stokes formula to obtain a very explicit description of the currents $[\Gamma^{\pm\infty}]$. This will require some more terminology.

For every $k\in \{0,\dotsc, m\}$, denote by $W^\pm_k$ the stable/unstable variety of the stationary point $\be_k$ of the flow $\Phi$.  If  $(t_0,\dotsc, t_m)$ denote the barycentric coordinates on $\Delta_m$ then
\[
W^+_k=\bigl\{ (t_0,\dotsc, t_m);\;\;t_j=0,\;\;\forall j>k,\;\;t_i<1,\;\;\forall i<k\,\bigr\}=[\be_k,\dotsc,\be_m]\setminus [\be_{k+1},\dotsc,\be_m],
\]
\[
  W^-_k =\bigl\{ (t_0,\dotsc, t_m);\;\;t_i=0,\;\;\forall i<k,\;\;t_j<1,\;\;\forall j>k\,\bigr\}=[\be_k,\dotsc,\be_m]\setminus [\be_0,\dotsc,\be_{k-1}].
\]

\begin{proposition}
\begin{subequations}
\begin{equation}
\supp \Gamma^\infty\subset \bigcup_{\ell\geq k} W^+_\ell\times W^-_k,
\label{eq: sup+}
\end{equation}
\begin{equation}
\supp \Gamma^{-\infty}\subset \bigcup_{k\leq \ell} W^-_k\times W^+_\ell,
\label{eq: sup-}
\end{equation}
\end{subequations}
\label{prop: supp}
\end{proposition}

\proof    The  inclusion (\ref{eq: sup-}) follows from (\ref{eq: sup+}) by time reversal so it suffices to prove (\ref{eq: sup+}). Suppose $(x_\infty,y_\infty)\in \Gamma^{\infty}$.  From the curve selection property, we can find continuous  definable paths
\[
[0,\infty)\ni s\longmapsto  t_s\in\bR, \;\;x_s\in\Delta_m,
\]
such that as $s\ra \infty$ we have
\[
t_s\ra \infty,\;\; x_s\ra x_\infty, \;\; \Phi_{t_s}x_s\ra y_\infty.
\]
If $x_s$ is a stationary point for all sufficiently large $s$ then  $x_\infty=y_\infty$ and the conclusion is immediate. We assume that $x_s$ is not a stationary point for any $s\geq 0$.

Denote by  $C_s$ the portion of trajectory
\[
C_s=\bigl\{\Phi_t x_s;\;\;t\in [0,t_s]\,\bigr\}
\]
and form the strip
\[
\Sigma=\bigcup_{s\geq 0} C_s.
\]
We set
\[
C_\infty:=\cl(\Sigma)\setminus \Sigma.
\]
Observe that $C_\infty$ is a \emph{compact, $\Phi$-invariant, tame}  subset of $\Delta_m$. Moreover $x_\infty\in C_\infty$.

Denote by $f:\Delta_m\ra \bR$ the affine function uniquely determined   by the conditions
\[
f(\be_i)=i,\;\;\forall  i=0,\dotsc, m.
\]
For $\ve>0$  sufficiently small define
\[
E_i:=\bigl\{ p\in \Delta_m;\;\; |f(p)-f(\be_i)|< \ve\,\bigr\}.
\]
$E_i$ is an  open tame neighborhood of $\be_i$ and if $\ve<\frac{1}{2}$ we have
\[
E_i\cap E_j=\emptyset,\;\;\forall i\neq j.
\]
For every $i=0,\dotsc, m$ we  set
\[
A_i(s):=\bigl\{\, t\in [0,t_s];\;\;\Phi_tx_s\in E_i\,\bigr\}.
\]
Note that because $f$ is a Lyapunov function for $f$ the set $A_i(s)$ is a (possible empty)  connected subset,  for every $i$ and $s$. We  have $(m+1)$  definable families of definable sets
\[
\bigl(\,A_0(s)\,\bigr)_{s\in[0,1)},\dotsc,\bigl(\, A_m(s)\,\bigr)_{s\in[0,1)}.
\]
For every $i=0,\dotsc, m$ and every $s\geq 0$ we denote by $L_i(s)$ the length of the interval $A_i(s)$ Define the \emph{relevant set}
\[
R:=\bigl\{\, i=0,\dotsc, m;\;\;\lim_{s\ra \infty} L_i(s)=\infty \bigr\}.
\]
Note that $R\neq \emptyset$. Indeed, if $R=\emptyset$, using the fact that $x_s$ is not a stationary point, we deduce
\[
C_\infty\cap \bigl\{ \be_0,\dotsc, \be_m\,\bigr\}=\emptyset.
\]
This is impossible  since $C_\infty$ is a compact invariant subset so it must contain stationary points of $\Phi$.

\begin{figure}[h]
\centerline{\epsfig{figure=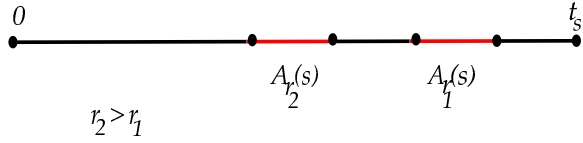, height=0.7in,width=3in}}
\caption{\sl  The relevant intervals}
\label{fig: 20}
\end{figure}

Fix $s_0>0$ such that
\[
A_r(s)\neq \emptyset,\;\;\forall s >s_0,\;\;r\in R.
\]
Since the flow $\Phi$ admits a Lyapunov function $f$ we deduce that, for  every $s>s_0$  and $r_1, r_2\in R$ such that $r_2>r_1$, the interval  $A_{r_2}(s)$  is situated to the left of the interval $A_{r_1}(s)$ (see Figure \ref{fig: 20}).

More precisely, this means
\[
t_2 <t_1,\;\;\forall t_1\in A_{r_1}(s),\;\; t_2\in A_{r_2}(s),\;\;s>s_0,\;\;r_2>r_1.
\]
 Now define
\[
\ell=\max R,\;\;k=\min R.
\]
We  deduce that $\Phi_\infty x_\infty=\be_{\ell}$, i.e. $x_\infty\in W^+_\ell$,   and $\Phi_{-\infty}y_\infty=\be_{k}$, i.e. $y_\infty\in W^-_k$.\qed

Observe that
\[
\dim  W^+_\ell \times W^-_k=  (m-\ell) +k
\]
to that  the $m$-dimensional strata of $\Gamma^\infty$ are contained in
\[
\bigcup_{k=0}^m W^+_k\times W^-_k.
\]
Hence, if we fix orientations  $\ori_k$ on $W^+_k\times W^-_k$  we obtain  an equality of the form
\begin{equation}
[\Gamma^\infty]=\sum_k \eps_{m,k}[ W^+_k\times W^-_k,\ori_k]
\label{eq: ginfty}
\end{equation}
where $\eps_{m,k}$ are some integers.      Our next goal will be to show  that we can choose the orientations $\ori_k$ in a natural way so that all the integers $\eps_{m,k}$ are equal to $1$.  This will require  a few more additional steps.

The key step towards achieving our goal is a remarkable property of the simplicial flow $\Phi_t$.   Denote  by $\eP_m$ the projection
\[
\eP: \Delta_{m}\setminus \{\be_{m}\}\ra[\be_0,\dotsc,\be_{m-1}]
\]
defined by
\[
\eP(x):=\mbox{the intersection of the line $\be_{m}x$ with the face $[\be_0,\dotsc,\be_{m-1}]$}.
\]

\begin{lemma}[Conservation of parallelism]  Suppose the two distinct points $x_0,x_1\in\Int\Delta_m$ determine a line  parallel to the face  $[\be_0,\dotsc \be_{m-1}]$.  i.e.,  they lie in a hyperplane $\{t_m=const\}$. Then for every $t\in \bR$, the line   determined by the points $\Phi_t(x_0)$ and $\Phi_t(x_1)$ is parallel with the line determined by the points $x_0,x_1$ and with the line determined  by $\eP_m(\Phi_t(x_0)\,)$ and $\eP_m(\Phi_t(x_1))$.
\label{lemma: para}
\end{lemma}

\proof  We argue by induction over $m$. For $m=0,1$ the statement is trivially true. We assume it is true for $\Delta_m$ and we prove its validity for $\Delta_{m+1}$. We denote by $S$ the set $\{x_0,x_1\}$, and we set for simplicity $\eP=\eP_{m+1}$.

The set  $S\subset \Int\Delta_{m+1}\setminus \{\be_{m+1}\}$  is contained in a hyperplane $\{t_{m+1}=c\}$, where  $c\in [0,1)$.       The restriction of $\eP$ to  $\Int\Delta_{m+1}\cap \{t_{m+1}=c\}$ defines an affine  map
\[
\Int\Delta_{m+1}\cap \{t_{m+1}=c\}\ra \Int[\be_0,\dotsc,\be_m],
\]
such  that   for any $y_0,y_1\in  \Int\Delta_{m+1}\cap \{t_{m+1}=c\}$, the line determined  by $y_0,y_1$ is  parallel with the line determined by $\eP(y_0)$ and $\eP(y_1)$, and
\begin{equation}
\dist(y_0,y_1)=(1-c)\dist\bigl(\,\eP(y_0),\eP(y_1)\,\bigr).
\label{eq: dilate}
\end{equation}
From the iterated cone   description   of $\Phi$ we deduce that  $\eP\circ\Phi_t=\Phi_t\circ\eP$, $\forall t\in \bR$.  The lemma now follows from the inductive assumption.
\qed

\begin{figure}[h]
\centerline{\epsfig{figure=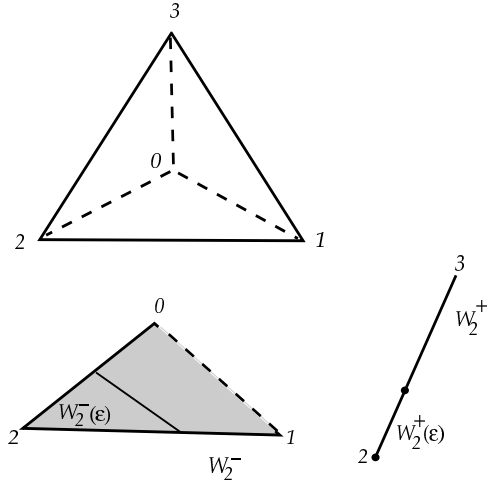, height=2.5in,width=2.5in}}
\caption{\sl  Organizing the (un)stable varieties of a simplicial flow.}
\label{fig: 19}
\end{figure}

For $\ve\in (0,1)$ we define an $\ve$-neighborhood of $\be_k \in W^\pm_k$ (see Figure \ref{fig: 19})
\[
W^\pm_k(\ve):=\bigl\{\, w\in W^\pm_k;\;\; |t_k(w)-1|<\ve \,\bigr\}.
\]

Let  $k=0,\dotsc, m$, and consider a point $w_+$ in the relative interior   of the   stable variety of $\be_k$,
\[
w_+\in \Int W^+_k=\Int[\be_k,\dotsc,\be_m]
\]
For $\ve_->0$ we  denote by $\eN_k^-(w_+,\ve_+)$ the translate of $W^-_k(\ve_-)$ at $w$
\[
\eN_k^-(w_+,\ve_-):=(w_+-\be_k)+ W^-_k(\ve_-).
\]
For $\ve_->0$ sufficiently small,   this  set is contained in $\Delta_m$. We denote by $\eN_k^-(w_+,\ve_-)_{reg}$ the regular (top dimensional) part of $\eN_k^-(w_+)$.

If we denote by $V_k^-(w_+)$ the  affine $k$-plane   through $w_+$, and parallel to the face $[\be_0,\dotsc,\be_k]$, then $V_k^-(w_+)$ intersects   $\Int W^+_k$ transversally at $w_+$, and for $\ve>0$ sufficiently small, $\eN_k^-(w_+,\ve_-)$ is a neighborhood of $w_+$ in $V_{k+1}(w)\cap \Delta_m$.

Similarly, for $w_-\in \Int W^-_k$, and $\ve_+>0$ we denote by $\eN^+_k(w_-,\ve_+)$ the translate of $W^+_k(\ve_+)$ at $w_-$,
\[
\eN_k^+(w_-,\ve_+):=(w_--\be_k)+ W^+_k(\ve_+).
\]
If we denote by $V^+_k(w_-)$ the affine $(m-k)$-plane  through $w$ and parallel to the face $[\be_k,\dotsc,\be_m]$, then $V^+_k(w_-)$ intersects  $\Int W^-_k$ transversally at $w$,  and for $\ve_+>0$ sufficiently small, $\eN^+_k(w_-,\ve_+)$ is a neighborhood of $w_-$ in $V^+_k(w_-)\cap \Delta_m$. We denote by $\eN^+_k(w_-,\ve_+)_{reg}$ the regular (top dimensional) part of $\eN^+_k(w_-,\ve_+)$.

\begin{proposition}   Let $k\in \{1,\dotsc, m-1\}$. Then there exist   a definable function
\[
T_k: \Int W^+_k\times \Int W^-_k\times (0,1)\times (0,1)\ra \bR,
\]
\[
(w_+,w_-,\ve_-,\ve_+)\mapsto T_k(w_+,w_-,\ve_-,\ve_+),
\]
such that, for all $(w_+,w_-,\ve_-,\ve_+)\in \Int W^+_k\times \Int W^-_k\times (0,1)\times (0,1)$,  and all \newline $t>T_k(w_+,w_-,\ve_-,\ve_+)$, the normal slice
\[
\eN_k(w_+,w_-,\ve):=\eN_k^-(w_+,\ve_-)_{reg}\times \eN^+_k(w_-,\ve_+)_{reg}
\]
intersects $\Gamma^t_{reg}$, the regular part of the graph of $\Phi_t$,   at a unique point. Moreover, the intersection at that point is transversal in $\Int\Delta_m\times\Int\Delta_m$.
\end{proposition}

\proof Observe that
\[
(x,y)\in \Gamma_t\cap \eN_k^-(w_+,\ve_-)_{reg}\times \eN^+_k(w_-,\ve_+)_{reg}
\]
if and only if
\[
 y\in \Phi_t\bigr(\eN_k^-(w_+,\ve_-)_{reg}\,\bigr)\cap \eN_k^+(w_+,\ve_+)_{reg},\;\;x=\Phi_{-t}y.
\]
Moreover
\[
\Gamma_t\pitchfork \eN_k^-(w_+,\ve_-)_{reg}\times \eN^+_k(w_-,\ve_+)_{reg} \Longleftrightarrow \Phi_t\bigr(\,\eN_k^-(w_+,\ve_-)_{reg}\,\bigr)\pitchfork \eN_k^+(w_-,\ve_+)_{reg}.
\]
Set $w_+(t):=\Phi_t w_+$. From the conservation of parallelism we deduce that  the set $\Phi_t\bigr(\,\eN_k^-(w_+,\ve)_{reg}\,\bigr)$ is an open subset of the affine plane $V^-_k(\, w_+(t)\,)$. In particular, if $\Phi_t\bigr(\,\eN_k^-(w_+,\ve)_{reg}\,\bigr)$ intersects  $\eN_k^+(w_+,\ve)_{reg}$, it does so transversally.

\begin{figure}[h]
\centerline{\epsfig{figure=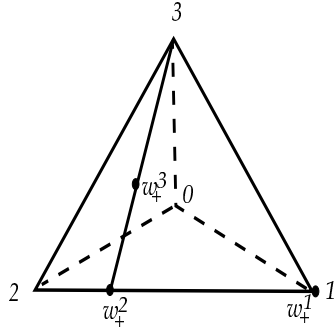, height=1.7in,width=1.7in}}
\caption{\sl  The sequence of shadows $w^m_+,\dotsc, w^k_+$, when $m=3$ and $k=1$.}
\label{fig: 21}
\end{figure}

 To understand the region $\Phi_t\bigr(\,\eN_k^-(w_+,\ve_-)\,\bigr)$ better, consider the projections
\[
\eP_j:[\be_0,\dotsc,\be_j]\setminus\{ \be_j\}\ra [\be_0,\dotsc,\be_{j-1}],
\]
\[
\eP_j(x):=\mbox{the intersection of the line $\be_j x$ with the face $[\be_0,\dotsc,\be_{j-1}]$}.
\]
We obtain a   sequence of points
\[
w_+^m,\dotsc, w^{k+1}_p,\;\;w^k_+
\]
defined inductively as
\[
w^m_+=w^+,\;\; w_+^{j-1}=\eP_j(w_+^j).
\]
Observe that (see Figure \ref{fig: 21})
\[
w^j_+\in \Int[\be_k,\dotsc,\be_j],\;\;\forall j>k,\;w^k_+=\be_k.
\]
Denote by $\eS$ the composition
\[
\eS=\eP_{k+1}\circ\cdots \circ \eP_m.
\]

\begin{figure}[h]
\centerline{\epsfig{figure=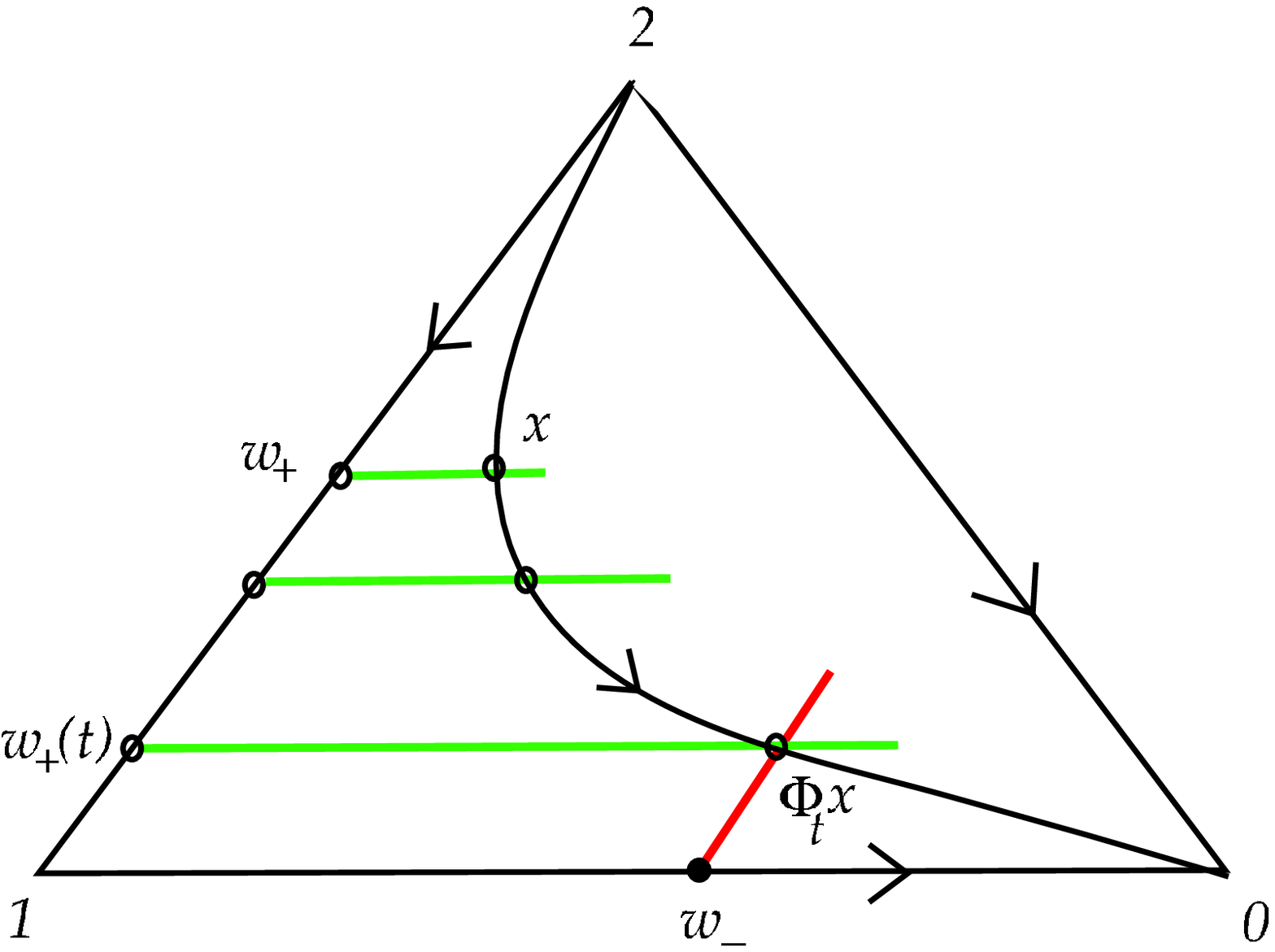, height=2.3in,width=2.3in}}
\caption{\sl  $\Phi_t(\eN^-_k(w_+,\ve_-)\,\bigr)$  is   depicted as the moving horizontal segment that is increasing in length.}
\label{fig: 22}
\end{figure}

From the conservation of parallelism and the iterated cone  description of $\Phi$ we deduce
\[
\eS\bigl(\,\eN^-_k(w_+,\ve_-)\,\bigr)=W^-_k(c\ve_-),
\]
for some $c>1$.  We set $w_+(t)=\Phi_tw_+$.  Note that
\[
\Phi_t\circ \eS=\eS\circ \Phi_t
\]
Using the conservation of parallelism  we deduce  that the map
\[
\eS: \Phi_t\bigl(\,\eN^-_k(w_+,\ve_-)_{reg}\,\bigr)\ra \eS\Phi_t\bigl(\,\eN^-_k(w_+,\ve_-)_{reg}\,\bigr)
\]
is a homothety. Now observe that
\[
\eS\Phi_t\bigl(\,\eN^-_k(w_+,\ve_-)_{reg}\,\bigr)=\Phi_t\eS\bigl(\,\eN^-_k(w_+,\ve_-)_{reg}\,\bigr)=\Phi_tW^-_k(c\ve_-)_{reg}.
\]
We conclude that
\[
\Phi_t\bigl(\,\eN^+_k(w_+,\ve)\,\bigr)=\eN^+_k(w_+(t), \ve_-(t)),
\]
and $\ve_-(t)\ra 1$ as $t\ra \infty$ (see Figure \ref{fig: 22}).

Denote by   $T_k=T_k(w_+,w_-,\ve_-,\ve_+)$  the   smallest real number $T>-1$  with the property that
\[
w_+(t)\in W^+_k(\ve_+)\;\; \mbox{and}\;\;\ve_-(t) > 1-t_k(w_-),\;\;\forall t>T.
\]
If $t>T_k(w_+,w_-,\ve_-,\ve_+)$, then affine   $k$-dimensional piece $\Phi_t\bigl(\,\eN^-_k(w_+,\ve_-)_{reg}\,\bigr)$ intersects the  affine $(m-k)$-dimensional piece  $\eN^+_k(w_-,\ve_+)$ at a unique point (see Figure \ref{fig: 22})
\[
y_t= (w_+(t)-\be_k)+ (w_--\be_k) +\be_k.
\]
If we think of $\be_k$ as the origin of our affine space then we can  rewrite the above equality in the simpler form
\[
y_t=w_+(t)+ w_-.
\]
Hence,   the normal slice  $\eN_k(w_+,w_-,\ve):=\eN_k^-(w_+,\ve_-)_{reg}\times \eN^+_k(w_-,\ve_+)_{reg}$  intersects  $\Gamma^t_{reg}$ at a single point
\begin{equation}
(x, y)= \bigl(\,x(w_+,w_-,t), \,y(w_+,w_-,t)\,\bigr)= \bigl(\,\Phi_{-t}(\,w_+(t)+ w_-),\, w_+(t)+ w_-\,\bigr),
\label{eq: inter}
\end{equation}
and the intersection is transversal. \qed

Observe  that the map
\[
(w_+,w_-,\ve_-,\ve_+)\mapsto T_k(w_+,w_-,\ve_-,\ve_+)
\]
is upper semicontinuous in the variables  $(w_-,w_+)$,  i.e. if
\[
T_k(w_+,w_-,\ve_-,\ve_+)<T,
\]
then there exist open neighborhoods $U_\pm$ of $w_\pm$ in $\Int W^\pm_k$ such  that, for all
\[
T_k(u_+,u_-,\ve_-,\ve_+) <T,\;\;\forall (u_+,u_-)\in U_+\times U_-.
\]
Given $T$ and $U_\pm$ as above we obtain  for every $t>T$ a tame continuous map given by (\ref{eq: inter}),
\[
U_+\times U_-\stackrel{\psi_t}{\longmapsto} \bigl( x(u_+,u_-,t), y(u_+,u_-,t)\,)\in \Gamma^t_{reg},
\]
which is a homeomorphism onto its image. $\Gamma^t_{reg}$  admits a natural orientation induced by the homeomorphism
\[
\Gamma^t_{reg}(x,\Phi_t x)\mapsto x\in\Int\Delta_m.
\]
We conclude that the homeomorphism $\psi_t$  induces an orientation $\ori=\ori_t$ on $U_+\times U_-$  which is independent of $t>T$.  For a different pair of points $(w_+',w_-')$, and corresponding neighborhood $U'_+\times U'_-$, the orientation $\ori'$  on $U_+'\times U_-'$  obtained by the above procedure coincides with $\ori$ on the overlap. We obtain in this fashion an orientation $\ori_k$ on $\Int W^+_k\times\Int  W^-_k$. We would like to give a more explicit description of $\ori_k$.

To achieve this, we place $w_\pm$ very close to $\be_k\in W^\pm_k$, and we choose $\ve_\pm$ relatively large, say $\ve_\pm=\frac{1}{2}$.   Then $T_k(w_+,w_-,\ve_-,\ve_+)<0$, and it suffices to understand the homeomorphism $\psi_t$, $t=0$.  In this case  the  equation (\ref{eq: inter}) takes the simple form
\[
(u_-,u_+)\mapsto (u_-+u_+, u_-+u_+).
\]
Thus   $\ori_k$ is the orientation with the property that the map
\[
T_{\be_k} W^+_k\times T_{\be_k} W^-_k\ra  T\Delta_m,\;\;(u_+,u_-)\mapsto u_++u_-
\]
is orientation preserving, where we recall that  we have fixed an orientation $\ori_m$ on $\Int \Delta_m$.

Let us observe that we have a natural tube  $\eT$  around $\{1\}\times \Int W^+_k \times \Int W^-_k$ inside $\bR\times \Int \Delta_m\times \Delta_m$  defined as  follows.

\smallskip

\noindent $\bullet$ Fix continuous definable  functions $\ve_\pm: \Int W^\mp_k\ra (0,1)$ such that
\[
\eN^\mp_k(w_\pm,\ve_\mp)\subset \Delta_m.
\]
\noindent $\bullet$  Set $d: W^+_k \times \Int W^-_k\ra (0,1)$,
\[
d(w_+,w_-)=\frac{2T\bigl(\,w_+,w_-,\ve_-(w_+),\ve_+(w_-)\,\bigr)}{\sqrt{1+4T\bigl(\,w_+,w_-,\ve_-(w_+),\ve_+(w_-)\,\bigr)^2}}=\tau\bigl(\,2T\bigl(\,w_+,w_-,\ve_-(w_+),\ve_+(w_-)\,\bigr)\,\bigr) .
\]
\noindent $\bullet$  Define
\[
\eT=\bigcup_{(w_-,w_+)\in \Int W^+_k \times \Int W^-_k}[d(w_-,w_+),2]\times \eN^-_k(w_+,\ve_-(w_+)\,)\times \eN^+(w_-,\ve_+(w_-)\,).
\]
\noindent $\bullet$  Define $\pi: \eT\ra \{1\} \Int W^+_k \times \Int W^-_k$ by $\pi(t, x)=\bigl(\,1 \pi_k^+(x),\pi_k^-(x)\,\bigr)$  where $\pi_k^\pm$ is the projection onto the  affine plane spanned by $W^\pm_k$   and parallel with the plane spanned by $W^\mp_k$. The fiber of $\pi$ over $(w_+,w_-)$ is the PL ball
\[
\eB(w_+,w_-):=[d(w_-,w_+),2]\times \eN^-_k(w_+,\ve_-(w_+)\,)\times \eN^+(w_-,\ve_+(w_-)\,)
\]

\smallskip

Then $\Gamma^t_{reg}$ intersects  the  boundary of the ball $\eB(w_+,w_-)$  exactly once, in the region
\[
\{d(w_+,w_-)\}\times \eN^-_k(w_+,\ve_-(w_+)\,)_{reg}\times \eN^+(w_-,\ve_+(w_-)\,)_{reg}.
\]
That intersection is transversal.  Using the generalized Stokes formula, Remark \ref{rem: stokes}(a), and the equality (\ref{eq: ginfty}) we obtain the following result.

\begin{theorem} Consider and affine $m$-simplex $\Delta_m=[\be_0,\dotsc,\be_m]$, and an orientation $\ori_m$ on its relative interior. Denote by $\Phi$ the  simplicial flow    determined by the above ordering of the vertices of $\Delta_m$.  Equip  the cartesian product $W^+_k\times W^-_k$ with the orientation $\ori_k^+$  defined by the property  that the map
\[
W_k^+\times W_k^-\ni(w_+,w_-)\mapsto w_+ +w_--\be_k
\]
is an orientation preserving map from $W_k^+\times W_k^-$ to the affine plane spanned by $\Delta_m$ and equipped with the orientation $\ori_m$.  Then
\[
\sum_k[W_k^+\times W_k^-,\ori_k^+]-[\Gamma^0] =\pa [\Gamma^{[0,\infty)}\,] +[\Gamma^{[0,\infty)}_\pa].
\]
Similarly,  we define an orientation $\ori_k^-$ on $W_k^-\times W_k^+$ with the property that the switch map
\[
(W_k^-\times W_k^+,\ori_k^-)\ra (W_k^+\times W_k^-,\ori_k^+)
\]
is orientation preserving.  Then
\[
[\Gamma^0] -\sum_k[W_k^-\times W_k^+,\ori_k^-]=\pa [\Gamma^{(-\infty,0]}\,] +[\Gamma^{(-\infty,0]}_\pa].\proofend
\]
\label{th: asy-symp}
\end{theorem}

We would like to use the above result, and the technique of  Harvey-Lawson \cite{HL}  to construct  a canonical chain homotopy between the DeRham complex of a compact, real analytic manifold,  and the simplicial chain complex associated to a tame triangulation of the manifold. Before we do this we would like to clarify  a few issues.

Suppose $M$ is a compact, orientable, real analytic manifold without boundary. We assume  $M$ is embedded in some Euclidean space $E$. Let $m:=\dim M$. We fix an orientation $\ori_M$ on $M$ and a tame triangulation of $M$,  which  is a pair   $(\eK, \Delta)$, where $\eK$ is a CSC and $\Delta$ is a tame homeomorphism
\[
\Delta:|\eK|\ra M.
\]
We assume that the restriction of $\Delta$ on the relative interiors of the faces of $\eK$ is $C^2$.

For every (combinatorial)  face $S\in\eK$  we denote by $\Delta_S$  the image of the closed face $|S|\subset|\eK|$ via the homeomorphism $\Delta$, and by $\Delta_S^\circ$  the image via $\Delta$ of the relative interior of $|S|$. We fix orientations $\ori_S$ on $\Delta_S^\circ$ so that  the orientations on the top dimensional faces   coincide with the  orientations induced by the orientation of $M$.

We denote by $C_j(\eK, M)$ the subgroup of tame integral currents $\eT_j(M)$ spanned by $[\Delta_S, \ori_S]$, $\# S= j+1$.  The  chain complex $(C_\bullet(\eK, M),\pa)$ is isomorphic to the  simplicial chain complex  associated to $\eK$. We form a cochain complex $(C^k_\bR(\eK), \delta)$  by setting
\[
C^k_\bR(\eK):= C_{m-k}(\eK)\otimes \bR,\;\;\delta= \pa.
\]
We see that this cochain complex is naturally  isomorphic to the  simplicial chain complex with real coefficients  determined by $\eK$.

\begin{figure}[h]
\centerline{\epsfig{figure=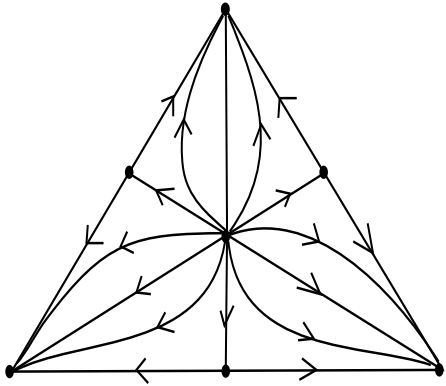, height=2.0in,width=2.3in}}
\caption{\sl  The Stieffel flow on a triangle.}
\label{fig: 23}
\end{figure}

Consider the  barycentric subdivision $D\eK$ of $\eK$.  We denote by $b_S$ the vertices of $D\eK$. We have a canonical homeomorphism $|D\eK|\ra |\eK|$, and we  thus  a tame homeomorphism
\[
\Delta':|D\eK|\ra M.
\]
We set
\[
x_S:=\Delta'(b_S)\in \Delta_S\subset M.
\]
 The simplicial complex $D\eK$ is the nerve of the poset $(\eK, \subset)$. We have a natural admissible  function on the poset $\eK$,
\[
f:\eK\ra \bZ,\;\; f(S)=\dim S.
\]
This defines a dynamical ordering of $D\eK$, and thus a  tame flow $\Psi$ on $|D\eK|$ and, via $\Delta'$,  a conjugate tame flow $\Phi$ on $M$. We will refer to these flows as the \emph{Stieffel flows} determined by a triangulation of $M$.  The  simplices $\Delta_S$ are invariant subsets of the Stieffel flow on $M$. The phase portrait of the Stieffel flow on a $2$-simplex $\Delta_S$ is depicted  in Figure \ref{fig: 23}.

From the definition, it follows immediately that the only stationary points of the flow $\Phi$ are the barycenters $x_S$, and the unstable variety of  $x_S$ is the open face $\Delta_S^\circ$. It is equipped with the orientation $\ori_S$.

If $\Delta_S$ is a face of dimension $k$ with barycenter  $x_S$, then we define the \emph{normal star} of $x_S$ to be the  union  of all $(m-k)$ simplices of the barycentric subdivision whose vertices are barycenters of faces  $T\supseteq S$. We denote by $\St^\perp(x_S)$  the normal  star. It is a tame $(m-k)$-manifold with boundary. Its boundary is called the normal link of $x_S$, and it is denoted by $\Lk^\perp(x_S)$.

For a barycenter $S$, we denote by $\eM_S$ the collection of  maximal faces of $\eK$ which contain $S$.  Each $\Sigma\in \eM_S$ determines a top dimensional face $\Delta_\Sigma\supset \Delta_S$.  $\Delta_\Sigma$ is a $\Phi$ invariant set and we denote by $\Phi^\Sigma$ the restriction of $\Phi$ to $\Delta_\Sigma$.   The barycenter $x_S$ is a  stationary point of $\Phi^\Sigma$.    We denote by $W_{S,\Sigma}^+$ the  stable variety  of $x_S$ in $\Delta_\Sigma$  with respect to $\Phi^\Sigma$.        We have the equality
\[
W^+_{S,\Sigma}=\bigl(\, \St^+(x_S)\setminus \Lk^\perp(x_S)\,\bigr)\cap \Delta_\Sigma.
\]
We deduce that the stable variety of  $x_S$ in $M$ with respect to  the flow $\Phi$ is
\[
W_S^\perp= \bigcup_{\Sigma\in \eM_S}W^+_{S,\Sigma}=\St^\perp(x_S)\setminus \Lk^\perp(x_S).
\]
We have a natural homeomorphism
\[
h_S: \St^\perp(x_S)\times  \Delta^\circ_S\ra \eT_S,
\]
where $\eT_S$ is  a tubular  neighborhood of $\Delta_S^\circ$.   Using the orientation $\ori_S$ on $\Delta_S^\circ$, and the  orientation $\ori_M$ on $\eT_S$, we obtain an orientation $\ori^\perp_S$ on $\St^\perp(x_s)$ such that
\[
\ori_S^\perp\times  \ori_S\stackrel{h_S}{\longmapsto}\ori_M.
\]
This defines  an  orientation $\ori_S^\perp$ on $W_S^+$.

Define again $\tau(t)=\frac{t}{\sqrt{1+t^2}}$.
\[
\Gamma_M^\pm=\bigl\{\, (\tau(t),x,\Phi_tx\,)\in [-1,1]\times M\times M;\;\;\pm t\geq 0,\,\bigr\},
\]
\[
\Gamma_M^t=\bigl\{\,(x,\Phi_tx)\in M\times M\,\bigr\}.
\]
Denote by $\eM$  the set of maximal simplices of $\eK$.  For $\Sigma \in\eM$ we define
\[
\Gamma_\Sigma^\pm=\bigl\{\, (\tau(t),x,\Phi_tx\,)\in [-1,1]\times \Delta_\Sigma\times \Delta_\Sigma;\;\;\pm t\geq 0,\,\bigr\},
\]
\[
\Gamma_\Sigma^t=\bigl\{\,(x,\Phi_tx)\in M\times M\,\bigr\}.
\]
As before, these tame sets are equipped with natural orientations and define currents $[\Gamma_M^\pm]$, $[\Gamma_\Sigma^\pm]$. Moreover
\[
[\Gamma_M^\pm]=\sum_{\Sigma\in \eM}[\Gamma_\Sigma^\pm].
\]
Using Theorem \ref{th: asy-symp}  and the fact that
\[
\sum_{\Sigma\in \eM}\pa [\Delta_\Sigma,\ori_M]=\pa [M,\ori_M]=0,
\]
we deduce
\[
\pa [\Gamma_M^+]=[\Gamma^\infty_M]-[\Gamma^0_M]=\sum_{S\in \eK}[W^+_S,\ori_S^\perp]\times [ W_S^-,\ori_S]-[\Gamma_0]
\]
\[
=\sum_{S\in\eK}[\St^\perp(x_S),\ori_S^\perp]\times [\Delta_S,\ori_S]-\Gamma^0,
\]
and similarly,
\[
\pa[\Gamma^-]=[\Gamma^0]-[\Gamma^{-\infty}=[\Gamma^0_M]-\sum_{S\in \eK}(-1)^{\dim S(m-\dim S)} [\Delta_S,\ori_S]\times [\St^\perp(x_S),\ori_S^\perp].
\]
Now we can start using the  formalism of kernels developed  by Harvey-Lawson in \cite{HL}.  For the reader's convenience we briefly  recall  it here.

Suppose that we are given a roof, i.e a diagram of the form
\[
\begin{diagram}
\node{}\node{Y}\arrow{sw,t}{f_0}\arrow{se,t}{f_1}\node{}\\
\node{X_1}\node{}\node{X_0}
\end{diagram},
\]
where $X_0,X_1,Y$  are oriented smooth manifolds, and $f_0, f_1$ are smooth maps. Assume $K$ is a $k$-dimensional \emph{kernel} for this  roof, i.e.  a $k$-dimensional current in $Y$  such that $f_0$ is proper over $\supp K$. Then we obtain a linear map
\[
K_\#: \Omega^m(X_1)\ra \Omega_{k-m}(X_0),\;\;K_\#\alpha= (f_0)_*\bigl(\,(f_1^*\alpha)\cap K\,\bigr).
\]
We have the following homotopy formula
\begin{equation}
(\pa K)_\# \alpha=K_\#(d\alpha)+ (-1)^m\pa K_\sharp\alpha,\;\;\forall \alpha\in \Omega^m(X_1).
\label{eq: homotop1}
\end{equation}
Indeed,
\[
K_\#(d\alpha)= (f_0)_* \bigl(\, d (f_1^*\alpha)\cap K\,\bigr)\stackrel{(\ref{eq: homotop})}{=} (f_0)_*\bigl(\, f_1^*\alpha\cap\pa K -(-1)^m\pa( f_1^*\alpha\cap K)\,\bigr)
\]
\[
=(\pa K)_\# \alpha -(-1)^m \pa K_\# \alpha.
\]
We can rewrite this in operator form
\begin{equation}
(\pa K)_\#= K_\#\circ d+(-1)^m \pa \circ K_\#.
\label{eq: homot-oper}
\end{equation}

We apply this formalism to the roof
\[
\begin{diagram}
\node{}\node{\bR\times M\times M}\arrow{sw,t}{\pi_0}\arrow{se,t}{\pi_1} \node{}\\
\node{M}\node[2]{M}
\end{diagram},\;\; \pi_i(t,x_0,x_1)= x_i,\;\;i=0,1,
\]
and the currents
\[
[\Gamma^-_M]\in \Omega_{m+1}(\bR\times M\times M),\;\;[\Gamma^0_M],[\Gamma^{-\infty}_M]\in \Omega_m(\bR\times M\times M).
\]
Clearly $[\Gamma^-_M]$, $[\Gamma^0_M]$  and $[\Gamma^{-\infty}_M]$ are  kernels for this roof, and
\[
\pa[\Gamma^-_M]=[\Gamma^0_M]-[\Gamma_M^{-\infty}].
\]
Since $M$ does not have boundary we deduce $\pa [\Gamma_M^{-\infty}]=0$. We  obtain operators
\[
[\Gamma^0_M]_\#,\;\;[\Gamma^{-\infty}_M]_\#:\Omega^j(M)\ra \Omega_{m-j}(M),
\]
and
\[
[\Gamma_M^-]_\#:\Omega^j(M)\ra \Omega_{m+1-j}(M),
\]
satisfying for every $\alpha\in \Omega^j(M)$ the equalities
\begin{equation}
[\Gamma_M^0]_\# \alpha-[\Gamma_M^{-\infty}]_\#\alpha= [\Gamma^-_M]_\#d\alpha+(-1)^j\pa [\Gamma^-_M]_\#\alpha,
\label{eq:  homtop}
\end{equation}
and
\begin{equation}
[\Gamma_M^0]_\sharp\alpha=(-1)^{j+1}[\Gamma_M^{-\infty}]_\#d\alpha=(-1)^{j+1}\pa [\Gamma_M^{-\infty}]_\#\alpha.
\label{eq: chain-morph}
\end{equation}
Observe that
\[
[\Gamma_M^0]_\#\alpha=\alpha\cap[M,\ori_M],
\]
and
\[
[\Gamma_M^{-\infty}]_\#\alpha= \sum_{S\in\eK,\,\dim S=m-k}(-1)^{k(m-k)}\bigl\lan\alpha,\, [\St^\perp(x_S),\ori_S^\perp]\,\bigr\ran [\Delta_S,\ori_S].
\]
The equality (\ref{eq: chain-morph}) shows that the maps
\[
[\Gamma_M^0]_\#,\;\; [\Gamma_M^{-\infty}]_\#:(\Omega^\bullet(M),d)\ra (\Omega_{m-\bullet}(M),\pa)
\]
are morphisms of chain complexes, while the equality (\ref{eq: homtop}) shows that they are chain homotopic. The morphism $[\Gamma_M^0]_\sharp$ is one-to-one, while  the image of the morphism  $[\Gamma_M^{-\infty}]_\#$ is the simplicial complex $(C^\bullet_\bR(\eK),\delta)$.  We have  obtained the following result.

\begin{theorem} The Stieffel flow associated  to a tame triangulation of a compact, real analytic manifold without boundary determines a chain homotopy between the DeRham complex and the   simplicial chain complex  with real coefficients associated to that triangulation. \qed
\label{th: derham-simpl}
\end{theorem}

\begin{remark} In the above proof, the tameness assumption is needed only to guarantee that the flow $\Phi$ is a finite volume flow on $M$.     We can reach this conclusion under weaker assumptions. We know that  the flow $\Psi$ on the geometric realization $|\eK|$ is tame, and thus has finite volume. If  the homeomorphism $\Delta: |\eK|\ra M$ happens to be bi-Lipschitz  then the flow  $\Phi$ will also     have finite volume.

If $M$ is only a smooth, then the triangulation procedure     employed by H. Whitney in \cite[Chap. IV.B]{Whit2} produces   triangulations with this property. In this case, for every $t\in \bR$ the map $\Phi_t:M\ra M$ is  bi-Lipschitz because  the conjugate map $\Psi_t:|\eK|\ra |\eK|$ is such. Then, for every smooth form $\alpha\in \Omega^k(M)$  the pullback  $\Phi_t^*\alpha$ is  defined almost everywhere  and it is a form with $L^\infty$-coefficients. Moreover
\[
[\Gamma_M^t]_\#\alpha =\Phi_t^*\alpha\cap[M,\ori_M].
\]
The current $[\Gamma_M^t]$ converges in the flat norm to $[\Gamma^{-\infty}]$ as $t\ra -\infty$, and we deduce that  $\Phi_t^*\alpha\cap[M,\ori_M]$ converges in the sense of currents to
\[
[\Gamma_M^{-\infty}]_\# \alpha= \sum_{S\in\eK,\,\dim S=m-k}(-1)^{k(m-k)}\bigl\lan\alpha,\, [\St^\perp(x_S),\ori_S^\perp]\,\bigr\ran [\Delta_S,\ori_S].
\]
Intuitively, this means that as $t\ra -\infty$ the form $\Phi_t^*\alpha$   begins to concentrate near the barycenters $x_S$ and in the normal planes to the face $\Delta^\circ_S$.\qed

\end{remark}

%% file: tameflow-a.tex
In this  appendix we want to present   a proof  of  the Stokes formula (\ref{eq: pa-symp}) which  does not use  the    advanced results of geometric measure theory in \cite{Hardt, Hardt1}.  We continue to use the notations in the proof of Theorem \ref{th: curr-bound}.

 We begin by  constructing a system of tubes  $(T_I,\pi_I,\rho_I \ve_I)$ around the  open faces  $D_f(I)$ of $D$.  As in the proof of \cite[Prop. 7.1]{BK}, for every $\theta_0\in (0,\frac{\pi}{2})$ we can choose  the  tube system  so that the following additional  conditions are satisfied
 \[
 T_I\cap T_J \neq \emptyset \Longleftrightarrow I\subset  J\;\;{\rm or}\;\; J\subset I,
 \]
 \[
 \forall y\in D_f\cap T_I\cap T_J: \bigl|\,\measuredangle\bigl(\,\nabla\rho_I(y),\nabla\rho_J(y)\,\bigr)-\frac{\pi}{2}\,\bigr|<\theta_0.
 \]
 Define
 \[
\bar{\ve}_I: T_I\ra (0,\infty),\;\;\bar{\ve}_I(y):=\ve_I\bigl(\,\pi_I(y)\,\bigr),\;\; \eT:=\bigcup_{\# I\leq m} T_I,
 \]
so that $\eT$ is an open neighborhood of $\bd(D)$.

 As in the proof of  \cite[Prop. 7.1]{BK}, we fix  a $C^3$  definable function
 \[
 h:[0,\infty]\ra [0,1],\;\;h(t)=\begin{cases}
 t &t\leq 1/3\\
 1 & t>2/3
 \end{cases}
 \]
 and define $\ell_I, \ell: \bR^n\ra [0,\infty)$  by
 \[
 \ell(x)= \begin{cases}  h\bigl(\,\frac{\rho_I}{\bar{\ve}_I}\,\bigr) & y\in T_I\\
 &\\
 1 &y\not\in T_I
 \end{cases},\;\;\ell=\prod_{\# I\leq m} \ell_I.
 \]
 We will say that $\ell$ is the \emph{boundary profile} associated to the isolating system of tubes.  As explained in the proof of \cite[Prop. 7.1]{BK} the profile $\ell$ satisfying the following properties.

 \begin{enumerate}

 \item[(${\bf P}_1$)] $\ell^{-1}(0)= \bd(D)$.

 \item[(${\bf P}_2$)] $\ell$ is $C^3$ on $\bR^n\setminus  \bd(D)$.

 \item[(${\bf P}_3$)] For every open neighborhood $U$ of $\bD$ there exists $\ve>0$ such that $f^{-1}([0,\ve])\subset U$.

 \item[(${\bf P}_4$)] There exists $\delta>0$ such that  any $t\in (0,\delta)$ is a  regular value of $\ell$.

 \item[(${\bf P}_5$)] If $(x_k)\in \bR^n\setminus \bd(D)$ is a sequence which converges to a point $x\in D_f(I)$,  and  if the line spanned by $\nabla\ell(x_k)$ converges to a line $L_\infty$, then the limit line $L_\infty$ is perpendicular to the tangent space $T_xD_f(I)$.

\end{enumerate}

\begin{figure}[h]
\centerline{\epsfig{figure=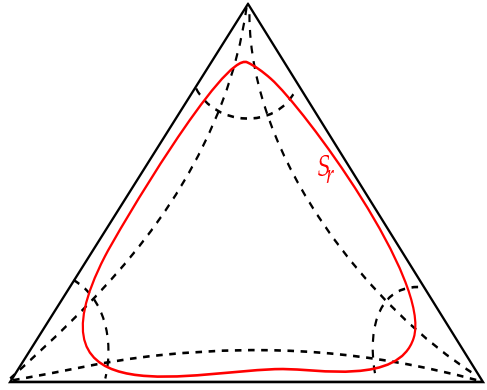,height=2.0in,width=2.5in}}
\caption{\sl   A tube system around the boundary of a tame simplex and its associated  profile.}
\label{fig: 18}
\end{figure}

We have depicted in Figure \ref{fig: 18} a tame $2$-simplex, with a tube system and the associated  boundary profile.

For every $r\in(0,1)$ and $I\subsetneq\{0,\dotsc, m\}$   we denote by $T_I(r)$ the closed tube
\[
T_I(r):=\bigl\{\, y\in T_I;\;\;\rho_I(y)\leq r\bar{\ve}_I({y})\,\bigr\}.
\]
We set
\[
\eT(r):=\bigcup_{\# I\leq m} T_I(r),\;\; \eT^0(r):=\bigcup_{\# I< m} T_I(r).
\]
Note that $\eT^0(r)$ is a neighborhood of the $(m-2)$-dimensional skeleton of $D_f$. Moreover, if $r\leq \frac{1}{3}$ then
\[
T_I(r)=\bigl\{ \ell_I\leq r\,\bigr\}.
\]
We can find a  tame function $\tau: (0,1)\ra (0,1)$, $r\mapsto \tau(r)$  satisfying the following conditions.
\begin{itemize}

\item $\tau(r)$ is a regular value of $\ell|_{D^f}$.

\item $D^f\cap \{\ell=\tau(r)\}\subset \eT(r)$.

\end{itemize}

We set $S_r:=D^f\cap \{\ell=\tau(r)\}$.  $S_r$ is a  compact, tame oriented $C^3$-submanifold  of $D_f$ of dimension $(m-1)$ which  approaches  $\bd(D^f)$ as $r\ra 0$ (see Figure \ref{fig: 18} for a $2$-dimensional rendition of $S_r$).  $S_r$   has a natural orientation  as boundary of
\[
D_r:=D\cap \{\ell\geq \tau(r)\,\}.
\]
We will prove  that
\begin{equation}
\forall \eta\in\Omega^{k-1}_{cpt}(\bR^n):\;\;\lim_{r\searrow 0}\int_{S_r}\eta=\sum_{k=0}^m(-1)^k\lan\eta, [D^k_f,\ori_k]\ran.
\label{eq: pa-curr}
\end{equation}
Clearly, (\ref{eq: pa-curr}) implies (\ref{eq: pa-symp}).

For every $x\in  D_f^k$,  we denote by $C_x$ the fiber of the  projection
\[
\pi_{I_k}: D\cap T_{I_k}\ra D_f^k.
\]
This fiber is a $C^3$-curve, and the map
\[
C_x\ra (0,1),\;\; y\mapsto s_x(y):= \frac{1}{\ve_i(x)}\rho(y)
\]
is a $C^3$-diffeomorphism. We  think of $s_x$ as a parameter along $C_x$ so  that the restriction of $\ell$ to $C_x$  can be regarded as a function of one variable $s=s_x$.

\begin{lemma} There exists a definable function  $\delta\mapsto r_1=r_1(\delta)$, such that, for all $x\in D_f^k\setminus \eT^0(\delta)$,  and all $r<r_0$,  the equation $\ell(y)=\tau(r)$ has at \emph{exactly one solution}   $y(r,\delta)\in C_x$. In other words, for any $x\in D_f^k$, and any   $r<r_1(\delta)$ the manifold $S_r$ intersects the fiber $C_x$  at a single point $y(r,\delta)$.
\label{lemma: profile}
\end{lemma}

\proof Let $x\in D_f^k\setminus \eT^0(\delta)$. Then, along $C_x$ we can use the parameter $s=s_x$, and we can think of the restriction of $\ell$ to $C_x$ as a $C^3$ function of a single variable $s$. Then
\[
\ell|_{C_x}=\ell_{I_k}|_{C_x}\prod_{I\neq I_k} \ell_{I_k}|_{C_x}.
\]
Observe that
\[
\ell_{I_k}|_{C_x}=s,\;\;\ell(s)=\ell|_{C_x}(s)= su(s),\;\; u(s):=\prod_{I\neq I_k} \ell_{I_k}|_{C_x}(s).
\]
There exists a constant $\nu=\nu(\delta)>1$, \emph{independent}  of $x\in D_f^k\setminus \eT^0(\delta)$, such that
\begin{subequations}
\begin{equation}
u(s)> \frac{1}{\nu},
\label{eq: lower}
\end{equation}
\begin{equation}
|u'(s)|< \nu,\;\;\forall s\in (0,1).
\label{eq: upper}
\end{equation}
\end{subequations}
To find one solution  of the equation $su(s)=\tau(r)$ we regard it as a fixed point problem
\[
s= g(s)= \frac{\tau(r)}{u(s)}.
\]
Observe that $u(0)>0$ so it suffices to  have $g(s)\leq 1$, i.e.,
\[
\tau(r)\leq u(s),\;\;\forall s\in[0,1].
\]
Using (\ref{eq: lower}) we deduce that  if
\[
\tau(r) \leq \frac{1}{\nu(\delta)},
\]
then we have at least one solution.   In  fact any, solution $\si$ of this equation must  satisfy the inequality
\[
0 <\si < \nu(\delta)\tau(r)
\]
To prove  the uniqueness,  it suffices  show that the derivative of $s\mapsto su(s)$ is positive in the interval
\[
0 <s < \min\bigl\{\,1,\nu(\delta)\tau(r)\,\bigr\}.
\]
We have
\[
\ell'(s)= u(s) + su'(s)\stackrel{(\ref{eq: upper})}{\geq} u(s) -\nu s \stackrel{(\ref{eq: lower})}{\geq}\frac{1}{\nu}-\nu^2\tau(r).
\]
Hence, if $\tau(r)<\frac{\delta}{\nu(\delta)^3}$, we also have uniqueness.   Moreover, the unique solution satisfies
\[
s<\min\{\delta, \nu(\delta)\tau(r)\},
\]
 i.e.,
\begin{equation}
y(r,\delta)\in T_{I_k}(\delta).
\label{eq: dist1}
\end{equation}
By definable selection, we can find a definable function $r_1(\delta)$ such that for all $r<r_1(\delta)$ we have $\tau(r)<\frac{\delta}{\nu^3}$.\qed

Set
\[
D_f^k(\delta)= D_f^k\setminus \eT^0(\delta),\;\;S^k_r(\delta)=\bigl\{\, y\in T_{I_k}\cap D;\;\;\pi_{i_k}(y)\in D_f^k(\delta)\,\bigr\}.
\]
\[
S_r(\delta)=S_r\setminus \bigcup_k S^k_r(\delta).
\]
Let $y\in S_r$.       For every  \emph{oriented, orthonormal} frame $\boldsymbol{f}=(f_1,\dotsc, f_{m-1})$ of $T_yS_r$ we get a scalar
\[
\eta(\boldsymbol{f},y)=\eta_y(f_1,\dotsc, f_{m-1}).
\]
This scalar  is \emph{independent} of the  frame $\boldsymbol{f}$, and thus defines a $C^3$-function $\eta_r$ on $S_r$.  Moreover, there exists $C_1>0$ such that
\[
|\eta_r(y)|\leq C_1,\;\;\forall r,\;\;\forall y\in S_r^k.
\]
Denote by $\eH^{m-1}$ the $(m-1)$-dimensional  Hausdorff measure.  We have
\[
\int_{S_r}\eta=\int_{S_r} \eta_k(y) d\eH^{m-1}(y).
\]
In particular
\[
\Bigl|\int_{S_r(\delta)}\eta\Bigr|\leq C_1\eH^{m-1}(S_r(\delta)).
\]
For $r<r_1(\delta)$ we have
\[
\int_{S_r}\eta= \int_{S_r(\delta)}\eta+\sum_k\int_{S_r^k(\delta)} \eta.
\]
Hence
\[
\Bigl| \int_{S_r}\eta- \sum_{k=0}^m(-1)^k\int_{D_f^k} \eta\,\Bigr|
\]
\[
\leq \underbrace{\Bigl|\int_{S_r(\delta)}\eta\Bigr|}_{T_1(r,\delta)}+\underbrace{\sum_{k=0}^m\Bigl|\int_{S_r^k(\delta)}\eta-(-1)^k\int_{D_f^k(\delta)} \eta\Bigr|}_{T_2(r,\delta)}+\underbrace{\sum_{k=0}^m\Bigl|\int_{D_f^k}\eta-\int_{D_f^k(\delta)}\eta\Bigr|}_{T_3(\delta)}
\]
We will prove the following things.

\begin{lemma} There exists  a  function $\ve\mapsto \delta_1(\ve)$ such that if  $\delta<\delta_1(\ve)$ and $r<r_1(\delta)$  then
\[
T_1(r,\delta)<\frac{\ve}{3}.
\]
\label{lemma: T1}
\end{lemma}

 \begin{lemma} There exists a  function $\ve\mapsto \delta_3(\ve)$ then
\[
T_3(\delta)<\frac{\ve}{3}.
\]
\label{lemma: T3}
\end{lemma}

\begin{lemma} There exists  a  function $(\ve,\delta)\mapsto r_2(\delta,\ve)$ such that if  $r< r_2(\delta,\ve)$  we have
\[
T_2(r,\delta) <\frac{\ve}{3}.
\]
\label{lemma: T2}
\end{lemma}

Assuming the above results, the   equality (\ref{eq: pa-curr}) is proved as follows. Fix $\ve>0$.  Choose  $\delta<\min\{\delta_1(\ve), \delta_2(\ve)\}$. Then, if $r<\min\{ r_1(\delta), r_2(\delta,\ve)\}$ we have
\[
T_1(r,\delta)+T_2(r,\delta)+T_3(\delta)<\ve.
\]

Using (\ref{eq: dist1}) we deduce that if $r<r_1(\delta)$ then
\[
S_r(\delta)\subset S_r\cap \eT^0( 2\delta).
\]
Lemma \ref{lemma: T1} and Lemma \ref{lemma: T3} are both consequences of the following result.

\begin{lemma} Suppose $X$ is a tame $C^3$-manifold of dimension $(m-1)$. Then
\[
\lim_{\delta\searrow 0} \eH^{m-1}(X\cap \eT^0(\hbar)\,\bigr)=0.
\]
\label{lemma: hausdorf}
\end{lemma}

\proof Denote by $\Graff^{m-1}$ the Grassmannian of affine planes in $\bR^n$ of codimension $(m-1)$. Denote by $\mu_{m-1}$ and invariant measure on $\Graff^{m-1}$, and set
\[
X_\hbar:=\cl\bigl( X \cap \eT^0(\hbar)\,\bigr).
\]
Then,  from Crofton formula (see \cite{BK, Feder}) we deduce
\[
\eH^{m-1}(X\cap \eT^0(\hbar)\,\bigr)= \eH^{m-1}(X_\hbar)=\int_{\Graff^{m-1}}\chi(L\cap X_\hbar) d\mu_{m-1}(L).
\]
The function
\[
(0,1)\times \Graff^{m-1}\ni (\hbar, L)\mapsto \chi(L\cap X_\hbar)
\]
is definable and thus  its range is finite.  From  dominated convergence theorem we deduce
\[
\lim_{\hbar\searrow 0}\eH^{m-1}(X_\hbar)=\int_{\Graff^{m-1}}\lim_{\hbar\searrow 0} \chi(L\cap X_\hbar).
\]
Suppose $L\in \Graff^{m-1}$ is such that
\[
\chi_0(L):= \lim_{\hbar\searrow 0} \chi(L\cap X_\hbar)\neq 0.
\]
Then the definable set $L\cap X_\hbar$ is nonempty for all $\hbar$ sufficiently small.  In particular, we can find a definable function
\[
\hbar \mapsto x_\hbar\in L\cap X_\hbar
\]
defined in a neighborhood of $0$.    Then the limit  $x_0=\lim_{\hbar\searrow 0} x_\hbar$ exists and it is a point in the intersection of $L$ with the  $(m-2)$-skeleton of $D$. We denote this skeleton by $ D^{(m-2)}$. Thus
\[
\chi_0(L)\neq 0\Longrightarrow L\cap D^{(m-2)} \neq\emptyset.
\]
The function
\[
\Graff^{m-1}\ni L\mapsto \chi_0(L)\in \bZ
\]
is definable and thus bounded. Hence
\[
\int_{\Graff^{m-1}}|\chi_0(L)|d\mu_{m-1}(L)\leq C\mu_{m-1}\bigl(\,\{ L\in\Graff^{m-1};\;\;L\cap D^{(m-2)}\neq\emptyset\} \,\bigr)
\]
By Sard's theorem, the  definable set
\[
\{ L\in\Graff^{m-1};\;\;L\cap D^{(m-2)}=\emptyset\}
\]
is \emph{dense}\footnote{A typical codimension $(m-1)$ affine plane will not intersect a  manifold of dimension $\leq (m-2)$, and $D^{(m-2)}$ is a finite union of such manifolds, $D_f(I)$, $\#I\leq m-1$.} in $\Graff^{m-1}$.

Hence, if $d=d(m,n)$ denotes the  dimension of $\Graff^{m-1}$, then
\[
\dim \{ L\in\Graff^{m-1};\;\;L\cap D^{(m-2)}\neq\emptyset\} <d.
\]
Up to a multiplicative constant $c>0$ we have
\[
\mu_{m-1}=c\eH^d,
\]
from which we deduce
\[
\mu_{m-1}\bigl(\,\{ L\in\Graff^{m-1};\;\;L\cap D^{(m-2)}\neq\emptyset\} \,\bigr)=\eH^d\bigl(\,\{ L\in\Graff^{m-1};\;\;L\cap D^{(m-2)}\neq\emptyset\} \,\bigr)=0.
\]
This completes the proof of Lemma \ref{lemma: hausdorf}. \qed

Lemma \ref{lemma: T3} is clearly a special case of Lemma \ref{lemma: hausdorf}.   To see that Lemma \ref{lemma: T1} is also a special case of Lemma \ref{lemma: hausdorf} observe that  for every $\hbar>0$ there exists $\delta_1(\hbar)$ such that if  $\delta<\delta_1(\hbar)$ and $r<r_1(\delta)$ then $S_r(\delta)\subset \eT^0(\hbar)$.

\noindent {\bf Proof of Lemma \ref{lemma: T2}.}      Fix $\delta<\frac{1}{2}$. We have to prove that for every $k=0,\dotsc, m\}$ we have
\[
\lim_{r\searrow 0}\int_{S_r^k(\delta)} \eta =(-1)^k\int_{D_f^k(\delta)}\eta.
\]
By Lemma \ref{lemma: profile},    for $r<r_1(\delta)$, the projection $\pi_{I_k}$ induces a homeomorphism
\[
S_r^k(\delta)\ra D_f^k(\delta).
\]
For simplicity we write $\bar{y}=\pi_{I_k}(y)$. We want to prove that
\begin{equation}
\lim_{r\ra 0} \sup_{y\in S_r^k(\delta)}\dist(\, T_yS_r^k(\delta), T_{\bar{y}}D_f^k(\delta)\,)=0.
\label{eq: diffeo1}
\end{equation}
We argue by contradiction. We can therefore find a constant $c>0$ and definable map $r\mapsto y_r\in S_r(\delta)$  such that
\begin{equation}
\dist(\, T_{y_r}S_r^k(\delta), T_{\bar{y}_r}D_f^k(\delta)\,)>c,\;\;\forall r<r_(\delta).
\label{eq: contr1}
\end{equation}
 Both limits $\lim_{r\ra 0}y_r$ and $\lim_{r\ra 0}\bar{y}r$  exist and they  coincide with a point $y_0\in\cl(D_f^k(\delta))\subset D^k_f$.   From the Whitney  regularity condition (a) and the property (${\bf P}_5)$ of the boundary profile $\ell$ we deduce
 \[
 \lim_{r\ra 0} T_{y_r} S_r^k(f)= T_{y_0} D_f^k.
 \]
 Clearly
 \[
\lim_{r\ra 0}T_{\bar{y}_r}D_f^k(\delta)= T_{y_0} D_f^k.
\]
This contradicts (\ref{eq: contr1})  and completes the proof of (\ref{eq: diffeo1}).

The equality (\ref{eq: diffeo1}) show that the map
\[
S_r^k(\delta)\ni y\mapsto \bar{y}\in  D_f^k(\delta)
\]
is a $C^3$-diffeomorphism for $r$ sufficiently small, and  changes the orientation by a factor of $(-1)^k$.   For every $y\in D^k_f(\delta)$ we denote by $\eta(y)$ the pairing between $\eta$ and oriented orthonormal frame of $T_yD_f^k$.

Using the change in variables formula we can write
\[
\int_{S_r^k(\delta)} \eta =(-1)^k\int_{D_f^k(\delta)}J_r(\bar{y})\eta_r(\bar{y})d\eH^{m-1}(\bar{y}),
\]
where $\eta_r(\bar{y})$ is the pullback of the function $\eta_r|_{S_r^k(\delta)}$ to $D_f^k(\delta)$, and  $J_r(\bar{y})$ is the Jacobian of the change in variables.  The  equality (\ref{eq: diffeo1}) and the continuity of the form $\eta$ imply that
\[
\lim_{r\ra 0} J_r(\bar{y})=1\;\;{\rm and}\;\; \lim_{r\ra 0}\eta_r(\bar{y})=\eta(\bar{y})
\]
  uniformly  on $D_f^k(\delta)$.  This completes the proof of  Lemma \ref{lemma: T2}   and of Theorem \ref{th: curr-bound}.\qed

%% file: tameflow-b.tex
We would like  to include  a few topological facts concerning tame sets. These are not needed in the main body of the paper, yet they may shed some light on  the subtleties of tame topology.

As we  mentioned in Section \ref{s: 1}, any compact tame set $S$ can be triangulated, i.e., there exists an affine  finite simplicial complex, and a tame homeomorphism  $\vfi:K\ra  S$.

Clearly, if $\vfi_i: K_i\ra S$, $i=0,1$, are two tame triangulations, then   the map
\[
\vfi_1\circ\vfi_0^{-1}: K_0\ra K_1
\]
is a tame  homeomorphism.      It turns out that the existence of a tame homeomorphism between two compact PL spaces imposes a severe restriction on these spaces. More precisely, M. Shiota has proved (see \cite[Chap. IV]{Shi}) the \emph{tame Haupvermutung}, namely  that   two compact PL spaces are  PL-homeomorphic if and only if  they are  tamely homeomorphic. Given this result, we can define the link of a point in a compact tame space to be its PL link as defined  e.g. in \cite[Chap. 2]{RS}.

To appreciate the strength of this result, consider the classical example of Cannon-Edwards \cite{Can, Lat}, the double suspension of a non-simply connected homology $3$-sphere, say the Poincar\'{e} sphere $\Sigma(2,3,5)$. This is a simplicial complex $K$  which is  homeomorphic, \emph{but not $PL$-equivalent} to the $5$-sphere, equipped  with the triangulation as boundary of a $6$-simplex.   The Shiota \emph{Hauptvermutung} implies that $X$ and $S^3$ are not tamely homeomorphic. In particular, there cannot exist a semi-algebraic homeomorphism from the round $5$-sphere to $X$.

In the above paragraphs, we have defined  the link of a point in a compact tame space indirectly, via triangulations and the tame Hauptvermutung.    We can  attempt a more intrinsic  approach, namely given a point $p_0$ in a compact tame set $X$, and a  tame continuous function $w:X\ra [0,\infty)$,  such that $w^{-1}(0)=\{p_0\}$, we can  define the link of $p_0$ as the level set $\{w=\ve\}$, where $\ve>0$ is sufficiently small. The homeomorphism type of this set is independent of $\ve>0$,  but  at this point, we do not know how to eliminate the dependence on $\ve>0$.

To understand   the subtleties of this question consider a closely related problem.

\smallskip

\emph{Suppose $w:\bR^n\ra [0,\infty)$ is a tame continuous function such that $w^{-1}(0)=\{0\}$. Then there exists $r_0>0$ such that  for every $x\in  \bR^n$, $|x|=1$, the function
\[
[0,r_0]\ni t\mapsto w(tx)
\]
is strictly decreasing, i.e.,  in a neighborhood of $0$ the function $w$ is a Lyapunov function for the radial flow.}

\smallskip

If $w$ is a $C^1$-function, then      this fact  follows from the non-depravedness arguments in  \cite[Sec. 2.4]{GM}.  When $w$ is merely continuous (and tame), this seems to be a rather slippery problem.

Let us observe that if $\Phi$  is a gradient-like  tame flow on the compact tame space, $x_0$ is a n isolated stationary point of $\Phi$ and $u, v: X\ra \bR$ are two Lyapunov functions such that $u(0)=v(0)=0$,, then  unstable links
\[
\eL^-_u(x_0)=W^-(x_0)\cap \{u=-\ve\},\;\;\eL_v^-(x_0)=W^-(x_0)\cap \{v=-\ve\},
\]
are \emph{tamely} homeomorphic for $\ve$ sufficiently small. In other words, the tame homeomorphism type of the unstable link, is a \emph{dynamical} invariant of the  stationary point.

%% file: tameflow.bbl
\begin{thebibliography}{XXXXXX}




\bibitem{Ar2} V.I. Arnold:  {\sl Geometrical Methods in the Theory of Ordinary Differential Equations}, 2nd Edition,  Springer Verlag 1988.






\bibitem{Bjorn} A. Bjorner: {\sl Topological methods},  Hadbook of Combinatorics, R. Graham, M. Gr\"{o}tschel and L. Lov\'asz, (Eds), North-Holland, Amsterdam, 1995, p. 1819-1872.


\bibitem{BCR} J. Bochnak, M. Coste, M.-F. Roy: {\sl Real Algebraic Geometry}, Springer Verlag, 1998.


\bibitem{BK} L. Br\"{o}cker, M. Kuppe: {\sl Integral Geometry of tame sets},  Geom. Ded.  {\bf 82}(2000), 285-323.

\bibitem{Can} J.W. Cannon: {\sl  Shrinking cell-like decompositions of manifolds. Codimension three} Ann. of Math., {\bf 110}(1979), 83-112.

\bibitem{Con}  C. Conley: {\sl Isolated Invariant Sets and the Morse Index},  CBMS Regional Conf. Series, vol. 38, Amer. Math. Soc., 1978.


\bibitem{Co} M. Coste: {\sl An Introduction to $o$-minimal geometry}, Real Algebraic and Analytic Geometry  Network,

\verb+ http://www.ihp-raag.org/publications.php+

\bibitem{CH} R. Courant, D. Hilbert: {\sl  Methods of Mathematical Physics},  vol.II,  John Wiley Classics, 1989.


\bibitem{Dr1} L.van der Dries: {\sl  Tame Topology and $o$-minimal Structures}, London Math. Soc. Series, vol. 248, Cambridge University Press, 1998.


\bibitem{DrMM} L.van der Dries, A. Macintyre,  D. Marker: {\sl  The elementary theory of restricted analytic fields with exponentiation}, Ann. Math. {\bf 140}(1994), 183-205.


\bibitem{DrMi2} L.van der Dries, C. Miller:  {\sl Geometric categories and $0$-minimal structures} Duke Math. J. {\bf 84}(1996), 497-540.

\bibitem{Feder} H. Federer: {\sl Geometric Measure Theory},  Springer Verlag, 1969.

\bibitem{For} R. Forman: {\sl Morse theory for cell complexes}, Adv. in Math., {\bf 134}(1998), 90-145.



\bibitem{GWPL} C.G. Gibson, K. Wirthm\"{u}ller, A.A. du Plesis, E.J.N. Looijenga: {\sl  Topological Stability of Smooth Mappings},  Lect. Notes in Math., vol. 552, Springer Verlag, 1976.

\bibitem{GM} M. Goresky, R. MacPherson: {\sl Stratified Morse Theory},  Ergebnisse der Mathematik, vol. 14,  Springer Verlag, 1988.

\bibitem{Gra} H. Grauert: {\sl On Levi's problem and the embedding of real analytic manifolds}, Ann. of Math., {\bf 68}(1958), 460-472.

\bibitem{GJ} R.E. Greene, H. Jacobowitz: {\sl Analytic isometric embeddings}, Ann. of Math.,  {\bf 93}(1971), 189-204.


\bibitem{Hardt}  R. Hardt: {\sl Slicing and intersection theory for chains associated with real analytic varieties}, Acta Math., {\bf 129}(1972), 57-136.

\bibitem{Hardt1} R. Hardt: {\sl  Topological properties of subanalytic sets}. Trans. Amer. Math. Soc., {\bf 211}(1975), 57-70.



\bibitem{HL} F.R. Harvey, H.B. Lawson: {\sl Morse theory and Stokes' theorem},  Surveys in differential geometry,  259--311, Surv. Differ. Geom., VII, Int. Press, Somerville, MA, 2000.

\verb+ http://www.math.sunysb.edu/~blaine/+

\bibitem{Hir} M. Hirsch: {\sl Differential  Topology}, Springer-Verlag, 1976.



\bibitem{KaSch}  M. Kashiwara, P. Schapira: {\sl Sheaves on Manifolds},   Gr\"{u}ndlehren der mathematischen Wissenschaften, vol. 292, Springer Verlag, 1990.

\bibitem{Kato} T. Kato: {\sl Perturbation Theory for Linear Operators}, Classics in Mathematics, Springer Verlag, 1995.

\bibitem{Kho} A. Khovanskii: {\sl Fewnomials},  Transl. Math. Monogr., {\bf 88}, Amer.Math. Soc., 1991.

\bibitem{Lauden} F. Laudenbach: {\sl On the Thom--Smale complex}, Ast\'erisque, vol. 205(1992), Soc. Math. France.

\bibitem{Lat} F. Latour: {\sl Double suspension d'une  sph\`{e}re d'homologie [D'apr\`{e}s R. Edwards]}, S\'{e}minaire Bourbaki, Exp. 515, Lecture Notes in Math., vol. 710, Springer Verlag, 1979.

\bibitem{LS} J.-M. Lion, P. Speissegger: {\sl Analytic stratification in the pfaffian  closure of an $o$-minimal structure}, Duke Math. J., {\bf 103}(2000), 215-231.

\bibitem{Loi}  T.L. Loi: {\sl Verdier and Thom strict stratifications in $o$-minimal structure}, Ill. J. Math., {\bf 42}(1998), 345-356. {\sf math.DG/9704232}.


\bibitem{LW}  A.T. Lundell, S. Weingram,  {\sl The Topology of CW Complexes}, Van Nonstrand, 1969.



\bibitem{Morg} F. Morgan: {\sl Geometric Measure Theory. A  Beginner's Guide}, Academic Press, 1995.


\bibitem{N2} L.I. Nicolaescu: {\sl An Invitation to Morse Theory}, Universitext, Springer Verlag, 2007.

\bibitem{Pliss} V.A. Pliss: {\sl On the reduction of an analytic system of differential equations},  Diff. Uravneniya, {\bf 1}(1965), 153-161.

\bibitem{RS}  C.P. Rourke, B.J. Sanderson: {\sl  Introduction to Piecewise-Linear Topology}, Springer Verlag, 1982.

\bibitem{Roy} H.L. Royden: {\sl  The analytic approximation of differentiable mappings}, Math. Ann., {\bf 139}(1960), 171-179.

\bibitem{Shi}  M.  Shiota: {\sl Geometry of Subanalytic and Semialgebraic Sets}, Birkh\"{a}user, 1997.

\bibitem{Sal} D. Salamon: {\sl Connected  simple systems and the Conley index of isolated invariant sets}, Trans. A.M.S., {\bf 291}(1985), 1-41.

\bibitem{Sp1}  P. Speissegger: {\sl The Pfaffian closure of an $o$-minimal  structure}, J. Reine. Angew. Math., {\bf 508}(1999), 189-211.


\bibitem{Trot1} D.  Trotman: {\sl  Geometric versions of Whitney regularity conditions}, Ann. Scien. \'{E}c. Norm. Sup., {\bf 4}(1979), 453-463.


\bibitem{Trot2} D. Trotman: {\sl Une version microlocal de la condition $(w)$ de Verdier}, Ann. Inst. Fourier, {\bf 39}(1989), 825-829.

\bibitem{Verd}  J.L. Verdier: {\sl  Stratification de Whitney et th\'{e}or\`{e}me  de Bertini-Sard},  Invent. Math., {\bf 36}(1976), 295-312.


\bibitem{White} J.H.C. Whitehead: {\sl  Simplicial spaces,   nuclei and $m$-groups},  Proc. London Math. Soc., {\bf 45}(1939), 243-327.


\bibitem{Whit2} H. Whitney:  {\sl Geometric Integration Theory},  Princeton University Press, 1957.


\end{thebibliography}
